\documentclass{amsart}
\usepackage{enumitem, subcaption, color, mathtools, url, mathrsfs, amsmath,amssymb,amsthm,color,mathrsfs, amsrefs, nicefrac, dutchcal, array, url}
\usepackage[utf8]{inputenc}
\usepackage{physics}
\usepackage[T1]{fontenc}
\usepackage{tcolorbox}
\usepackage{bbm}
\usepackage{esint}

\newtheorem{theorem}{Theorem}[section]
\newtheorem{lemma}[theorem]{Lemma}

\theoremstyle{definition}
\newtheorem{definition}[theorem]{Definition}

\theoremstyle{remark}
\newtheorem{remark}[theorem]{Remark}
\newtheorem{Notation}[theorem]{Notation} 
\newtheorem{Convention}[theorem]{Convention} 

\numberwithin{equation}{section}

\newcommand{\nc}{\newcommand}
\nc{\RR}{\mathbb{R}}
\nc{\R}{\mathbb{R}}
\nc{\CC}{\mathbb{C}}
\nc{\C}{\mathbb{C}}
\nc{\mrm}{\mathrm}
\nc{\mL}{\mrm{L}}
\nc{\mC}{\mrm{C}}
\nc{\mA}{\mrm{A}}
\nc{\mH}{\mrm{H}}
\nc{\mW}{\mrm{W}}
\nc{\mF}{\mrm{F}}
\nc{\mM}{\mrm{M}}
\nc{\mK}{\mrm{K}}
\nc{\mD}{\mrm{D}}
\nc{\mX}{\mrm{X}}
\nc{\mY}{\mrm{Y}}
\nc{\mS}{\mrm{S}}
\nc{\Ec}{\mrm{E_c}}
\nc{\calL}{L}
\nc{\loc}{\mrm{loc}}
\nc{\comp}{c}
\nc{\supp}{\mrm{supp}}
\nc{\Hardy}{\mathfrak{H}}
\nc{\calH}{H}
\nc{\ctru}{\mathfrak{u}}
\nc{\ctrv}{\mathfrak{v}}
\nc{\be}{\boldsymbol{e}}
\nc{\bx}{\boldsymbol{x}}
\nc{\by}{\boldsymbol{y}}
\nc{\bz}{\boldsymbol{z}}
\nc{\bc}{\boldsymbol{c}}
\nc{\lbr}{\lbrack}
\nc{\rbr}{\rbrack}
\nc{\dsp}{\displaystyle}
\nc{\vphi}{\varphi}

\begin{document}
	
	\title[Modified-operator method for band diagrams]{Modified-operator method for the calculation of \\ band diagrams of crystalline materials}

	\author{Eric Canc\`es}\address{CERMICS, \'Ecole des Ponts and Inria Paris, 6 and 8 avenue Blaise Pascal, 77455 Marne-la-Vall\'ee, France.} \email{eric.cances@enpc.fr} 
	\author{Muhammad Hassan}\address{Sorbonne Universit\'e, CNRS, Universit\'e Paris Cit\'e, Laboratoire Jacques-Louis Lions (LJLL), F-75005 Paris, France.} \email{hassan@ljll.math.upmc.fr}
	\author{Laurent Vidal}\address{CERMICS, \'Ecole des Ponts and Inria Paris, 6 and 8 avenue Blaise Pascal, 77455 Marne-la-Vall\'ee, France.} \email{laurent.vidal@enpc.fr}

	\begin{abstract} 
		In solid state physics, electronic properties of crystalline materials are often inferred from the spectrum of periodic Schr\"odinger operators. As a consequence of Bloch's theorem, the numerical computation of electronic quantities of interest involves computing derivatives or integrals over the Brillouin zone of so-called energy bands, which are piecewise smooth, Lipschitz continuous periodic functions obtained by solving a parametrized elliptic eigenvalue problem on a Hilbert space of periodic functions. Classical discretization strategies for resolving these eigenvalue problems produce approximate energy bands that are either non-periodic or discontinuous, both of which cause difficulty when computing numerical derivatives or employing numerical quadrature. In this article, we study an alternative discretization strategy based on an ad hoc operator modification approach. While specific instances of this approach have been proposed in the physics literature, we introduce here a systematic formulation of this operator modification approach. We derive a priori error estimates for the resulting energy bands and we show that these bands are periodic and can be made arbitrarily smooth (away from band crossings) by adjusting suitable parameters in the operator modification approach. Numerical experiments involving a toy model in 1D, graphene in 2D, and silicon in 3D validate our theoretical results and showcase the efficiency of the operator modification approach.
	\end{abstract}
	
	
	
	\maketitle
	
	\section{Introduction}
	
	In solid state physics, macroscopic properties such as the electrical and thermal conductivities, heat capacity, magnetic susceptibility, and optical absorption of crystalline materials are often explained through the use of an independent electron model (see, e.g., \cite[Part II, Chapter 5]{harrison1980solid}, \cite[Chapter 12]{martin2020electronic}, \cite[Chapter 7]{kittel1996introduction}, and \cite[Chapter 5]{grosso2013solid}). This model consists of treating the crystalline material as an infinite, perfect crystal and modeling the electrons as independent of each other (quasiparticle approach) and evolving under the influence of an \emph{effective} periodic potential. The behavior of each electron is thus determined by the spectrum of an unbounded, self-adjoint, periodic Schr\"odinger operator acting on $L^2(\R^3)$ (see, e.g., \cite[Chapter XIII]{MR0493421}). Although the independent electron assumption might seem naive, this model has achieved great success in explaining basic phenomena such as the difference between conductors, semi-conductors and insulators, as well as describing the electronic properties of many ubiquitous non-strongly correlated materials (see, e.g., \cite[Part III]{harrison1980solid}, \cite[Part V]{martin2020electronic}, \cite[Chapters 10-12]{grosso2013solid}, and \cite[Chapters 6]{kittel1996introduction}). In addition, Kohn-Sham Density Functional Theory (DFT) provides a method to parameterize this independent-electron model and obtain {\em quantitatively} accurate results for a very large class of materials of practical interest (see, for instance, \cite{kaxiras_2003, dreizler2012density}).
	
	In the independent electron model, the practical computation of electronic quantities of interest is based on the use of the Bloch-Floquet transform (see, e.g., \cite[Chapter XIII]{MR0493421}). The Bloch transform essentially yields an explicit block-diagonal decomposition of the underlying Schr\"odinger operator into so-called Bloch \emph{fibers}, which are bounded-below, self-adjoint operators  acting on a space of periodic square-integrable functions. Thus, the problem of computing the spectrum of the periodic Schr\"odinger operator is reduced to one of calculating the low-lying eigenvalues of the Bloch fibers. These Bloch fibers are typically indexed by a parameter $\bold{k}$ that belongs to a $d$-dimensional torus (the Brillouin zone), and therefore each resulting eigenvalue (often referred to as an energy) can be viewed as a periodic function on the $d$-dimensional Brillouin Zone. It is thus common in the solid-state physics literature to speak of \emph{energy bands}.
	
	Energy bands provide both qualitative and quantitative information about the electronic properties of the crystalline material being studied (see, e.g., the references quoted above). Insulators and conductors for instance, are characterized by the presence or absence, respectively, of an energy \emph{band gap}. Other electronic quantities of interest can be expressed in terms of integrals (over the Brillouin zone) or derivatives involving the energy bands (see, e.g., \cite{MR4071823}). In order to estimate important quantities such as the \emph{integrated density of states} or the \emph{integrated density of energy} (see Section \ref{sec:2} for precise definitions of these quantities), it is therefore necessary to
	\begin{itemize}
		\item sample the energy bands at different $\bold{k}$-points which corresponds to solving approximately the $\bold{k}$-fiber eigenvalue problems posed on a periodic domain;
		\item use suitable numerical quadrature to approximate integrals involving these energy bands.
	\end{itemize}
	
	Concerning the first step, the famous Monkhorst-Pack numerical scheme \cite{MR418801} is widely used to select the specific $\bold{k}$-points at which the eigenvalue problem is to be solved. For the second step, a number of numerical quadrature methods for integration in the Brillouin zone have been proposed including the well-known linear tetrahedron method (see, e.g., \cite{lehmann1972numerical}) and the improvement due to Bl\"ochl et al. \cite{blochl1994improved}, and smearing methods (see, e.g., \cites{morgan2018efficiency, pickard1999extrapolative, henk2001integration, methfessel1989high}). 
	
	From a mathematical and computational point of view, two natural questions now arise. First, which discretization method should be employed in the actual numerical resolution of the $\bold{k}$-fiber eigenvalue problems, and second, what can be said about the convergence rate of the various numerical quadrature methods that are in use? For technical reasons, these questions become particularly relevant for \emph{metallic} systems (see, e.g., \cite{MR3554547} for an analysis involving insulators and semi-conductors), and in this case, the latter question has recently been addressed by the first author and coworkers in \cite{MR4071823}. The analysis carried out in \cite{MR4071823} revealed that the \emph{periodicity} (with respect to the Brillouin zone) and \emph{regularity} properties of the energy bands play a crucial role in the quadrature convergence rates, which of course is consistent with the experience of classical integration schemes in numerical analysis. Given that different eigenvalue discretization methods can conceivably produce (and in fact \emph{do} produce, as we show in Section \ref{sec:3}) energy bands that possess different regularity properties or may be altogether aperiodic, the choice of discretization scheme becomes vitally important. This article is concerned precisely with the study of approximation strategies for energy bands in the Brillouin zone. 
	
	The remainder of this article is organized as follows: We begin in Section \ref{sec:2} by introducing our notation and stating precisely the problem setting and governing equations. We then present in Section \ref{sec:3} two classical Galerkin discretization strategies for approximating the $\bold{k}$-fiber eigenvalue problems, and we show the problems associated with the energy bands produced by these classical approaches. Next, in Section \ref{sec:4}, we present an alternative discretization scheme, systematizing ideas first introduced in the physics literature (see Remark \ref{rem:literature} below), which is based on modifying in a controlled manner the underlying $\bold{k}$-fiber operator. In Section \ref{sec:5}, we present our two main results on the analysis of this alternative approach: we derive a priori error estimates with respect to a discretization cutoff for the modified energy bands, and we show that these bands are periodic and can be made arbitrarily smooth (away from band crossings) by adjusting suitable parameters in the operator modification approach. Numerical experiments in Section \ref{sec:6} involving a 1D toy model, and two real materials (graphene and face-centered cubic silicon), validate our theoretical results and showcase the efficiency of the operator modification approach. Finally, in Section \ref{sec:7}, we present the proofs of our main results.

	\section{Problem Formulation and Setting} \label{sec:2}~

	Perfect crystals are structures composed of a periodic arrangement of atoms. Such structures can therefore be described very conveniently through the use of a suitable lattice: assuming a $d$-dimensional lattice with~$d\in \mathbb{N}^*=\left\{1,2, 3, \ldots \right\}$, we denote by $\{\bold{a}_i\}_{i=1}^d$ a collection of~$d$ linearly independent primitive vectors in $\mathbb{R}^d$, and we denote by $\{\bold{b}_i\}_{i=1}^d \subset \mathbb{R}^d$ the corresponding reciprocal vectors, i.e., vectors in $\mathbb{R}^d$ that satisfy $\bold{a}_i \cdot \bold{b}_j = 2\pi\delta_{ij} ~ \forall i, j \in \{1, \ldots , d\}$. The primitive lattice $\mathbb{L} \subset \mathbb{R}^d$ and reciprocal lattice $\mathbb{L}^* \subset \mathbb{R}^d$ are then defined as 
	\begin{align*}
		\mathbb{L} := \left\{ \mathbb{Z}\bold{a}_1+ \ldots + \mathbb{Z} \bold{a}_d\right\} \qquad \text{and} \qquad \mathbb{L}^* := \left\{ \mathbb{Z}\bold{b}_1+ \ldots + \mathbb{Z} \bold{b}_d \right\}.
	\end{align*}

	We denote by $\Omega\subset \mathbb{R}^d$ and $\Omega^* \subset \mathbb{R}^d$ the first Wigner-Seitz unit cell of the primitive and reciprocal lattice respectively. Recall that the first Wigner-Seitz unit cell of a lattice in $\mathbb{R}^d$ is the locus of points in $\mathbb{R}^d$ that are closer to the origin of the lattice than to any other lattice point. The first Wigner-Seitz cell $\Omega^*$ of the reciprocal lattice is called the (first) \emph{Brillouin zone}.
	
	Finally for clarity of the subsequent exposition, let us introduce the so-called translation operator and the related notion of lattice periodicity: Given any~$\bold{y} \in \mathbb{R}^d$ and denoting $\mathcal{D}(\R^d):= \mathscr{C}^{\infty}_{\rm c}(\R^d)$ the space of complex-valued smooth compactly-supported functions on $\R^d$, we define the translation operator $\tau_{\bold{y}} \colon \mathcal{D}(\R^d) \rightarrow \mathcal{D}(\R^d)$ as the mapping with the property that 
	\begin{align*}
		\forall \Phi \in \mathcal{D}(\R^d)\colon  \qquad \big(	\tau_{\bold{y}} \Phi \big)(\bold{x}) := \Phi(\bold{x} - \bold{y})  \qquad \text{for a.e. } \bold{x} \in \mathbb{R}^d.
	\end{align*}
	It follows that for any $\bold{y} \in \mathbb{R}^d$, the translation operator extends by duality as a mapping $\tau_{\bold{y}} \colon \mathcal{D}'(\R^d) \rightarrow \mathcal{D}'(\R^d)$.
	
	Given now some $\Phi \in \mathcal{D}'(\R^d)$, we will say that $\Phi$ is $\mathbb{L}^*$-periodic (resp. $\mathbb{L}$-periodic) if $\tau_{\bold{G}}\Phi = \Phi$ for all~$\bold{G} \in \mathbb{L}^*$ (resp. $\tau_{\bold{R}}\Phi = \Phi$ for all~$\bold{R} \in \mathbb{L}$).

	\subsection{Function spaces and norms}~
	
	We define the function space $L^2_{\rm per}(\Omega)$ as the set of (equivalence classes of) functions given by
	\begin{align*}
		L^2_{\rm per}(\Omega):= \left\{ f \in  L^2_{\rm loc}(\mathbb{R}^d) \text{ such that } f \text{ is } \mathbb{L}\text{-periodic}\right\},
	\end{align*}
	equipped with the inner-product
	\begin{align*}
		\forall f, g \in L^2_{\rm per}(\Omega) \colon \quad (f, g)_{L^2_{\rm per}(\Omega)}:=  \int_{\Omega} \overline{f(\bold{x})} g(\bold{x})\, d\bold{x},
	\end{align*}
	where $L^2_{\rm loc}(\mathbb{R}^d)$ denotes the space of complex-valued, locally square-integrable functions on $\mathbb{R}^d$, and $\overline{f(\cdot)}$ indicates the complex conjugate of $f(\cdot)$. The spaces $L^p_{\rm per}(\Omega)$, $p \in [1,2) \cup (2 ,\infty]$ are defined analogously.
	
	We denote by $\mathcal{B}$, the orthonormal Fourier basis of $L^2_{\rm per}(\Omega)$, i.e., 
	\begin{align*}
		\mathcal{B} := \left\{e_{\bold{G}}(\bold{x}):=\frac{1}{\vert \Omega \vert^{\frac{1}{2}} }e^{\imath \bold{G} \cdot \bold{x}}\colon \quad  \bold{G}\in \mathbb{L}^* \right\}.
	\end{align*}
	Given any $f \in L^2_{\rm per}(\Omega)$, we will frequently express $f$ in the form
	\begin{align*}
		f= \sum_{\bold{G} \in \mathbb{L}^*} \widehat{f}_\bold{G} e_{\bold{G}}, \qquad \text{where }~ \widehat{f}_{\bold{G}} &:= \int_{\Omega} {f(\bold{x})}  \overline{e_{\bold{G}}(\bold{x})}\, d \bold{x}, \qquad \text{and we have }\\ 
		  \sum_{\bold{G} \in \mathbb{L}^*} \big\vert\widehat{f}_\bold{G}\big\vert^2&= \int_{\Omega} \vert f(\bold{x}) \vert^2 \; d\bold{x} < \infty.
	\end{align*}
	
	Periodic Sobolev spaces of positive orders are constructed analogously. Indeed, we define for each~$s > 0$ the set
	\begin{align*}
		H^s_{\rm per}(\Omega):= \left\{ f \in L^2_{\rm per}(\Omega) \colon  \sum_{\bold{G} \in \mathbb{L}^*} \left(1+\vert \bold{G}\vert^2\right)^s \big\vert\widehat{f}_\bold{G}\big\vert^2 < \infty\right\},
	\end{align*}
	equipped with the inner-product
	\begin{align*}
		\forall f, g \in H^s_{\rm per}(\Omega) \colon \quad (f, g)_{H^s_{\rm per}(\Omega)}:=  \sum_{\bold{G} \in \mathbb{L}^*} \left(1+\vert \bold{G}\vert^2\right)^s \overline{\widehat{f}_\bold{G}}\widehat{g}_\bold{G}.
	\end{align*}
	Naturally, we have $H^0_{\rm per}(\Omega):= L^2_{\rm per}(\Omega)$, and we define periodic Sobolev spaces of negative orders through duality, i.e., for each $s>0$ we define $H^{-s}_{\rm per}(\Omega) := \left(H^{s}_{\rm per}(\Omega)\right)'$, and we equip $H^{-s}_{\rm per}(\Omega) $ with the canonical dual norm. 
	
	Finally, given a Banach space $X$, we will write $\mathcal{L}(X)$ to denote the Banach space of bounded linear operators from $X$ to $X$, equipped with the usual operator norm.

	\subsection{Governing operators and quantities of interest}\label{sec:2b}~
	
	In this section, we assume that the electronic properties of the crystal that we study are encoded in an effective one-body Schr\"odinger operator  
	\begin{equation}\label{eq:Hamiltonian}
		\mH:= -\frac{1}{2} \Delta + V \qquad \text{acting on } L^2(\mathbb{R}^d) \quad \text{with domain } H^2(\mathbb{R}^d),
	\end{equation}
	where $V \in L^{\infty}_{\rm per}(\Omega)$ is an $\mathbb{L}$-periodic \emph{effective potential}. Many electronic properties of the crystal we study can be computed from the spectral decomposition of this one-body Hamiltonian operator $\mH$, and we are therefore interested in its analysis and computation. The classical approach to this problem relies on the use of the Bloch-Floquet transform (see, e.g.,~\cite[Chapter XIII]{MR0493421}), which we will now briefly present. The following exposition is based on the article \cite{MR4289497}.

	We begin by introducing for each $\bold{G} \in \mathbb{L}^*$, the unitary multiplication operator $T_{\bold{G}} \colon L^2_{\rm per}(\Omega) \rightarrow L^2_{\rm per}(\Omega)$ defined as
	\begin{equation*}
		\forall v \in L^2_{\rm per}(\Omega) \colon  \qquad \big(T_{\bold{G}}v \big)(\bold{x}) = e^{-\imath \bold{G} \cdot \bold{x}} v(\bold{x}) \hspace{1cm} \text{for a.e.}~ \bold{x} \in \mathbb{R}^d.
	\end{equation*}
	
	Next, we introduce the Hilbert space of $\mathbb{L}^*$-quasi-periodic, $L^2_{\rm per}(\Omega)$-valued functions on $\mathbb{R}^d$ as the vector space
	\begin{align*}
		L^2_{\rm qp} (\mathbb{R}^d; L^2_{\rm per}(\Omega) ):= \Big\{\mathbb{R}^d \ni \bold{k} \mapsto &u_{\bold{k}} \in L^2_{\rm per}(\Omega) \colon \int_{\Omega^*} \Vert u_{\bold{k}}\Vert^2_{L^2_{\rm per}(\Omega)}{\rm d}{\mathbf{k}} < \infty \quad \text{and}\\ &u_{\bold{k}+ \bold{G}}= T_{\bold{G}}u_{\bold{k}} ~ \forall ~\bold{G} \in \mathbb{L}^* ~\text{ and a.e. } \bold{k}\in \mathbb{R}^d\Big\},
	\end{align*}
	equipped with the inner product
	\begin{align*}
		\forall u,v  \in 	L^2_{\rm qp} (\mathbb{R}^d; L^2_{\rm per}(\Omega) ) \colon \qquad 	\left(u, v\right)_{L^2_{\rm qp} (\mathbb{R}^d; L^2_{\rm per}(\Omega) )} = \fint_{\Omega^*} \left(u_{\bold{k}}, v_{\bold{k}}\right)_{L^2_{\rm per}(\Omega) }   d\bold{k},
	\end{align*}
	where we have denoted $\fint_{\Omega^*} := \frac{1}{\vert \Omega^* \vert} \int_{\Omega^*} $ and we have used the subscript `${\rm qp}$' to highlight \emph{quasi-periodicity}.
	
	\sloppy The Bloch-Floquet transform is now the unitary mapping from $L^2(\mathbb{R}^d)$ to $L^2_{\rm qp} (\mathbb{R}^d; L^2_{\rm per}(\Omega) )$ with the property that any $u \in \mathscr{D}(\mathbb{R}^d)$ is mapped to the element of $L^2_{\rm qp} (\mathbb{R}^d; L^2_{\rm per}(\Omega) )$ defined as
	\begin{align*}
		\R^d \ni \bold{k} \mapsto  u_{\bold{k}}:=  \sum_{\bold{R} \in \mathbb{L}} u\big(\bullet + \bold{R}\big) e^{-\imath \bold{k} \cdot (\bullet + \bold{R})} \in L^2_{\rm per}(\Omega).
	\end{align*}
	
	Any bounded linear operator $\mA \colon L^2(\mathbb{R}^d) \rightarrow L^2(\mathbb{R}^d)$ that is $\mathbb{L}$-periodic, i.e., one that commutes with the translation operator $\tau_{\bold{R}}$ for all $\bold{R} \in \mathbb{L}$, is decomposed by the Bloch transform in the following sense: there exists a function $\bold{k} \mapsto \mA_{\bold{k}} $ in $L^{\infty}_{\rm qp}\left(\mathbb{R}^d; \mathcal{L}\big( L^2_{\rm per}(\Omega)\big)\right)$ such that for any $u \in L^2(\mathbb{R}^d)$, all $\bold{G} \in \mathbb{L}^* $ and a.e. $\bold{k} \in \mathbb{R}^d$ it holds that
	\begin{align}\label{eq:bloch_decomp}
		(\mA u)_{\bold{k}} = \mA_{\bold{k}}u_{\bold{k}}, \qquad \text{and} \qquad \mA_{\bold{k}+\bold{G}}= T_{\bold{G}} \mA_{\bold{k}}T_{\bold{G}}^*. 
	\end{align}
	where the operators $\left(\mA_{\bold{k}}\right)_{\bold{k} \in \mathbb{R}^d} \in \mathcal{L}\left(L^2_{\rm per}(\Omega)\right)$ are called the Bloch fibers of $\mA$.
	
	The Bloch decomposition \eqref{eq:bloch_decomp} can also be extended to \emph{unbounded}, $\mathbb{L}$-periodic self-adjoint operators such as the one-body electronic Hamiltonian defined through Equation \eqref{eq:Hamiltonian}. In this case, the fibers $\mH_{\bold{k}} , ~\bold{k}\in \mathbb{R}^d$ of the electronic Hamiltonian $\mH$ are \emph{unbounded} operators on $L^2_{\rm per}(\Omega)$ given by
	\begin{equation}\label{eq:fiber}
		\mH_\bold{k} :=  \frac{1}{2}\left(-\imath \nabla + \bold{k} \right)^2 + V, \qquad \text{with domain } H^2_{\rm per}(\Omega).
	\end{equation}
	A detailed proof of this technical result can be found in \cite[Chapter XIII]{MR0493421}.

	Thanks to the Bloch-Floquet decomposition \eqref{eq:bloch_decomp}, the spectral properties of the Hamiltonian $\mH$  can be deduced using properties of the fibers $\mH_{\bold{k}}, ~{\bold{k}}\in \mathbb{R}^d$. Indeed, it is a classical result (see, e.g., \cite[Chapter XIII]{MR0493421})  that
	
	\begin{itemize}
		\item each $\mH_\bold{k}$ is a self-adjoint operator on $L^2_{\rm per}(\Omega)$ with domain $H^2_{\rm per}(\Omega)$ and form domain $H^1_{\rm per}(\Omega)$. Additionally, each $\mH_\bold{k}$ is bounded below and has compact resolvent so that each $\mH_\bold{k}$ has a purely discrete spectrum with eigenvalues accumulating at $+\infty$ and eigenfunctions that form an orthonormal basis for~$L^2_{\rm per}(\Omega)$;
		
		\item $\mH$ is a bounded-below, self-adjoint operator on $L^2(\mathbb{R}^d)$ with domain $H^2(\mathbb{R}^d)$ and form domain $H^1(\mathbb{R}^d)$. Additionally, $\mH$ has a purely absolutely continuous spectrum, and it holds that $\sigma(\mH) = \sigma_{\rm ac}(\mH)= \underset{\bold{k} \in \overline{\Omega^*}}{\cup} \sigma(\mH_{\bold{k}})$.
	\end{itemize}
	
	From the point of view of applications, the Bloch-Floquet decomposition \eqref{eq:bloch_decomp} also allows for the calculation of electronic properties of interest of the underlying crystal using only spectral information from the fibers $\mH_\bold{k}, ~\bold{k} \in \mathbb{R}^d$ (see below and also, e.g., \cites{MR0493421, MR4071823} for details). As a consequence, it suffices to focus our attention purely on the resolution of the eigenvalue problem for the operators $\mH_\bold{k}, ~\bold{k}\in \mathbb{R}^d$ defined through Equation \eqref{eq:fiber}.

	Given a fiber $\mH_\bold{k}, ~\bold{k}\in \mathbb{R}^d$ defined through Equation \eqref{eq:fiber}, we seek an $L^2_{\rm per}(\Omega)$-orthonormal basis $\left\{\left(\varepsilon_{n, \bold{k}}, u_{n, \bold{k}}\right)\right\}_{n \in \mathbb{N}^*} \subset \big(\mathbb{R} \times L^2_{\rm per}(\Omega)\big)^{\mathbb{N}^*}$ of eigenmodes of $\mH_{\bold{k}}$:
	\begin{equation}\label{eq:eigenvalue_problem}
		\begin{split}
			\mH_\bold{k} u_{n, \bold{k}}&= \varepsilon_{n, \bold{k}} u_{n, \bold{k}} \qquad \text{and}\\
			\left(u_{n, \bold{k}}, u_{m, \bold{k}}\right)_{L^2_{\rm per}(\Omega)}&= \delta_{nm} \qquad  \qquad \forall n,m \in \mathbb{N}^*.
		\end{split}
	\end{equation}

	Equipped with such a basis, we can introduce several important electronic properties of interest.  To this end, we first require a convention and some notation.
	
	\begin{Convention} \label{Conv:1}
		Consider the setting of the eigenvalue problem \eqref{eq:eigenvalue_problem}. By convention, for every $\bold{k} \in \mathbb{R}^d$ we order the eigenvalues $\varepsilon_{n, \bold{k}}, ~n \in \mathbb{N}^*$ (counting multiplicities) such that
		\begin{align*}
			\varepsilon_{1, \bold{k}} \leq \varepsilon_{2, \bold{k}} \leq \varepsilon_{3, \bold{k}} \leq \varepsilon_{4, \bold{k}} \ldots.
		\end{align*}
		Moreover, for every $n \in \mathbb{N}^*$ we will write $\varepsilon_n\colon \mathbb{R}^d \rightarrow \mathbb{R}$ for the mapping $\bold{k} \mapsto \varepsilon_{n, \bold{k}}$, and we will call $\varepsilon_n$ the $n^{\rm th}$ energy band. Since the functions $\varepsilon_n, ~n \in \mathbb{N}^*$ are continuous (as a straightforward consequence of the Courant-Fisher min-max theorem), we have
		\begin{align*}
			\sigma (\mH) = \underset{n \in \mathbb{N}^*}{\cup} \text{\rm Ran}\;(\varepsilon_n),
		\end{align*}
		and for each $n \in \mathbb{N}^*$ it holds that $\text{\rm Ran}\;(\varepsilon_n)= [\min \varepsilon_n, \max \varepsilon_n]$ is an interval.
		
		We will use a similar convention for any subsequent eigenvalue problem that we introduce in the sequel.
	\end{Convention}

	The energy bands play a key role in the definition of various electronic properties of a perfect crystal. Indeed, given $\bold{k} \in \mathbb{R}^d$ and the Bloch fiber $\mH_{\bold{k}}$ defined through Equation \eqref{eq:fiber}, the $\bold{k}^{\rm th}$ fiber of the \textbf{one-body ground-state density matrix} at chemical potential $\mu\in \mathbb{R}$ is defined as the bounded self-adjoint operator $\gamma_{\bold{k}} \colon L^2_{\rm per}(\Omega) \rightarrow L^2_{\rm per}(\Omega)$ given~by
	\begin{align*}
		\gamma_{\bold{k}}:= \mathbbm{1}(\mH_{\bold{k}} \leq \mu) = \sum_{n \in \mathbb{N}^*} \mathbbm{1}\big(\varepsilon_{n}(\bold{k}) \leq \mu\big) \ket{u_{n, \bold{k}}}\bra{u_{n, \bold{k}}}.
	\end{align*}
	
	The \textbf{integrated density of states} is defined as the function $\mathcal{N} \colon \mathbb{R} \rightarrow \mathbb{R}_+$ with the property that
	\begin{align*}
		\forall \mu \in \mathbb{R} \colon \hspace{5mm} \mathcal{N}(\mu) := \sum_{n \in \mathbb{N}^*} \fint_{\Omega^*} \mathbbm{1}\big(\varepsilon_{n}(\bold{k}) \leq \mu\big)d \bold{k}.
	\end{align*}
	
	Lastly, the \textbf{integrated density of energy} is defined as the function $\mathcal{E} \colon \mathbb{R} \rightarrow \mathbb{R}$ with the property that
	\begin{align*}
		\forall \mu \in \mathbb{R} \colon \hspace{5mm} \mathcal{E}(\mu):= \sum_{n \in \mathbb{N}^*} \fint_{\Omega^*} \varepsilon_n(\bold{k})\mathbbm{1}\big(\varepsilon_{n}(\bold{k}) \leq \mu\big)d \bold{k}.
	\end{align*}
	
	Often the above quantities are computed for $\mu= \mu_{\rm F} \in \mathbb{R}$ where $\mu_{\rm F}$ is known as the \emph{Fermi level} and is defined through the relation $\mathcal{N}(\mu_{\rm F})= N$ with $N$ being the number of electrons (or electrons pairs if spin is taken into account) per unit cell. Naturally, computing any of these physical observables requires the approximation, through numerical quadrature, of integrals over the Brillouin zone $\Omega^*$ that involve the energy bands~$\{\varepsilon_{n}\}_{n \in \mathbb{N}^*}$. This is a highly non-trivial problem in the case of \emph{metallic} systems for which the Fermi level $\mu_{\rm F}$ is an interior point of $\sigma(\mH)$, and several numerical methods have been proposed for Brillouin zone integration (see, for instance, the previously cited articles \cites{lehmann1972numerical, blochl1994improved, morgan2018efficiency, pickard1999extrapolative, henk2001integration}). From the point of view of numerical analysis, it is natural to ask for error bounds for the various numerical methods in the literature, and such an error analysis has recently been carried out in \cite{MR4071823} under the assumption that the values of the functions $\varepsilon_n, ~n \in \mathbb{N}^*$ can be computed exactly at any $\bold{k} \in \Omega^*$.
	
	As is typically the case in the analysis of quadrature methods, the error analysis in \cite{MR4071823} makes use of functional properties of the energy bands $\{\varepsilon_{n}\}_{n \in \mathbb{N}^*}$. This analysis shows that there are two properties of these energy bands in particular that are necessary in order to deduce higher order convergence rates for numerical quadrature in the Brillouin zone. \vspace{4mm}

	\noindent\textbf{Property one} {(Periodicity of the eigenvalues)}\textbf{.}~
	
	Consider the setting of the eigenvalue problem~\eqref{eq:eigenvalue_problem} and let Convention \ref{Conv:1} hold.  Then for each $n \in \mathbb{N}^*$, the function $\varepsilon_{n}$ is~$\mathbb{L}^*$-periodic.
	
	The proof follows in view of Convention \ref{Conv:1} by recognizing that for any $\bold{k} \in \mathbb{R}^d$ and any $\bold{G} \in \mathbb{L}^*$, the operators $\mH_{\bold{k}}$ and $\mH_{\bold{k} +\bold{G}}$ are unitarily equivalent through the  unitary multiplication operator $T_{\bold{G}} \colon L^2_{\rm per}(\Omega) \rightarrow  L^2_{\rm per}(\Omega) $ defined in Section \ref{sec:2b}. \vspace{4mm}
	
	\noindent\textbf{Property two} (Continuity of the eigenvalues)\textbf{.}~
	
	Consider the setting of the eigenvalue problem \eqref{eq:eigenvalue_problem}, and let the maps $\{\varepsilon_{n}\}_{n \in \mathbb{N}^*}$ be defined according to Convention~\ref{Conv:1}. Then, each function $\varepsilon_{n}$ is Lipschitz continuous on $\mathbb{R}^d$. Additionally, if $\bold{k}_n \in \mathbb{R}^d$ is such that 
	\begin{align*}
		\varepsilon_{n,\bold{k}_n} \neq \varepsilon_{m,\bold{k_n}}  \quad\forall  ~\mathbb{N}^* \ni m\neq n,  \hspace{1cm} \big(\text{No energy band crossings at } (\bold{k}_n, \varepsilon_{n,\bold{k}_n})\big),
	\end{align*}
	then $\varepsilon_{n}$ is locally real-analytic at $\bold{k}_n$, i.e.,  $\exists \delta_n >0$ such that $\varepsilon_n$ is real-analytic on the open ball $\mathbb{B}_{\delta_n}(\bold{k}_n)$.
	
	A proof of this statement can, for instance, be found in \cite[Lemma 3.2]{MR4071823}.  \vspace{4mm}

	The $\mathbb{L}^*$-periodicity and real-analyticity away from crossings of the energy bands $\{\varepsilon_{n}\}_{n \in \mathbb{N}^*}$ has significant consequences for evaluating Brillouin zone integrals involving these functions. Indeed as can readily be deduced from~\cite{MR2664609}, for $d$-dimensional periodic integrands of class $\mathscr{C}^{r}$, the uniform grid quadrature rule converges as $\mathcal{O}\Big((\Delta x)^\frac{r}{d}\Big)$ when integrating over an entire period. For real-analytic periodic integrands, a uniform grid quadrature rule even recovers exponential convergence (see, e.g., \cite{MR3245858}). This fact is essential in understanding the approximability of the integrals appearing in the definitions of the various electronic properties of interest defined above, and is a key element of the higher order convergence rates for numerical quadrature obtained in~\cite{MR4071823}.
	
	Of course in practice, we typically do not have access to the \emph{exact} energy bands $\{\varepsilon_{n}\}_{n \in \mathbb{N}^*}$, these being solutions to infinite-dimensional eigenvalue problems. Instead, the eigenvalue problem \eqref{eq:eigenvalue_problem} is typically discretized in some $M$-dimensional basis for specific values of $\bold{k} \in \Omega^* \subset \mathbb{R}^d$ corresponding to the grid points of our chosen quadrature method. This naturally raises the question of how the resulting \emph{approximate} energy bands $\{\varepsilon_n^{\rm approx}\}_{i=1}^{M}$ compare to the exact bands $\{\varepsilon_{n}\}_{n \in \mathbb{N}^*}$, and in particular whether \textbf{Properties 1} and \textbf{2} also hold for the approximate bands $\{\varepsilon_n^{\rm approx}\}_{i=1}^{M}$, these properties being essential to the quadrature error analysis. This is the topic of the next section.
	
	\section{Classical Discretization Strategies}\label{sec:3}
	We will now describe two well-known discretization strategies for resolving the eigenvalue problem \eqref{eq:eigenvalue_problem}.  Throughout this section, we assume the setting of Section \ref{sec:2}, and we recall in particular that we will use Convention \ref{Conv:1} in ordering and labelling all eigenvalue problems that appear in this section.

	\begin{definition}[Uniform basis set]\label{def:basis_1}~
		
		Let ${\rm E_c}>0$ denote a scalar cutoff. We define the basis set~$\mathcal{B}_0^{\rm E_c} \subset H^1_{\rm per}(\Omega)$ as
		\begin{align*}
			\mathcal{B}_0^{\rm E_c}:= \left\{ e_{\bold{G}}\colon \hspace{2mm} \bold{G} \in \mathbb{L}^*  ~\text{ with} \quad \frac{1}{2} \vert \bold{G} \vert^2 < {\rm E_c}\right\},
		\end{align*}
		and we define the subspace spanned by this basis set as $\mX_0^{\rm E_c}:={\rm span}\;\mathcal{B}_0^{\rm E_c}$. 
	\end{definition}
	
	\begin{Notation}[Projections involving the uniform basis set]\label{def:proj_1}~
		
		Consider the setting of Definition \ref{def:basis_1}. We denote by $\Pi_{{\rm E_c}} \colon L^2_{\rm per}(\Omega) \rightarrow L^2_{\rm per}(\Omega)$ the $L^2_{\rm per}$-orthogonal projection operator onto $\mX_0^{{\rm E_c}}$, and we denote by $\Pi_{{\rm E_c}}^{\perp}$ its complement, i.e., $\Pi_{{\rm E_c}}^{\perp}:= \mathbb{Id}- \Pi_{{\rm E_c}}$.
		
		Additionally, for each $\bold{k} \in \mathbb{R}^d$ we denote by ${\breve{\mH}^{{\rm E_c}}_{\bold{k}}}$ the two-sided projection of the Hamiltonian fiber $\mH_\bold{k}$ in $\mX_0^{{\rm E_c}}$, i.e.,~${\breve{\mH}^{{\rm E_c}}_{\bold{k}}}:= \Pi_{{\rm E_c}}\mH_\bold{k}\Pi_{{\rm E_c}}$.
	\end{Notation}

	Equipped with the uniform basis sets defined through Definition \ref{def:basis_1}, we can propose the following elementary Galerkin discretization of Eigenvalue problem~\eqref{eq:eigenvalue_problem}.
	
	\vspace{0.4cm}
	
\noindent	\textbf{Uniform Galerkin discretion of the eigenvalue problem \eqref{eq:eigenvalue_problem}.}~
	
	Given a fiber $\mH_\bold{k}, ~\bold{k}\in \mathbb{R}^d$ defined through Equation \eqref{eq:fiber} and a scalar cutoff ${\rm E_c} >0$, we seek an orthonormal basis $\big\{\big(\breve{\varepsilon}^{\rm E_c}_{n, \bold{k}}, \breve{u}^{\rm E_c}_{n, \bold{k}}\big)\big\}\subset \mathbb{R} \times \mX_0^{\rm E_c}$ of eigenvectors of~$\breve{\mH}^{{\rm E_c}}_\bold{k}$:
	\begin{equation}\label{eq:Galerkin_1}
		\begin{split}
			\breve{\mH}^{{\rm E_c}}_\bold{k}  \breve{u}^{\rm E_c}_{n, \bold{k}}&= \breve{\varepsilon}^{\rm E_c}_{n, \bold{k}} \breve{u}^{\rm E_c}_{n, \bold{k}} \hspace{8mm} \text{and}\\[0.5em]
			\big(\breve{u}^{\rm E_c}_{n, \bold{k}}, \breve{u}^{\rm E_c}_{m, \bold{k}}\big)_{L^2_{\rm per}(\Omega)}&= \delta_{nm} \qquad \hspace{8mm}  \forall n,m \in \{1, \ldots, \text{dim}\;\breve{\mH}^{{\rm E_c}}_\bold{k} \}.
		\end{split}
	\end{equation}
	
	\vspace{4mm}
	
	An alternative to the uniform Galerkin discretization \eqref{eq:Galerkin_1} is provided by the use of so-called $\bold{k}$-dependent basis sets.

	\begin{definition}[$\bold{k}$-dependent basis set]\label{def:basis_2}~
		
		Let ${\rm E_c} >0$ denote a scalar cutoff, and let $\bold{k} \in \mathbb{R}^d$. We define the basis set~$\mathcal{B}_{\bold{k}}^{\rm E_c} \subset H^1_{\rm per}(\Omega)$ as 
		\begin{align*}
			\mathcal{B}_{\bold{k}}^{\rm E_c}:= \left\{ e_{\bold{G}}\colon \hspace{2mm} \bold{G} \in \mathbb{L}^* ~ \text{ with } \quad \frac{1}{2} \vert \bold{k}+\bold{G} \vert^2 < {\rm E_c}\right\},
		\end{align*}
		and we define the subspace spanned by this basis set as $\mX_{\bold{k}}^{\rm E_c}:={\rm span}\;\mathcal{B}_{\bold{k}}^{\rm E_c}$. Additionally, we write $M_{\rm E_c}(\bold{k})$ to denote the cardinality of $\mathcal{B}_{\bold{k}}^{\rm E_c}$, and we refer to $\mathcal{B}_{\bold{k}}^{E_c}$ as a $\bold{k}$-dependent basis set.
	\end{definition}
	
	\begin{remark}\label{rem:basis_2}
		Consider the setting of Definition \ref{def:basis_2}. It can readily be seen that for a fixed ${\rm E_c}$, the cardinality $M_{\rm E_c}(\bold{k})$ of the basis $\mathcal{B}_{\bold{k}}^{\rm E_c}$ is not fixed and depends indeed on $\bold{k} \in {\mathbb{R}^d}$. In the sequel, we will therefore regard $M_{\rm E_c}(\cdot)$ as a piecewise constant mapping from $\mathbb{R}^d$ to $\mathbb{N}^*$, which is moreover uniformly bounded below and above by optimal constants $M_{\rm E_c}^-$ and $M_{\rm E_c}^+$ respectively that depend on ${\rm E_c}$. 
	\end{remark}
	
	\begin{Notation}[Projections involving the $\bold{k}$-dependent basis set]\label{def:proj_2}~
		
		Consider the setting of Definition \ref{def:basis_2}. For each $\bold{k} \in \mathbb{R}^d$, we denote by $\Pi_{\bold{k}, {\rm E_c}} \colon L^2_{\rm per}(\Omega) \rightarrow L^2_{\rm per}(\Omega)$ the $L^2_{\rm per}$-orthogonal projection operator onto $\mX_{\bold{k}}^{{\rm E_c}}$, and we denote by $\Pi_{\bold{k}, {\rm E_c}}^{\perp}$ its complement, i.e., $\Pi_{\bold{k}, {\rm E_c}}^{\perp}:= \mathbb{Id}- \Pi_{\bold{k}, {\rm E_c}}$.
		
		Additionally, for each $\bold{k} \in \mathbb{R}^d$ we denote by $ {\mH^{{\rm E_c}}_{\bold{k}}}$ the two-sided projection of the Hamiltonian fiber $\mH_\bold{k}$ in $\mX_{\bold{k}}^{{\rm E_c}}$, i.e.,~$ {\mH^{{\rm E_c}}_{\bold{k}}}:= \Pi_{\bold{k},{\rm E_c}}\mH_\bold{k}\Pi_{\bold{k},{\rm E_c}}$.
	\end{Notation}

	\vspace{0.4cm}
	
\noindent	\textbf{$\bold{k}$-dependent Galerkin discretization of the eigenvalue problem \eqref{eq:eigenvalue_problem}.}~
	
	Given a fiber $\mH_\bold{k}, ~\bold{k}\in \mathbb{R}^d$ defined through Equation \eqref{eq:fiber} and a scalar cutoff ${\rm E_c} >0$, we seek an orthonormal basis $\big\{\big( {\varepsilon^{\rm E_c}_{n, \bold{k}}},  {u^{\rm E_c}_{n, \bold{k}}}\big)\big\}\subset \mathbb{R} \times \mX_{\bold{k}}^{\rm E_c}$ of eigenvectors of ${\mH^{{\rm E_c}}_\bold{k}}$:
	\begin{equation}\label{eq:Galerkin_2}
		\begin{split}
			{\mH^{{\rm E_c}}_\bold{k}}   {u^{\rm E_c}_{n, \bold{k}}}&=  {\varepsilon^{\rm E_c}_{n, \bold{k}}}  {u^{\rm E_c}_{n, \bold{k}}} \hspace{8mm} \text{and}\\[0.5em]
			\big( {u^{\rm E_c}_{n, \bold{k}}},  {u^{\rm E_c}_{m, \bold{k}}}\big)_{L^2_{\rm per}(\Omega)}&= \delta_{nm} \qquad \hspace{8mm}  \forall n,m \in \{1, \ldots, M_{\rm E_c}(\bold{k})\}.
		\end{split}
	\end{equation}
	
	\vspace{4mm}

	%

	The Galerkin discretizations \eqref{eq:Galerkin_1} and \eqref{eq:Galerkin_2} are both well-posed, and a straightforward analysis reveals the following error bound: for any fixed $\bold{k} \in \mathbb{R}^d$, any $n \in \mathbb{N}^*$, there exists ${\rm E^*} >0$ such that for scalar cutoffs ${\rm E_c}\geq {\rm E}^*$, we have eigenvalue bounds of the form:
	\begin{align}\label{eq:Error_1}
		\vert \breve{\varepsilon^{\rm E_c}_{n, \bold{k}}}  - \varepsilon_{n \bold{k}} \vert  \lesssim {\rm \big(E_c}\big)^{-s} \qquad \text{and} \qquad 	\vert  {\varepsilon^{\rm E_c}_{n, \bold{k}}}  - \varepsilon_{n \bold{k}} \vert  \lesssim {\rm \big(E_c}\big)^{-s},
	\end{align}
	where the convergence rate $s \geq 0$ depends on the regularity of the effective potential $V\in L^{\infty}_{\rm per}(\Omega)$.
	
\begin{figure}[h!]
		\centering
		\begin{subfigure}{0.495\textwidth}
			\centering
			\includegraphics[width=\textwidth, trim={0cm, 2cm, 0cm, 0cm},clip=true]{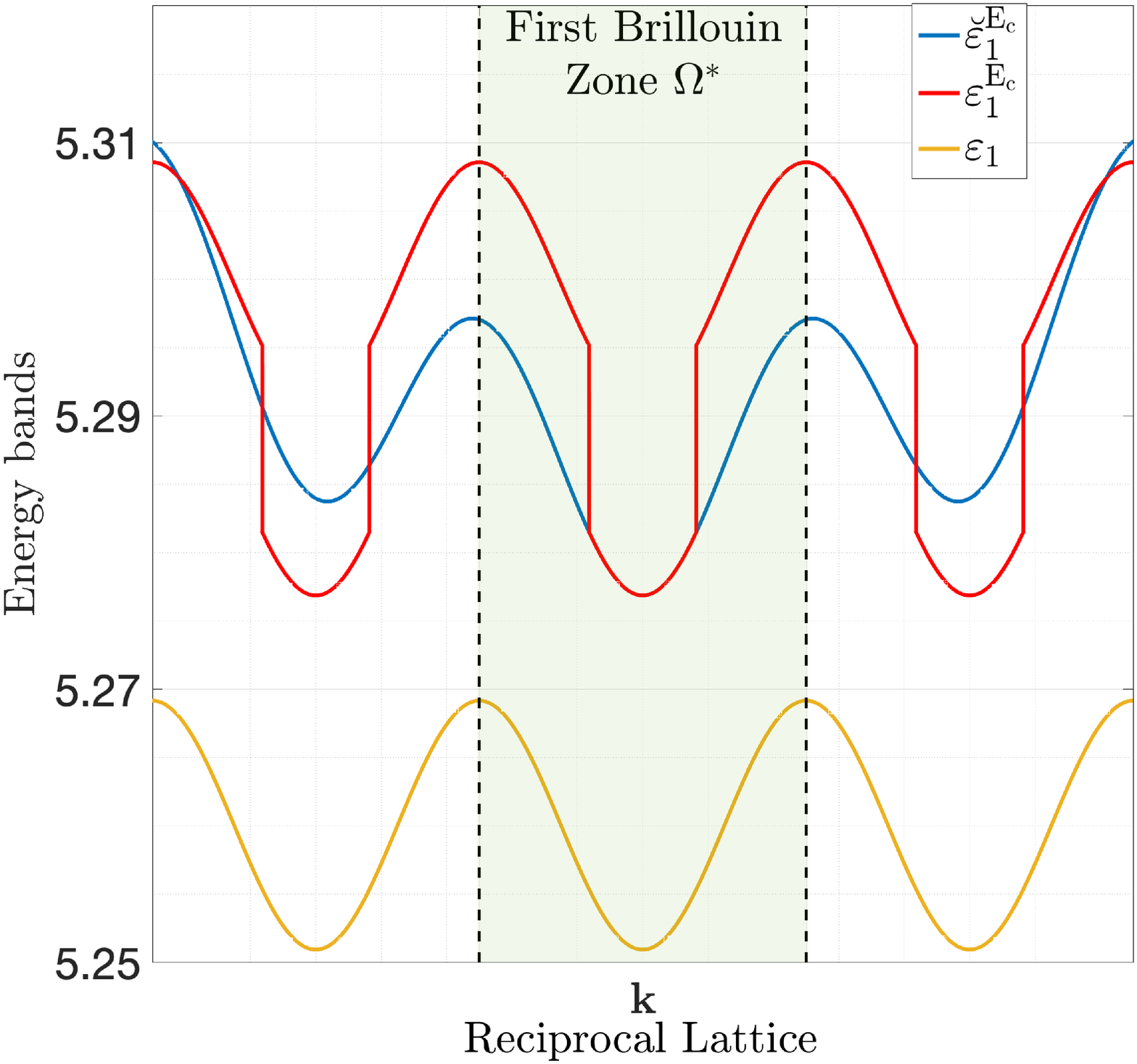} 
			\caption{The approximate and exact energy bands}		
			\label{fig:01a}
		\end{subfigure}\hfill
		\begin{subfigure}{0.495\textwidth}
			\centering
			\includegraphics[width=\textwidth, trim={0cm, 2cm, 0cm, 0cm},clip=true]{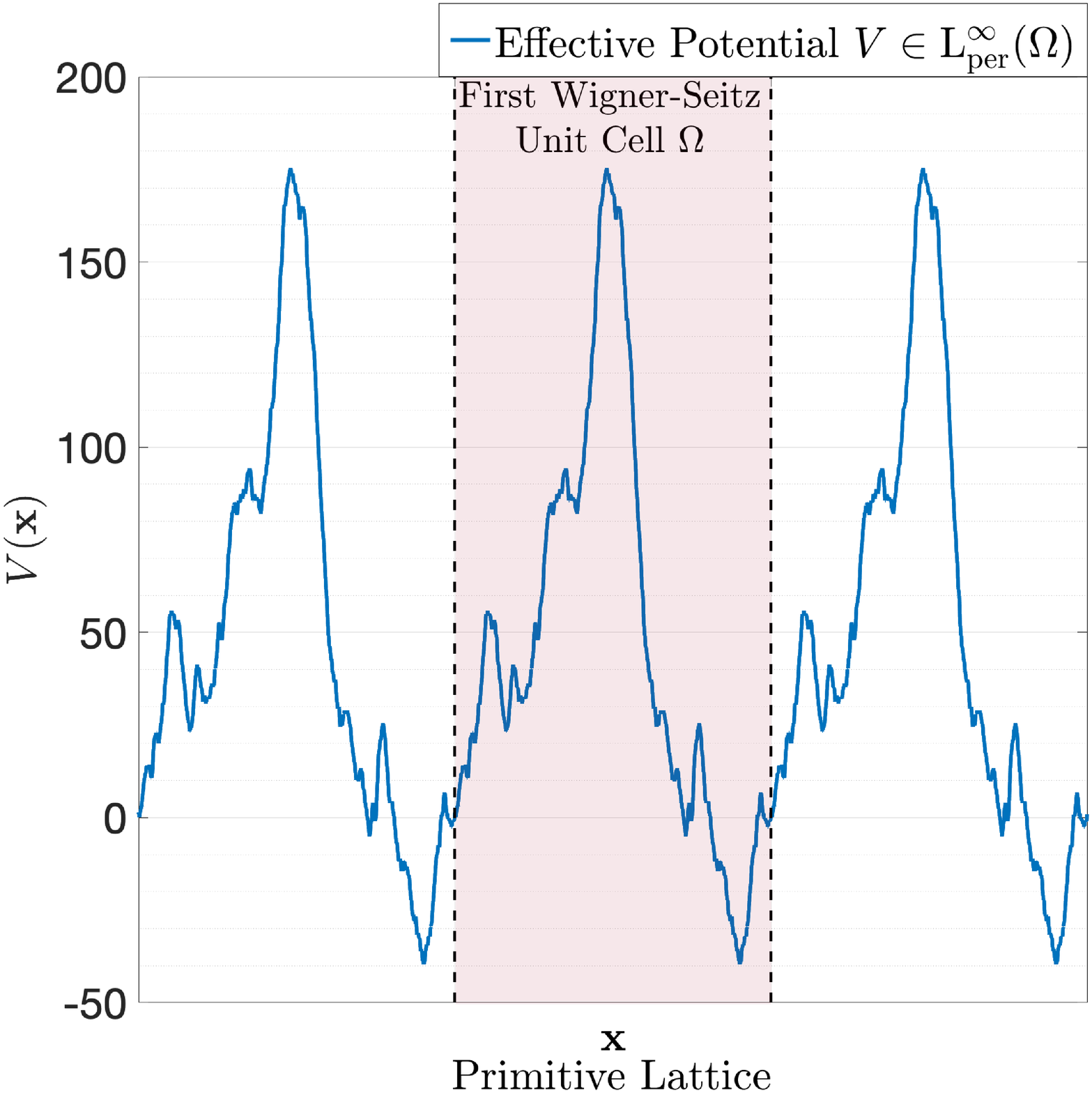} 
			\caption{The $\mathbb{L}$-periodic potential $V \in L_{\rm per}^{\infty}(\Omega)$}
			\label{fig:01b}
		\end{subfigure}
		\caption{Lowest energy bands for a simple 1-D example with effective potential $V$ as shown.}
		\label{fig:01}
	\end{figure}

	Unfortunately, in spite of the availability of the error estimate \eqref{eq:Error_1}, a closer study of the Galerkin discretizations~ \eqref{eq:Galerkin_1} and \eqref{eq:Galerkin_2} reveals a serious deficiency that may not have been immediately apparent: the approximate energy bands $\left\{ \breve{\varepsilon}^{\rm E_c}_{n}\right\}_{n \in M_{\rm E_c}(0)} $ and $\left\{  {\varepsilon^{\rm E_c}_{n}}\right\}_{n \in M_{\rm E_c}(k)} $ do not preserve \textbf{Properties 1} and \textbf{2} of the exact energy bands $\{\varepsilon_{n}\}_{n \in \mathbb{N}^*}$ respectively. An example of this phenomenon is displayed in Figure \ref{fig:01a} where we plot the exact and approximate ground state energy bands for a simple one-dimensional example.
	
	The core problem is that while the exact fibers $\mH_{\bold{k}}$ and $\mH_{\bold{k}+\bold{G}}$ are unitarily equivalent for all $\bold{k} \in \mathbb{R}^d$ and $\bold{G} \in \mathbb{L}^*$, the same is not always true for the uniform-basis projected fibers. Indeed, $\breve{\mH}_{\bold{k}}^{\rm E_c}=\Pi_{{\rm E_c}}\mH_\bold{k}\Pi_{{\rm E_c}}$ is not, in general, unitarily  equivalent to $\breve{\mH}_{\bold{k}+\bold{G}}^{\rm E_c}$ as can readily be verified by a direct calculation in the case $V\equiv 0$. Unitary equivalence is, conversely, preserved for the $\bold{k}$-basis projected fibers ${\mH}_{\bold{k}}^{\rm E_c}=\Pi_{\bold{k}, {\rm E_c}}\mH_\bold{k}\Pi_{\bold{k}, {\rm E_c}}$ but in this case, the rank of  ${\mH}_{\bold{k}}^{\rm E_c}$ changes as a function of $\bold{k}$ (recall Remark \ref{rem:basis_2}). This causes the continuity argument used to prove~\textbf{Property~2} to break down.
	
	In other words, the choice  of basis set (uniform or $\bold{k}$-dependent) represents a trade-off between the $\mathbb{L}^*$-periodicity and regularity of the resulting approximate energy bands. Both choices are obviously sub-optimal from the point of view of numerical quadrature in the Brillouin zone, and it is therefore of great interest to develop an alternative discretization scheme that might result in approximate energy bands that are both $\mathbb{L}^*$-periodic and of class $\mathscr{C}^r$ for some $r \geq 0$. One such methodology which has been proposed in the physics literature (see \cite{bernasconi1995first} for the original proposal and \cite{janssen2016precise} for a recent presentation) and also implemented in several quantum chemistry simulation softwares (see \cites{Abinit, Qbox}) relies on the idea of modifying the `diagonal' part of the fiber $\mH_{\bold{k}}$ in a controlled manner. A systematic description of this discretization method is the subject of the next section.

	\section{Operator Modification Approach}\label{sec:4}

	We begin this section by defining a one-dimensional ``blow-up'' function that will be central to the construction of a modified Hamiltonian operator. Throughout this section, we will assume the settings of Sections \ref{sec:2} and \ref{sec:3}.

	\begin{definition}[Blow-up function]\label{def:g}~
	
		Let $m \in \mathbb{N}$, let $\mathcal{f}_2$ denote the quadratic monomial, i.e., $\mathcal{f}_2(x)=x^2$ for all $x \in \mathbb{R}$, and denote by~$\mathcal{h}\colon \left[\frac{1}{2}, 1\right)\rightarrow \mathbb{R}$ a function with the following four properties: \vspace{2mm}
		
		\begin{enumerate}
			\item It holds that $\mathcal{h} \in \mathscr{C}^{m}\left(\frac{1}{2}, 1\right) $. \vspace{2mm}
			
			\item It holds that $\lim_{x \to 1^-}\;\left((1-x)^m\; \mathcal{h}(x)\right) = +\infty$.\vspace{2mm}
			
			\item It holds that $\mathcal{h}(x) = \mathcal{f}_2(x)$ on $ \left[0,\frac{1}{2}\right]$ and $ \mathcal{h}(x) \ge \mathcal{f}_2(x)$ for all $x \in \left(\frac{1}{2}, 1\right)$.\vspace{2mm}
			
			\item For all $j \in \{0, \ldots, m\}$ it holds that $\mathcal{h}^{(j)}\left(\frac{1}{2}\right)= \mathcal{f}_2^{(j)}\left(\frac{1}{2}\right)$, where $\mathcal{h}^{(j)}(\cdot )$ and $\mathcal{f}_2^{(j)}(\cdot)$ denote the $j^{\rm th}$ derivative of $\mathcal{h}$ and  $\mathcal{f}_2$ respectively. \vspace{2mm}
		\end{enumerate}
		Then we define the blow-up function $\mathscr{G} \colon \mathbb{R} \rightarrow \mathbb{R}$ as the mapping given by
		\begin{equation}\label{eq:g}
			\mathscr{G}(x) = \begin{cases}
				\mathcal{f}_2(x) \quad &\text{for } \vert x \vert \in [0, \nicefrac{1}{2}] ~\cup ~ [1, \infty) ,\\[0.5em]
				\mathcal{h}(x) \quad &\text{for } \vert x \vert \in\left (\frac{1}{2}, 1\right),
			\end{cases}
		\end{equation}
		where we have suppressed the dependency of $\mathscr{G}$ on ${m}$ and $\mathcal{h}$ by assuming once and for all that $m$ and $\mathcal{h}$ are fixed for the remainder of our analysis.
	\end{definition}

	Consider Definition \ref{def:g} of the blow-up function $\mathscr{G}$. We emphasize three properties of $\mathscr{G}$ that will be useful in the sequel: first, that it is of class $\mathscr{C}^m$ on the interval $[0, 1)\subset \mathbb{R}$; second, that it is \emph{point-wise} bounded below by the quadratic map $x \mapsto x^2$ on all of $\mathbb{R}$, and third that it blows up as $x\to 1^-$ at a rate greater than $\frac{1}{(1-x)^{m}}$.

	\begin{remark}[Blow-up functions in the physics literature]\label{rem:literature}~
		
		It is pertinent at this point to contrast our rather general definition of the blow-up function $\mathscr{G}$ with those that have been proposed in the literature (see \cites{bernasconi1995first, janssen2016precise}) and implemented in electronic structure calculation codes (see \cites{Abinit, Qbox}).
		
		In fact, the functions $\widetilde{\mathscr{G}}$ proposed in \cites{bernasconi1995first, janssen2016precise} are not `blow-up' functions at all, in the sense that $\lim_{x \to 1^-} \widetilde{\mathscr{G} } (x)\neq +\infty$. Instead, both papers propose the use of the error function to construct $\widetilde{\mathscr{G}}$ such that $\lim_{x \to 1} \widetilde{\mathscr{G} } (x)= c \gg 1$ but with $c < \infty$. The implementation in the quantum chemistry code QBOX \cite{Qbox} is based on similar ideas. In contrast, the software suite ABINIT \cite{Abinit} employs a true `blow-up' function that satisfies the conditions of Definition \ref{def:g} for $m=1$. 
	\end{remark}

	We will now propose a modified Galerkin discretization for the eigenvalue problem \eqref{eq:eigenvalue_problem}. To this end, we first require a definition, and we recall in particular Definition \ref{def:basis_2} of the $\bold{k}$-dependent basis set on $L^2_{\rm per}(\Omega)$ and Notation \ref{def:proj_2}. 
	
	\begin{definition}[Modified Hamiltonian operator]\label{def:Mod_Hamiltonian}~
		
		Let ${\rm E_c}>0$ and let the blow-up function $\mathscr{G}$ be defined according to Equation \eqref{eq:g}. For each $\bold{k} \in \mathbb{R}^d$, we define the operator $\widetilde{\mH}_{\bold{k}}^{\mathscr{G}, {\rm E_c}} \colon {\rm X}_{\bold{k}}^{\rm E_c} \rightarrow {\rm X}_{\bold{k}}^{\rm E_c}$ as the mapping with the property that
		\begin{equation}\label{eq:fiber_mod}
			\widetilde{\mH}^{\mathscr{G}, \rm E_c}_\bold{k}:= \Pi_{\bold{k}, {\rm E_c}} \left(\; {\rm E_c}\; \mathscr{G}\left(\frac{\vert-\imath \nabla + \bold{k}\vert}{\sqrt{2 {\rm E_c}}} \right) + V \right)\Pi_{\bold{k}, {\rm E_c}}.
		\end{equation}
	\end{definition}
	
	\begin{remark}\label{rem:fibers_2}~
		Consider the setting of Definition \ref{def:Mod_Hamiltonian}. In Equation \eqref{eq:fiber_mod}, the term $\mathscr{G}\left(\frac{\vert-\imath \nabla + \bold{k}\vert}{\sqrt{2 {\rm E_c}}} \right)$ should be understood in the sense of functional calculus. In particular, given some $\Phi \in {\rm X}_{\bold{k}}^{\rm E_c}\subset H^2_{\rm per}(\Omega)$, we have
		\begin{align*}
			\widetilde{\mH}^{\mathscr{G}, \rm E_c}_\bold{k} \Phi= \Pi_{\bold{k}, {\rm E_c}}\Bigg(\sum_{\substack{\bold{G}\in \mathbb{L}^*\\ \frac{1}{2}\vert \bold{k}+\bold{G}\vert^2 < {\rm E_c}}} \widehat{\Phi}_{\bold{G}} \; \left({\rm E_c}\mathscr{G}\left(\frac{\vert\bold{G} + \bold{k}\vert}{\sqrt{2{\rm E_c}}} \right)+V\right)\; e_{\bold{G}}\Bigg).
		\end{align*}
		
		Additionally, recalling the definition of the Bloch-Floquet fibers $\mH_{\bold{k}}, ~\bold{k} \in \mathbb{R}^d$ given by Equation \eqref{eq:fiber}, we notice that for each $\bold{k} \in \mathbb{R}^d$, thanks to the definition of the blow-up function $\mathscr{G}$, we have that $\widetilde{\mH}^{\mathscr{G}, \rm E_c}_\bold{k} \geq \mH^{\rm E_c}_{\bold{k}}:= \Pi_{\bold{k},{\rm E_c}} \mH_{\bold{k}} \Pi_{\bold{k},{\rm E_c}}$, i.e., for all $\Phi \in {\rm X}_{\bold{k}}^{\rm E_c} \subset H^1_{\rm per}(\Omega)$ it holds that
		\begin{align*}
			\left(\Phi, \widetilde{\mH}^{\mathscr{G}, \rm E_c}_\bold{k}\Phi\right)_{L^2_{\rm per}(\Omega)} \geq \left(\Phi, \mH^{\rm E_c}_{\bold{k}} \Phi\right)_{L^2_{\rm per}(\Omega)}.
		\end{align*}
		
	\end{remark}

	\vspace{4mm}	
	
	\noindent \textbf{$\bold{k}$-dependent modified Galerkin discretization of the eigenvalue problem \eqref{eq:eigenvalue_problem}}~
	
	Let ${\rm E_c}>0$, let the blow-up function $\mathscr{G}$ be defined according to Equation \eqref{eq:g}, and let the modified Hamiltonian operator $\widetilde{\mH}^{\mathscr{G}, {\rm E_c}}_\bold{k}, ~ \bold{k} \in \mathbb{R}^d$ be defined through Definition \ref{def:Mod_Hamiltonian}. We seek an $L^2_{\rm per}(\Omega)$-orthonormal basis $\big\{\big( {\widetilde{\varepsilon}^{\rm E_c}_{n, \bold{k}}},  \widetilde{u}^{\rm E_c}_{n, \bold{k}}\big)\big\}\subset \mathbb{R} \times \mX_{\bold{k}}^{\rm E_c}$ of eigenmodes of $\widetilde{\mH}^{\mathscr{G}, {\rm E_c}}_\bold{k}$:
	\begin{equation}\label{eq:Galerkin_3}
		\begin{split}
			\widetilde{\mH}^{\mathscr{G}, {\rm E_c}}_\bold{k}   \widetilde{{u}}^{\rm E_c}_{n, \bold{k}}&=  \widetilde{\varepsilon}^{\rm E_c}_{n, \bold{k}}  \widetilde{u}^{\rm E_c}_{n, \bold{k}} \hspace{8mm} \text{and}\\[0.5em]
			\big( \widetilde{u}^{\rm E_c}_{n, \bold{k}},  \widetilde{u}^{\rm E_c}_{m, \bold{k}}\big)_{L^2_{\rm per}(\Omega)}&= \delta_{nm} \qquad \hspace{8mm}  \forall n,m \in \{1, \ldots, M_{\rm E_c}(\bold{k})\}.
		\end{split}
	\end{equation}

\begin{figure}[ht]
		\centering
		\begin{subfigure}{0.495\textwidth}
			\centering
			\includegraphics[width=\textwidth, trim={0cm, 0cm, 0cm, 0cm},clip=true]{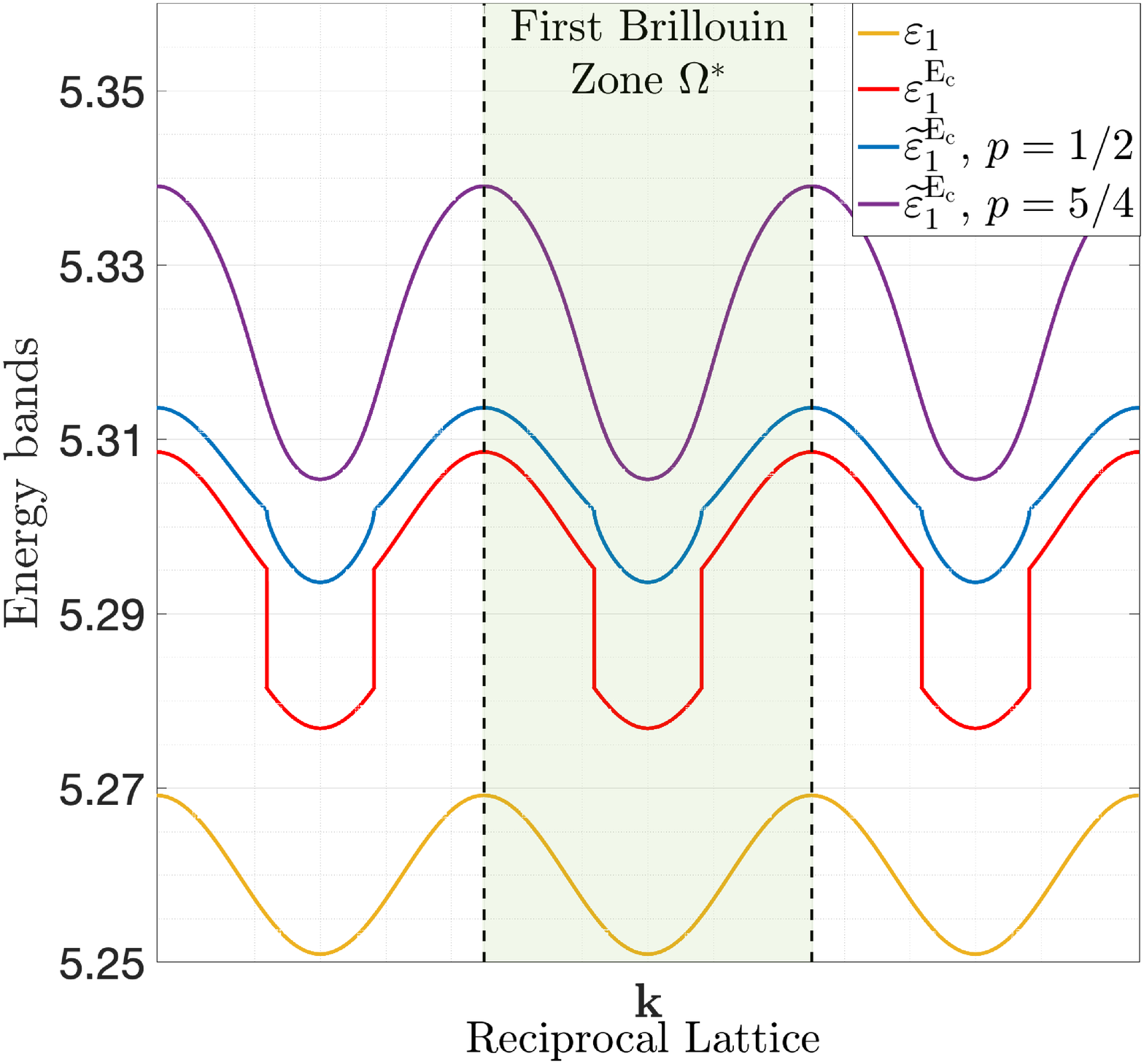} 
			\caption{The approximate and exact energy bands}		
			\label{fig:02a}
		\end{subfigure}\hfill
		\begin{subfigure}{0.495\textwidth}
			\centering
			\includegraphics[width=\textwidth, trim={0cm, 0cm, 0cm, 0cm},clip=true]{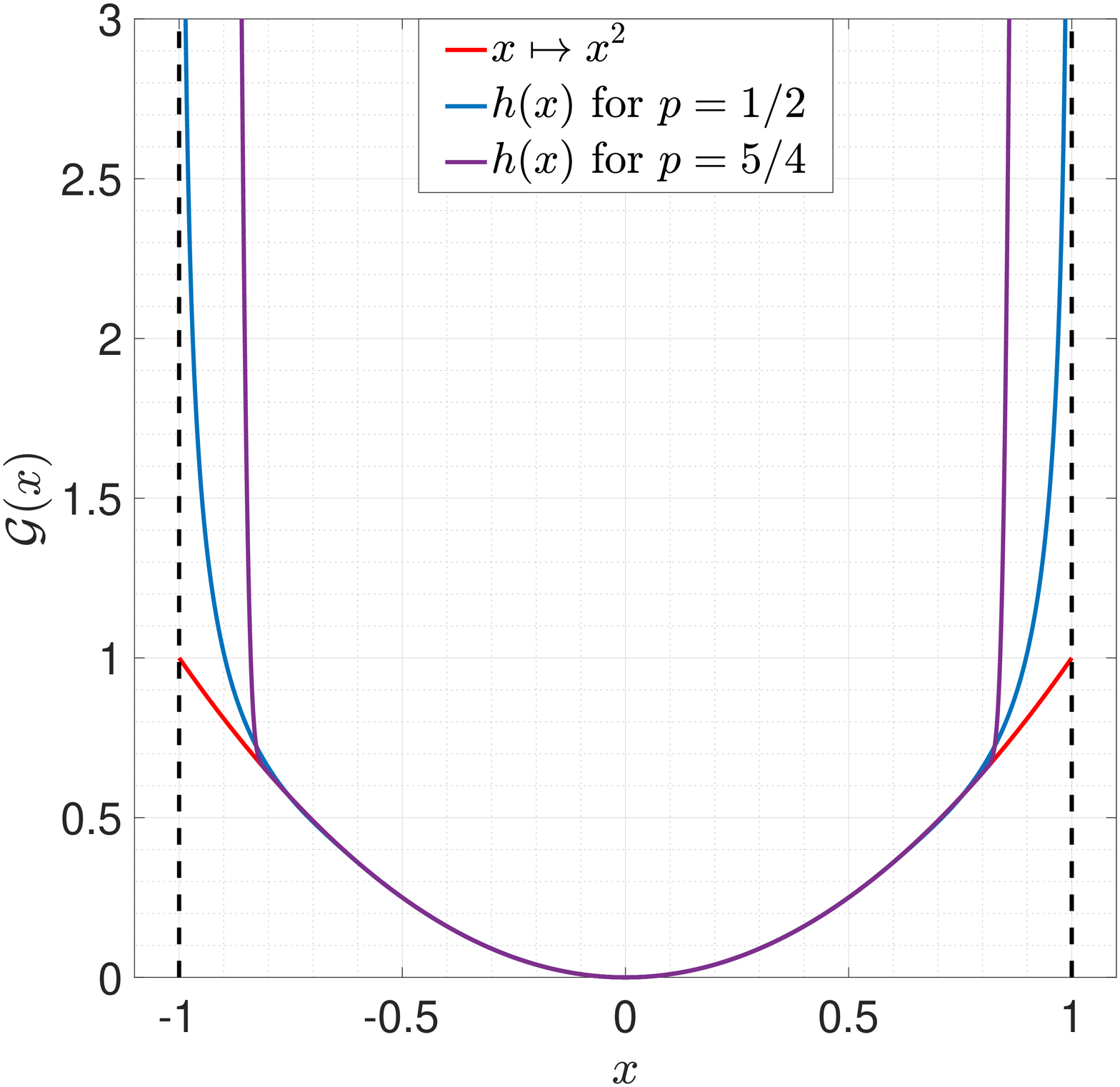} 
			\caption{Blow-up function $\mathcal{h}$ for different order singularities.}
			\label{fig:02b}
		\end{subfigure}
		\caption{Lowest energy bands for the same 1D effective potential $V$ used to produce Figure~\ref{fig:01}. The blow-up functions were of the form $\mathcal{h}(x)= {\rm C}(1-x)^{-p}$ in the vicinity of $1^-$.}
		\label{fig:2}
	\end{figure}

	The eigenvalue problem \eqref{eq:Galerkin_3} can now be solved for different choices of the parameter ${\rm E_c}>0$ and blow-up function~$\mathscr{G}$. Figure \ref{fig:2} displays the approximations of the lowest energy band ${\bold k} \mapsto \varepsilon_{1,\bold{k}}$ for two different choices of $\mathscr{G}$ and the same ${\rm E_c}>0$ and effective potential $V \in L_{\rm per}^{\infty}(\Omega)$ as chosen to produce Figure~\ref{fig:01}. The most interesting feature of Figure~\ref{fig:2} is the fact that-- in contrast to the approximate energy band $\varepsilon^{\rm E_c}_{1}$-- the approximate energy band $\widetilde{\varepsilon}^{\rm E_c}_{1} $ remains $\mathbb{L}^*$-periodic, while also appearing to no longer be discontinuous.
	
	In the next section, we will present our two main results on the analysis of the modified discretization \eqref{eq:Galerkin_3}. Our first result is on the error analysis of this modified discretization where we prove that the approximate modified energy bands $\{\widetilde{\varepsilon}_n^{\rm E_c}\}_{n \in \mathbb{N}^*}$ converge with the expected asymptotic rate to the exact energy bands $\{\varepsilon_n\}_{n \in \mathbb{N}^*}$ in a point-wise sense when the cutoff energy $\rm{E_c}$ goes to infinity. Our second result is a precise characterization of regularity properties of each energy band $\widetilde{\varepsilon}^{\rm E_c}_{n}$ with respect to blow-up functions of different singularity orders.
	
	\section{Main Results on the Analysis of the Operator Modification Approach} \label{sec:5}

	Throughout this section, we will assume the setting of Section \ref{sec:4}. We recall in particular Definition \ref{def:g} of the blow-up function $\mathscr{G}$ as well as the modified Hamiltonian matrices~$\{\widetilde{\mH}_{\bold{k}}^{\mathscr{G}, \rm E_c}\}_{{\rm E_c}>0, \bold{k} \in \mathbb{R}^d}$ defined through Equation \eqref{eq:fiber_mod}. 
	
	Our first main result concerns the error analysis of the $\bold{k}$-dependent modified Galerkin discretization \eqref{eq:Galerkin_3} of the exact eigenvalue problem \eqref{eq:eigenvalue_problem}. 

	\begin{theorem}[Error estimate for approximate, modified energy bands]\label{prop:1}~
		
		Consider the settings of the exact eigenvalue problem \eqref{eq:eigenvalue_problem} and the modified discretization \eqref{eq:Galerkin_3} with a blow-up function $\mathscr{G}$ satisfying the conditions in Definition~\ref{def:g}. Assume that the effective potential $V \in L^{\infty}_{\rm per}(\Omega) \cap H^r_{\rm per}(\Omega)$ for some $r >1-\frac{d}{4}$. Let $n \in \mathbb{N}^*$, and for each $\bold{k} \in \mathbb{R}^d$, let the subspace $\mY_{{\rm n}}^{\bold{k}} \subset H^2_{\rm per}(\Omega)$ be defined as the span of the first $n$-eigenfunctions of the exact fiber $\mH_{\bold{k}}$, i.e.,
		\begin{align*}
			\mY_{{\rm n}}^{\bold{k}}  := \text{\rm span} \{u_{j, \bold{k}} \colon j \in \{1, \ldots, n \} \}.
		\end{align*}
		Then there exists ${\rm E_c^*}>0$ and a constant ${\rm C}>0$ such that for every ${\rm E_c} \geq {\rm E_c^*}$ and all $\bold{k} \in \mathbb{R}^d$ it holds that 
		\begin{equation}\label{eq:Error_2}
			0\leq \widetilde{\varepsilon}_{n,\bold{k}}^{\rm E_c}- \varepsilon_{n,\bold{k}} \leq  \left(\frac{{\rm C}}{\rm E_c}\right)^{r+1-\frac{d}{4}} \max_{\substack{\Phi \in \mY_{{\rm n}}^{\bold{k}}\\[0.2em] \Vert\Phi \Vert_{L^2_{\rm per}(\Omega)} =1}} \Vert \Phi\Vert^2_{H^{r+2}_{\rm per}(\Omega)} .
		\end{equation}
	\end{theorem}
	
	An immediate consequence of Theorem \ref{prop:1} is that the modified energy bands $\{\widetilde{\varepsilon}_{n}^{\rm E_c}\}_{n \in \mathbb{N}^*}$ converge at the same asymptotic rate as the unmodified energy bands $\{{\varepsilon}_{n}^{\rm E_c}\}_{n \in \mathbb{N}^*}$, with respect to ${\rm E_c}$, to the exact bands $\{{\varepsilon}_n\}_{n \in \mathbb{N}^*}$. Additionally, Theorem \ref{prop:1} informs us that the approximate energy bands $\{\widetilde{\varepsilon}_{n}^{\rm E_c}\}_{n\in M_{\rm E_c}^-}$ are bounded functions of $\mathbb{R}^d$. This latter fact will be of use in the proof of our second main result (see Section \ref{sec:7}).
	
	Next, we present our second main result, which concerns the regularity properties of these energy bands.

	\begin{theorem}[Regularity of approximate, modified energy bands]\label{prop:2}~
		
		Consider the setting of the $\bold{k}$-dependent modified discretization \eqref{eq:Galerkin_3} with a blow-up function $\mathscr{G}$ satisfying the conditions in Definition~\ref{def:g}. Let ${\rm E_c}>0$ be such that $M_{\rm E_c}^->0$, and let $n \in \{1, \ldots, M_{\rm E_c}^- \}$. If $\bold{k}_0 \in \mathbb{R}^d$ is such that
		\begin{align*}
			\widetilde{\varepsilon}^{\rm E_c}_{n,\bold{k}_0} \neq \widetilde{\varepsilon}^{\rm E_c}_{\widetilde{n},\bold{k}_0} \quad \forall \widetilde{n} \in \{1, \ldots, M_{\rm E_c}^-\} \text{ with } \widetilde{n} \neq n \quad\hspace{2mm} \big(\text{no band crossings at } (\bold{k}_0, \widetilde{\varepsilon}^{\rm E_c}_{n,\bold{k}_0})\big),
		\end{align*}
		then the approximate energy band $\widetilde{\varepsilon}^{\rm E_c}_{n}$ is of class $\mathscr{C}^m$ in a neighborhood of $\bold{k}_0$.
		
		If on the other hand, $\bold{k}_0 \in \mathbb{R}^d$ is such that
		\begin{align*}
			\exists \widetilde{n} \in \{1, \ldots, M_{\rm E_c}^-\} \text{ with } \widetilde{n}\neq n \colon \quad  \widetilde{\varepsilon}^{\rm E_c}_{n,\bold{k}_0} = \widetilde{\varepsilon}^{\rm E_c}_{\widetilde{n},\bold{k}_0} \quad \big(\text{band crossing at } (\bold{k}_0, \widetilde{\varepsilon}^{\rm E_c}_{n,\bold{k}_0})\big),
		\end{align*}
		then the approximate energy band $\widetilde{\varepsilon}^{\rm E_c}_{n}$ is 
		\begin{align}
			\begin{cases}
				\text{Lipschitz continuous in a neighborhood of $\bold{k}_0$} & \text{ if } m\geq 1,\\[1em] 
				\text{continuous in a neighborhood of $\bold{k}_0$} & \text{ otherwise }.
			\end{cases}
		\end{align}

	\end{theorem}
	Theorem \ref{prop:2} indicates that by designing a blow-up function $\mathscr{G}$ that satisfies Properties (1)-(4) from Definition \ref{def:g}, and in particular has a blow-up singularity of the correct order, we can obtain modified energy bands $\{\widetilde{\varepsilon}^{\rm E_c}_n\}_{n \in \mathbb{N}^*}$ of arbitrarily high regularity away from band crossings. Moreover, thanks to Theorem \ref{prop:1}, we also have a precise {\it a priori} characterization of the error in a given band $\widetilde{\varepsilon}^{\rm E_c}_{n}$ with respect to varying cutoff energies ${\rm E_c}$.
	
	 In order to prove the continuity result in Theorem \ref{prop:2}, we will make use of the following general lemma, which is also valid for non-Hermitian matrices and could be used to extend the present analysis to more advanced electronic structure models such as GW~(see \cite{cances_gontier_stoltz} for a mathematical analysis of the latter model). Note however that the continuity result in Theorem~\ref{prop:2} can also be obtained by using arguments specific to Hermitian matrices based on spectral inequalities and residual estimates. 
 
		\begin{lemma}\label{lem:contin}~Let $M \in \mathbb{N}^*$, let $p \in\mathbb{N}^*$ be such that $1 \le p < M$, and let $(\mH_{n})_{n \in \mathbb{N}} \in \left(\mathbb{C}^{M \times M}\right)^{\mathbb{N}}$ be a sequence of matrices that admit the block decomposition
			\begin{align*}
				\mH_{n}=  \arraycolsep=6pt\def\arraystretch{2} \left[
				\begin{array}{c|c}
					{\rm A}_{n} & {\rm B}_{n} \\
					\hline
					\widetilde{{\rm B}_{n}} & {\rm C}_{n}
				\end{array}
				\right],
			\end{align*}
			where $	{\rm A}_{n}  \in \mathbb{C}^{p \times p}$, ${\rm B}_{n}  \in \mathbb{C}^{p \times (M-p) }$, $\widetilde{{\rm B}_{n}}  \in \mathbb{C}^{(M-p) \times p}$, and ${\rm C}_{n}  \in \mathbb{C}^{(M-p) \times  (M-p)}$ are sub-matrices such that
			\begin{align*}
				&\exists  {\rm A} \in \mathbb{R}^{p \times p } \quad \text{\rm such that} \quad \lim_{n \to \infty} \Vert {\rm A}_{n} -{\rm A}\Vert_2 =0,\\[1em]
				&\sup_{n \in \mathbb{N}} \Vert {\rm B}_{n} \Vert_2 < \infty, \quad \sup_{n \in \mathbb{N}} \Vert \widetilde{\rm B}_{n} \Vert_2 <  \hspace{1mm}\infty, \\[1em]
    &{\rm C}_n \text{\rm \; is invertible for each $n$ and} \quad \lim_{n \to \infty} \Vert {\rm C}_{n}^{-1} \Vert_2 =0,
			\end{align*}
			with $\Vert \cdot \Vert_2$ denoting the usual matrix $2$-norm. Then 
			
			\begin{enumerate}
				
				\item for every $\rho >0$ sufficiently small, there exists $N(\rho)\in \mathbb{N}$ such that for any eigenvalue $\lambda^{\rm A}$ of the matrix~${\rm A}$ with algebraic multiplicity $Q \in \mathbb{N}^*$ and all $n \geq N(\rho)$, the open disc $\mathbb{B}_{\rho}\big(\lambda^{\rm A}\big) \subset \mathbb{C}$ contains exactly $Q$ eigenvalues of the matrix ${\mH}_{n}$ counting algebraic multiplicities;

				\item for every $\Upsilon>0$ sufficiently large, there exists $\widetilde{N}(\Upsilon)\in \mathbb{N}$ such that for all $n \geq \widetilde{N}(\Upsilon)$, there are exactly $M-p$ eigenvalues of $\mH_{n}$ with magnitude larger than or equal to $\Upsilon$.
			\end{enumerate}
		\end{lemma}
		
		Before stating the proofs of Theorems \ref{prop:1} and \ref{prop:2} and Lemma \ref{lem:contin}, we will present some numerical results on the use of the operator modification approach that we have described. The aim of these numerical studies is to provide numerical support for the conclusions of our main results Theorems \ref{prop:1} and \ref{prop:2}. These numerical studies are the subject of the next section.
	
	\section{Numerical Results} \label{sec:6}

	Throughout this section, we assume the setting described in Sections \ref{sec:2}-\ref{sec:5}. Our goal is now two-fold. First, we wish to present numerical results supporting the conclusions of Theorem \ref{prop:1} and Theorem \ref{prop:2}. Second, we would like to demonstrate the effectiveness of the operator modification methodology described in Section \ref{sec:4} for computing the energy bands of realistic materials such as graphene and face-centered cubic (FCC) silicon crystals.\vspace{2mm}

	\subsection{Validation of theoretical results in one spatial dimension}~
	
	We begin by considering a simple one-dimensional geometric setting. We set $\mathbb{L}= \mathbb{Z}$ which results in $\mathbb{L}^*= 2\pi \mathbb{Z}$, $\Omega = [-\frac{1}{2}, \frac{1}{2})$ and $\Omega^*= [-\pi, \pi)$. The effective potential $V \in L^{\infty}_{\rm per}(\Omega)$ is chosen so that $V \in H^{1-\epsilon}_{\rm per}(\Omega)$ for every $\epsilon >0$. Figure \ref{fig:01b} in Section \ref{sec:3} displays a plot of the chosen potential. For all subsequent simulations, the eigenvalue solver tolerance was set to machine (double) precision and the reference eigenvalues $\{\varepsilon_{n, k}\}_{n \in \mathbb{N}^*, k \in \mathbb{R}}$ were computed using the uniform Galerkin discretization \eqref{eq:Galerkin_1} with ${\rm E_c}= 72,000$. Unless stated otherwise, the $\bold{k}$-point mesh-width is chosen equal to $\Delta = 10^{-3}$. All blow-up functions $\mathscr{G}$ have regularity $\mathscr{C}^6$ on the interval $(0, 1)$ and are of the form ${\rm C}(1-x)^{-p}$ in the vicinity of $1^-$. \vspace{4mm}

	\noindent \textbf{Error convergence with respect to $\rm E_c$}
	
	Our first set of numerical experiments is designed to demonstrate the dependence of the eigenvalue errors in the operator modification approach as a function of the discretization parameter $\rm E_c$. We compute the lowest energy bands $\widetilde{\varepsilon}_{1}^{\rm E_c}$ and ${\varepsilon}_{1}^{\rm E_c}$ for different values of the cutoff energy ${\rm E_c}$. We also compute the Fermi levels $\widetilde{\mu}_{\rm F}^{\rm E_c}$ and $\mu_{\rm F}^{\rm E_c}$ corresponding to the energy bands $\widetilde{\varepsilon}_{1}^{\rm E_c}$ and ${\varepsilon}_{1}^{\rm E_c}$ and one electron per unit cell, and we consider the errors 
	\begin{align*}
		\int_{\Omega^*} \left\vert \left(\varepsilon_{n,\bold{k}}- \mu_{\rm F}\right) - \left(\widetilde{\varepsilon}_{n,\bold{k}}^{\rm E_c}- \widetilde{\mu}^{\rm E_c}_{\rm F}\right)\right\vert\, d\bold{k} \quad \text{and} \quad \int_{\Omega^*} \left\vert \left(\varepsilon_{n,\bold{k}}- \mu_{\rm F}\right) - \left({\varepsilon}_{n,{\bold k}}^{\rm E_c}- \mu^{\rm E_c}_F\right)\right\vert\, d\bold{k}.
	\end{align*}
	
	\begin{figure}[ht]
		\centering
		\begin{subfigure}{0.495\textwidth}
			\centering
			\includegraphics[width=\textwidth, trim={0cm, 0cm, 0cm, 0cm},clip=true]{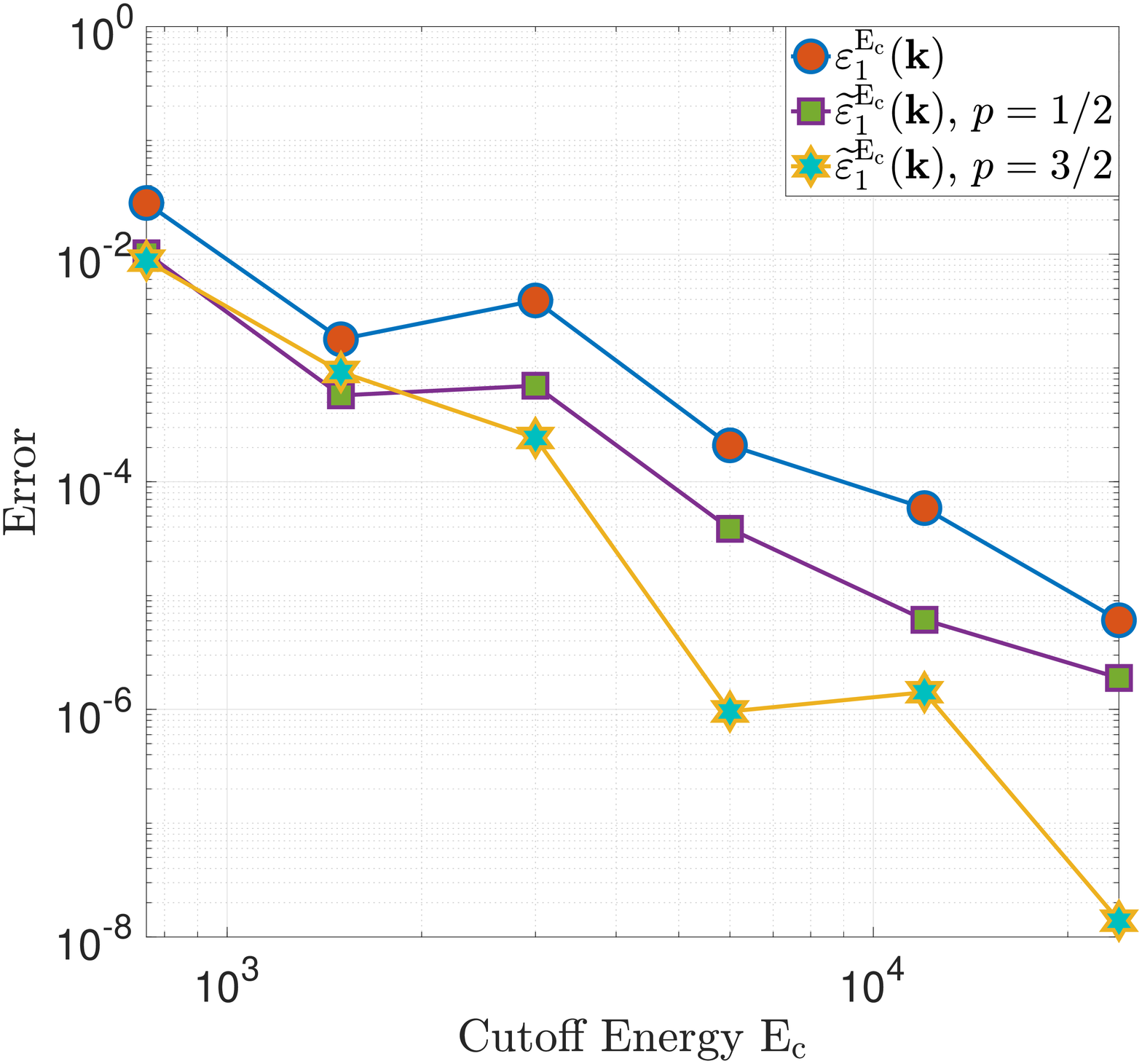} 
			\caption{Fermi level-adjusted error of the $\bold{k}$-dependent and modified energy bands as a function of ${\rm E_c}$.}		
			\label{fig:numerics_1a}
		\end{subfigure}\hfill
		\begin{subfigure}{0.495\textwidth}
			\centering
			\includegraphics[width=\textwidth, trim={0cm, 0cm, 0cm, 0cm},clip=true]{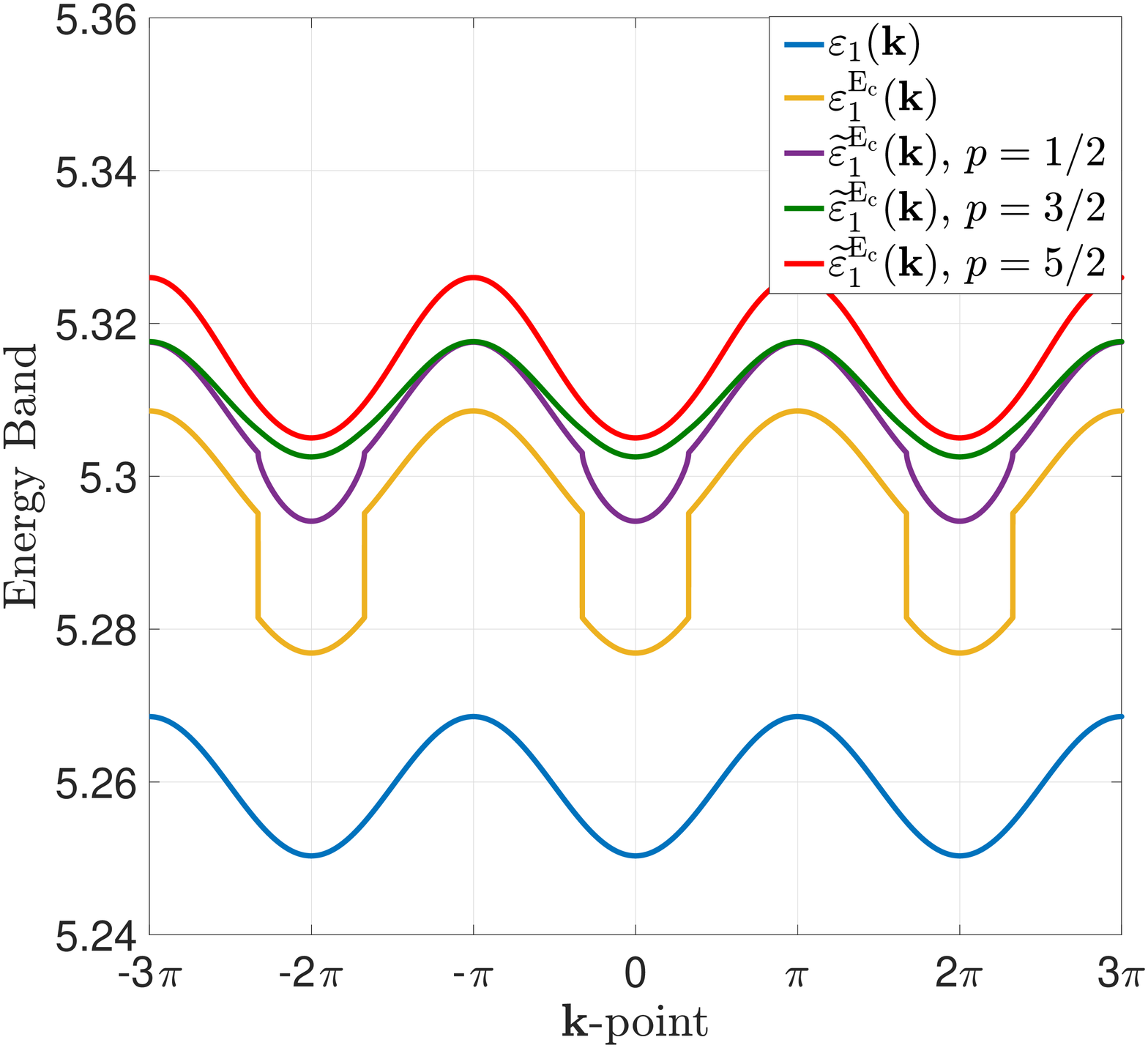} 
			\caption{Lowest modified energy bands for different choices of the blow-up function $\mathscr{G}$.}
			\label{fig:numerics_1b}
		\end{subfigure}
		\caption{Fermi level-adjusted error in lowest energy bands as a function of ${\rm E_c}$ (left) and the lowest energy bands for different blow-up functions of the form $\mathcal{h}(x)= {\rm C}(1-x)^{-p}$ in the vicinity of $1^-$. Notice that ${\varepsilon}_{1}^{\rm E_c}$ has a jump discontinuity, while the modified energy bands $\widetilde{\varepsilon}_{1}^{\rm E_c}$ are at least continuous.}
		\label{fig:numerics_1}
	\end{figure}

	\begin{figure}[ht]
		\centering
		\begin{subfigure}{0.495\textwidth}
			\centering
			\includegraphics[width=\textwidth, trim={0cm, 0cm, 0cm, 0cm},clip=true]{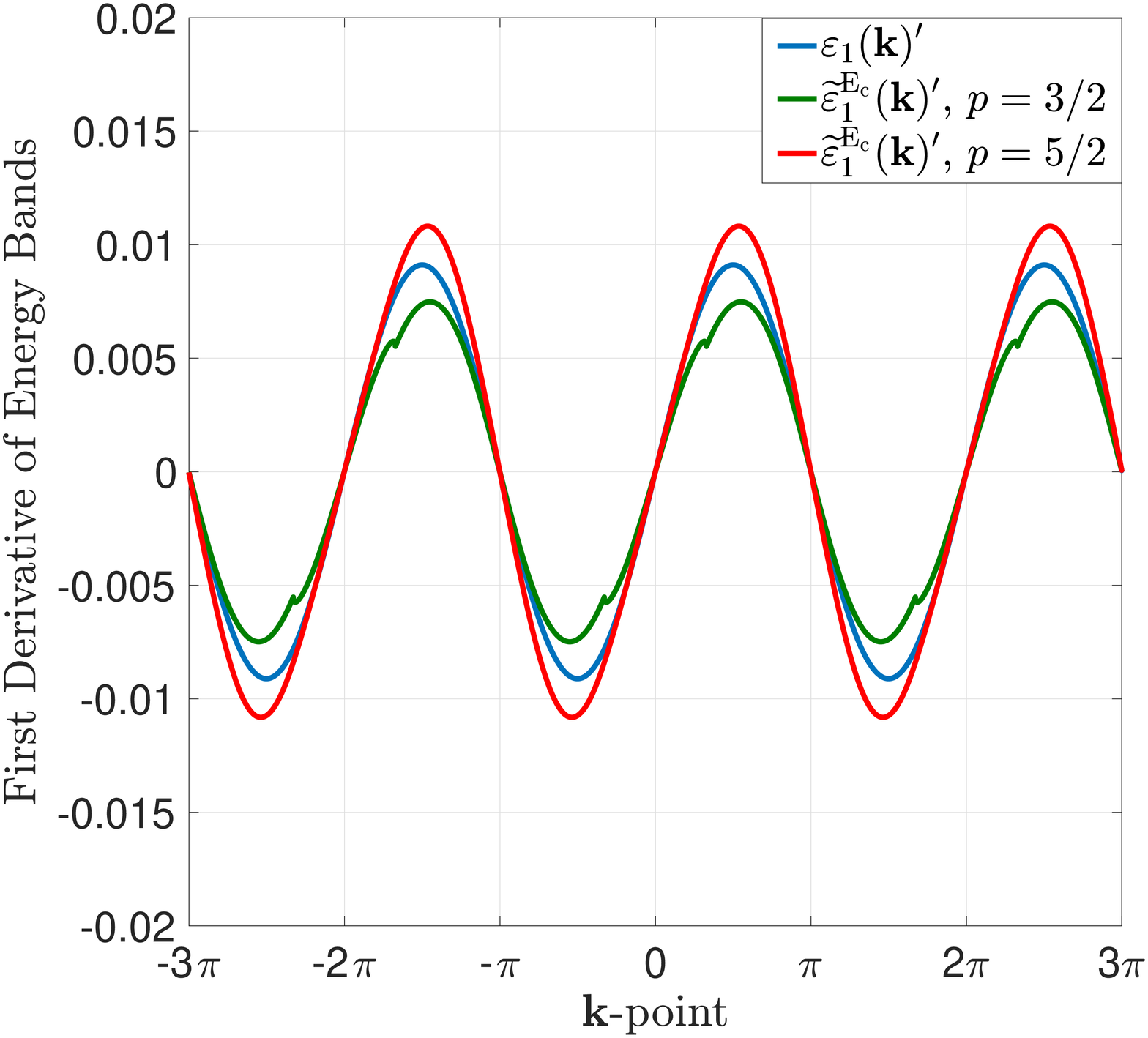} 
			\caption{First derivative of the lowest modified energy bands for different choices of the blow-up function $\mathscr{G}$.}		
			\label{fig:numerics_2a}
		\end{subfigure}\hfill
		\begin{subfigure}{0.495\textwidth}
			\centering
			\includegraphics[width=\textwidth, trim={0cm, 0cm, 0cm, 0cm},clip=true]{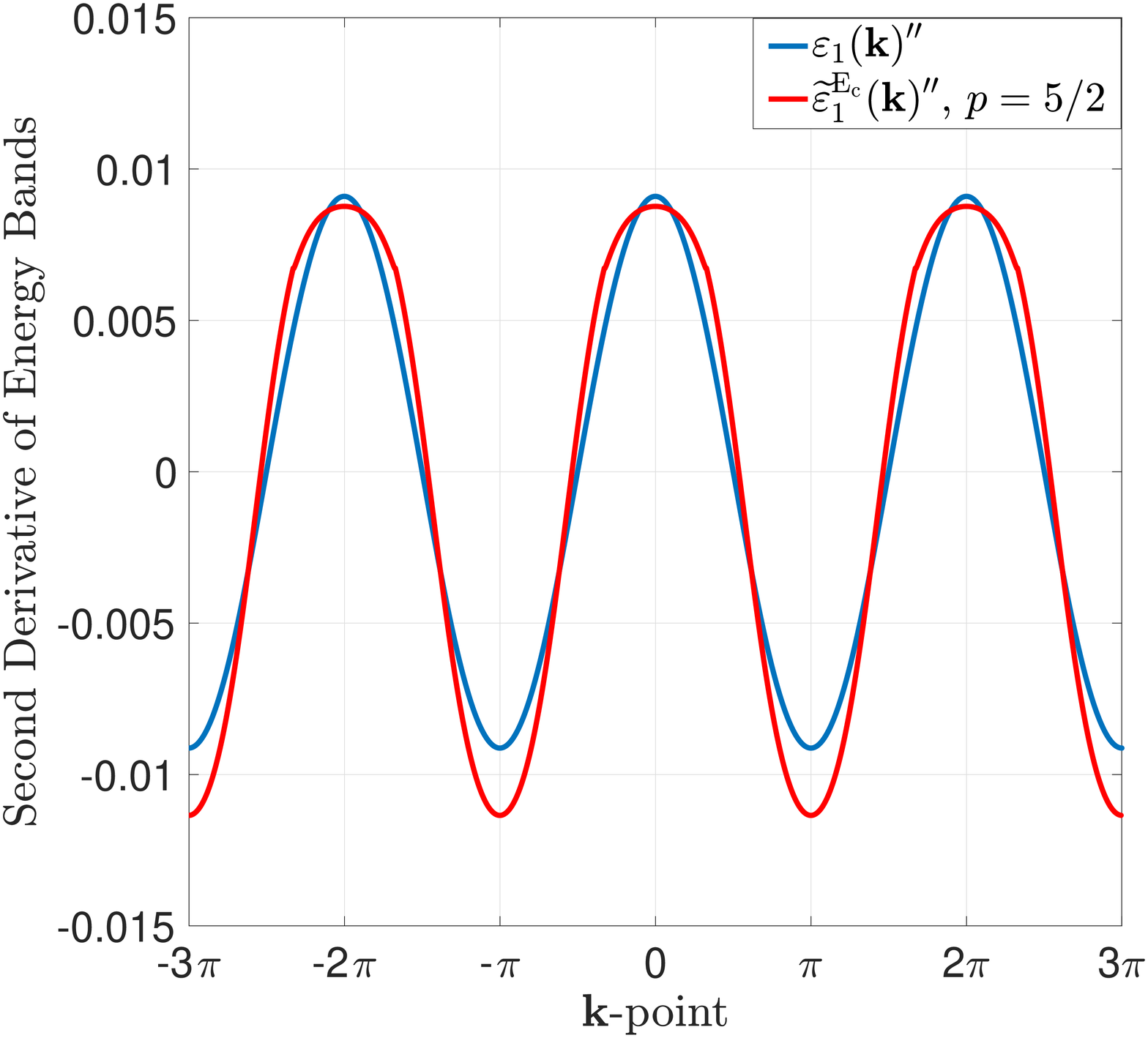} 
			\caption{Second derivative of the lowest modified energy bands for a higher order blow-up function $\mathscr{G}$.}
			\label{fig:numerics_2b}
		\end{subfigure}
		\caption{First and second derivatives of the lowest energy bands for different blow-up functions of the form $\mathcal{h}(x)= {\rm C}(1-x)^{-p}$ for $x \in [a, 1)$ with $a \in (\frac{1}{2}, 1)$ and ${\rm C}>0$. (Left) The energy band $\widetilde{\varepsilon}_{1}^{\rm E_c}$ produced using $p=\frac{3}{2}$ is of class~$\mathscr{C}^1(\mathbb{R})$ since it has a kink in the first derivative. (Right) The energy band $\widetilde{\varepsilon}_{1}^{\rm E_c}$ produced using  $p=\frac{5}{2}$ is of class~$\mathscr{C}^2(\mathbb{R})$ since it has a kink in the second derivative.}
		\label{fig:numerics_2}
	\end{figure}
	
		\begin{figure}[ht]
		\centering
		\begin{subfigure}{0.495\textwidth}
			\centering
			\includegraphics[width=\textwidth, trim={0cm, 0cm, 0cm, 0cm},clip=true]{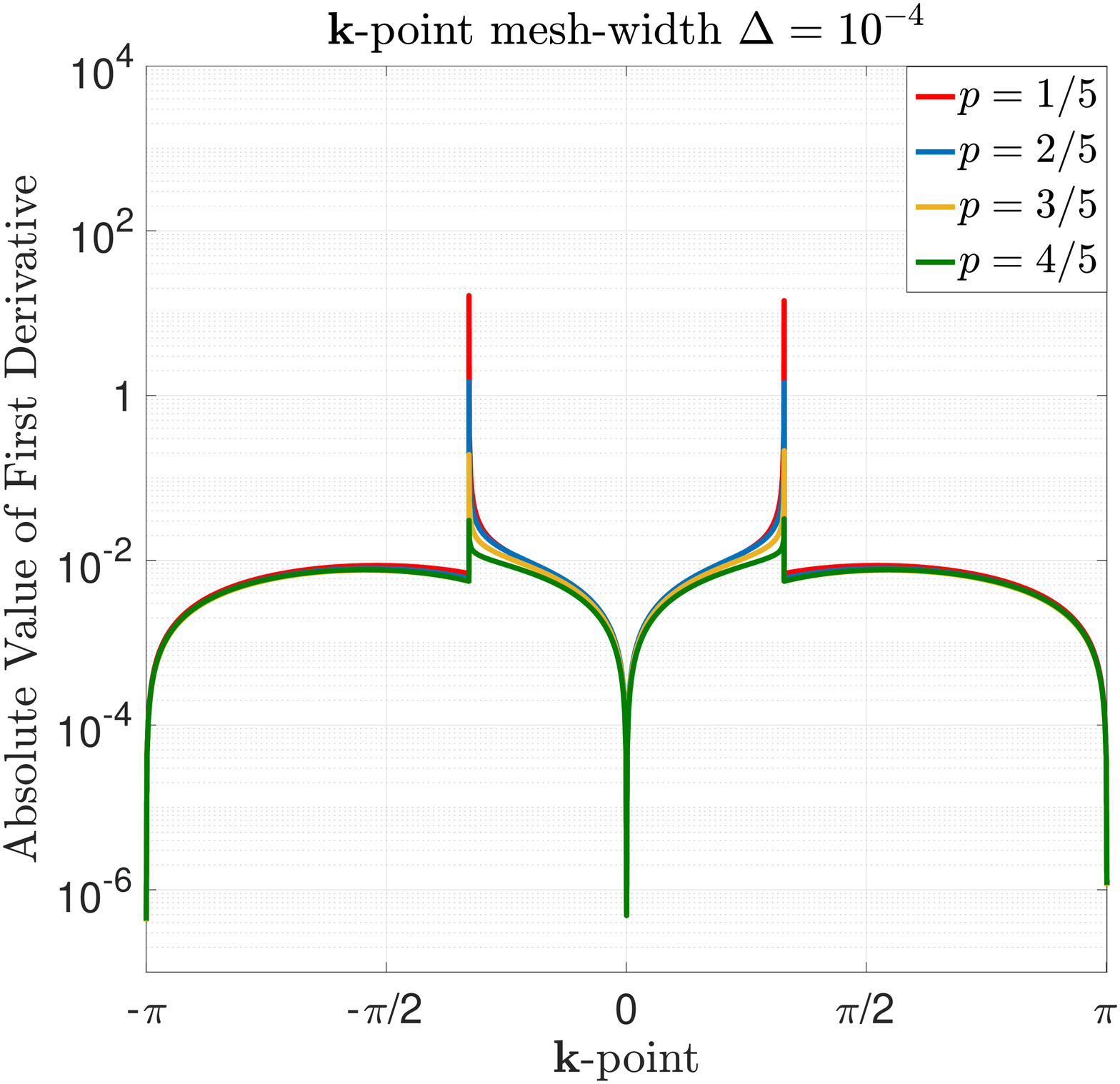} 
			\caption{First derivative of the lowest energy bands for different choices of blow-up function $\mathscr{G}$. }		
			\label{fig:numerics_3a}
		\end{subfigure}\hfill
		\begin{subfigure}{0.495\textwidth}
			\centering
			\includegraphics[width=\textwidth, trim={0cm, 0cm, 0cm, 0cm},clip=true]{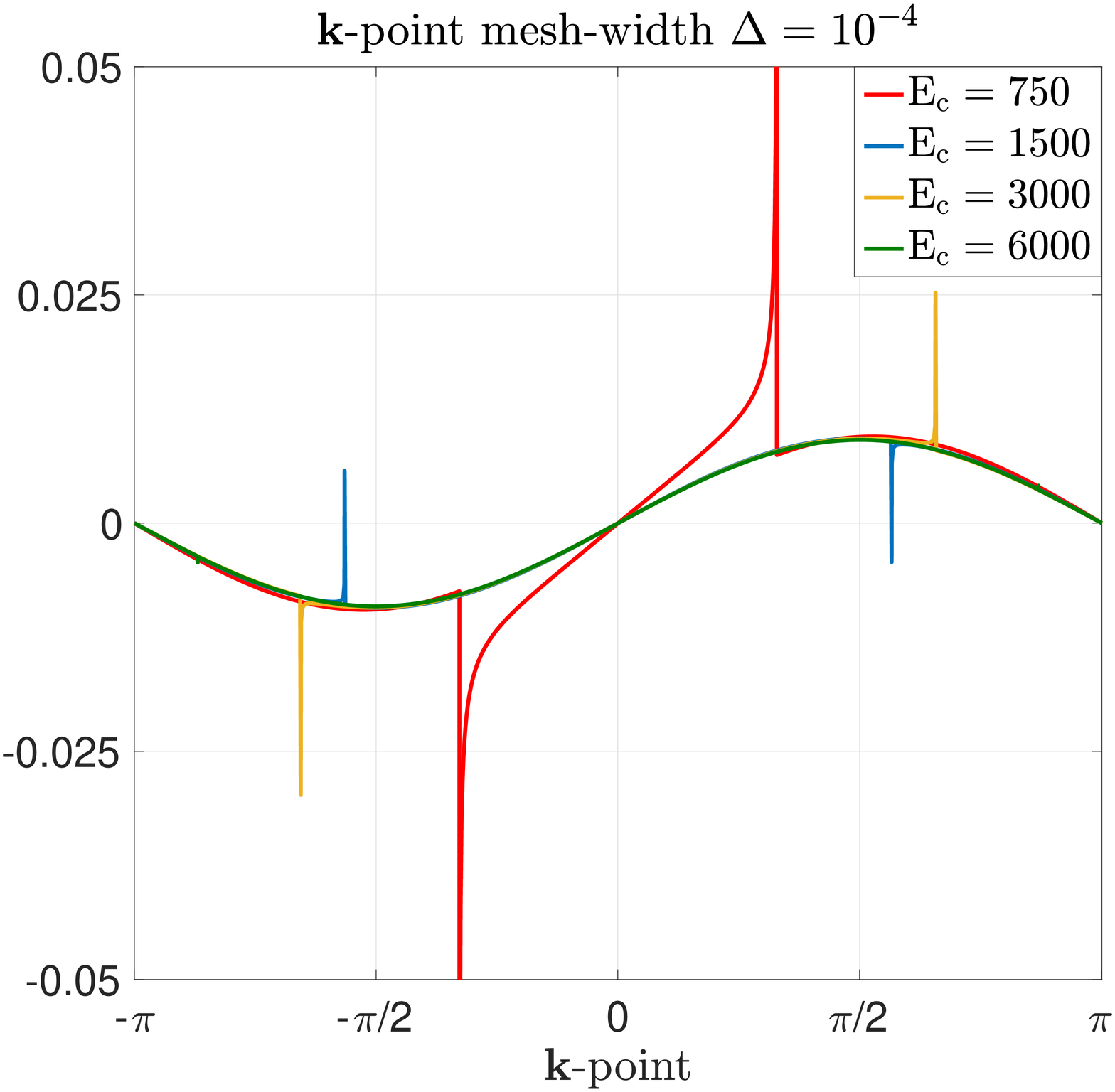} 
			\caption{First derivative of the lowest energy band for different choices of cutoff energies ${\rm E_c}$ with $\bold{k}$-point mesh-width $\Delta = 10^{-4}$.}
			\label{fig:numerics_3b}
		\end{subfigure}
		\caption{Non-Lipschitz energy bands ($m=0$): (Left) First derivative of the lowest energy band for different blow-up functions of the form $\mathcal{h}(x)= {\rm C}(1-x)^{-p}$ for $p <1$ in the vicinity of $1^-$. (Right) First derivative of the lowest energy band for different choices of cutoff energies ${\rm E_c}$ and with a blow-up function of the form $\mathcal{h}(x)= {\rm C}(1-x)^{-1/2}$ in the vicinity of $1^-$. Although not shown here, in both cases the magnitudes of the peaks increase when $\mathbf{k}$-point mesh-width is decreased indicating truly unbounded derivatives.}
		\label{fig:numerics_3}
	\end{figure}

	Figure \ref{fig:numerics_1a}  displays our results for two different choice of blow-up function $\mathscr{G} \colon \mathbb{R} \rightarrow \mathbb{R}$, one of which has a singularity blow-up of order $\vert \cdot \vert^{-\frac{1}{2}}$ and thus satisfies Properties (1)-(4) in Definition \ref{def:g} for $m=0$, and the other with a singularity blow-up of order $\vert \cdot \vert^{-\frac{3}{2}}$ which thus satisfies Properties (1)-(4) in Definition \ref{def:g} for $m=1$. We observe from Figure \ref{fig:numerics_1a} that the asymptotic convergence rate with respect to ${\rm E_c}$ of both the $\bold{k}$-dependent Galerkin discretization scheme \eqref{eq:Galerkin_2} and the modified discretization scheme \eqref{eq:Galerkin_3} are identical, and thus the use of the operator modification approach does not result in any asymptotic degradation of the discretization error. Additionally, we see that for a given cutoff energy ${\rm E_c}$, the error of the $\bold{k}$-dependent Galerkin discretization \eqref{eq:Galerkin_2}  is strictly larger than that of the modified discretization \eqref{eq:Galerkin_3}. \vspace{4mm}

	

	\noindent \textbf{Regularity of energy bands as a function of blow-up function singularity}
	
	Our next set of numerical simulations is designed to support the conclusion of Theorem \ref{prop:2} concerning the regularity of the energy bands produced by the modified Galerkin discretization \eqref{eq:Galerkin_3}. We consider the regularity of the lowest energy band $\widetilde{\varepsilon}_{1}^{\rm E_c}$ for cutoff energy ${\rm E_c} =750$ and three different choices of blow-up functions $\mathscr{G}\colon \mathbb{R} \rightarrow \mathbb{R}$. More precisely, we consider blow-up functions $\mathscr{G}$ that satisfy Properties (1)-(4) from Definition \ref{def:g} for $m = 0, 1,$ and $2$ respectively. Our theoretical results indicate that the resulting energy bands should be of class $\mathscr{C}^0(\mathbb{R})$, $\mathscr{C}^1(\mathbb{R})$ and $\mathscr{C}^2(\mathbb{R})$ respectively since there are no energy band crossings for non-trivial effective potentials in one dimension (see, e.g., ~\cite[Chapters 8-9]{Ashcroft76}). Figures \ref{fig:numerics_1b}, \ref{fig:numerics_2a}, and \ref{fig:numerics_2b} display our results and show perfect agreement with the conclusions of Theorem~\ref{prop:2}.

	Considering the energy bands displayed in Figure \ref{fig:numerics_1b}, it is natural to ask if the use of a blow-up function that satisfies Properties (1)-(4) from Definition \ref{def:g} only for $m = 0$ results in energy bands that are \emph{Lipschitz} continuous rather than simply continuous, and a similar question can be asked for the derivatives of the energy bands when using blow-up functions with stronger singularities. 
	
	In order to answer this question, we compute the lowest energy band $\widetilde{\varepsilon}_{1}^{\rm E_c}$ resulting from the modified discretization \eqref{eq:Galerkin_3} for cutoff energy ${\rm E_c} =750$ and four different choices of blow-up functions $\mathscr{G}$. The blow up functions $\mathscr{G}$ are constructed such that they satisfy Properties (1), (3) and (4) from Definition \ref{def:g} for $m= 6$, and such that the singularity of $\mathscr{G}(x)$ at $x=1$ is of order $\vert \cdot \vert^{-p}$ for $p= \frac{1}{5}, \frac{2}{5}, \frac{3}{5},$ and $\frac{4}{5}$ respectively. Figure \ref{fig:numerics_3a} displays the absolute values of the first derivatives of the resulting lowest energy band $\widetilde{\varepsilon}_{1}^{\rm E_c}$ for a $\bold{k}$-point mesh-width $\Delta = 10^{-4}$. The figure indicates that the derivative of $\widetilde{\varepsilon}_{1}^{\rm E_c}$ exhibits peaks at the two points of discontinuity, although the magnitude of the peaks at the points of discontinuity seems to decrease with increasing ${\rm E_c}$. In fact, although we do not display the plot here, the magnitudes of these peaks increase as the mesh width $\Delta$ is decreased, which indicates that the first derivative is truly unbounded at these points.  \vspace{4mm}
	
	

	\noindent \textbf{Regularity of energy bands as a function of the cutoff energy}
	
	The goal of the final set of numerical experiments in this subsection is to explore, for a fixed choice of blow-up function $\mathscr{G} \colon \mathbb{R} \rightarrow \mathbb{R}$, how the regularity of the modified energy bands $\{\widetilde{\varepsilon}_{n}^{\rm E_c}\}_{n \in \mathbb{N}^*}$ varies as a function of the cutoff energy ${\rm E_c}$. To this end, we consider once again the regularity of the lowest energy band $\widetilde{\varepsilon}_1^{\rm E_c}$ resulting from the modified Galerkin discretization \eqref{eq:Galerkin_3}. For these experiments, we set the $\bold{k}$-point mesh-width $\Delta = 10^{-4}$, and use a blow-up function $\mathscr{G}$ with blow-up of order $\vert \cdot \vert^{-\frac{1}{2}}$.


	Figure \ref{fig:numerics_3b} displays the first derivative of the energy band $\widetilde{\varepsilon}_1^{\rm E_c}$ for different values of ${\rm E_c}$. The blow-up function chosen for these simulations satisfies Properties (1)-(4) from Definition \ref{def:g} only for $m = 0$, so Theorem \ref{prop:2} indicates that the energy bands should be continuous and not differentiable. It is readily seen that this is indeed the case, although the first derivative of $\widetilde{\varepsilon}_1^{\rm E_c}$ is noticeably regularized by increasing the value of~${\rm E_c}$. In fact, the finite magnitude of the peak is a numerical artifact since (although not shown here) the magnitude of the peaks increases as the $\bold{k}$-point mesh-width is decreased, which indicates that the derivative is truly unbounded for finite ${\rm E_c}$ at these points. Note that the points of discontinuity occur at different $\bold{k}$-points depending on the chosen energy cutoff~${\rm E_c}$.


	
	\vspace{2mm}
	
	\subsection{Numerical experiments on real materials}~

	We now investigate the effectiveness of the operator modification approach introduced in Section \ref{sec:4} on two realistic crystalline systems: a single layer of graphene and FCC crystalline silicon. In what follows, all computations are performed using the plane-wave density functional theory package {\it DFTK.jl} \cite{DFTKjcon} in the {\it Julia} language~\cite{Julia-2017}. The code used to produce the following results has been integrated into the {\it DFTK.jl} package which is available at \cite{DFTK}.

\vspace{5mm}
	\begin{remark}\label{rem:graphene_and_co} 
		Two comments are in order:
  \begin{itemize}
      \item DFTK makes use of norm-conserving pseudopotentials consisting of a local component (a periodic multiplication operator) and a non-local component (a periodized finite-rank operator). The Kohn-Sham Hamiltonians obtained with DFTK are therefore not Schr\"odinger operators of the form \eqref{eq:Hamiltonian} with $V$ a periodic function. Our theoretical results can easily be extended to handle general pseudopotentials at the price of slightly more cumbersome proofs. For the sake of simplicity, we chose to restrict our analysis to the case of purely local potentials;
      \item graphene is a real material living in the three-dimensional physical space, but its Bravais lattice is the two-dimensional hexagonal lattice, hence the name 2D materials used to refer to graphene, hexagonal boron nitride, transition metal dichalcogenides, phosphorene, and other atomic-thin layered materials. The Bloch fibers of a periodic 2D material are labelled by a 2D quasi-momentum ${\mathbf k} \in \R^2$ and read as
      $$
      \mH_\bold{k} = \frac 12 (- i \nabla_{{\mathbf x}_\parallel}+ {\mathbf k})^2 - \frac 12 \partial^2_{x_3} + V
      $$
      where ${\mathbf x}_\parallel = (x_1,x_2)$ denotes the in-plane position, and $x_3$ the out-of-plane coordinate. The operators $\mH_\bold{k}$ act on $L^2_{\rm per}(\Omega)$ where $\Omega = \omega \times \R$ with $\omega \subset \R^2$ being the Wigner-Seitz cell of the 2D Bravais lattice. They do not have compact resolvents and do not admit spectral decompositions of the form \eqref{eq:eigenvalue_problem}. However, the Bloch fibers of the Kohn-Sham Hamiltonian of a real 2D material have discrete eigenvalues below the bottom of their essential spectrum forming the so-called valence bands and low-energy conduction bands. Our theoretical results can thus easily be extended to the case of 2D materials.
  \end{itemize}
	\end{remark}
	\vspace{1mm}

	\noindent {\bf Numerical setting}

\begin{figure}[ht]
		\centering
		\begin{subfigure}{0.45\textwidth}
			\includegraphics[width=\linewidth]{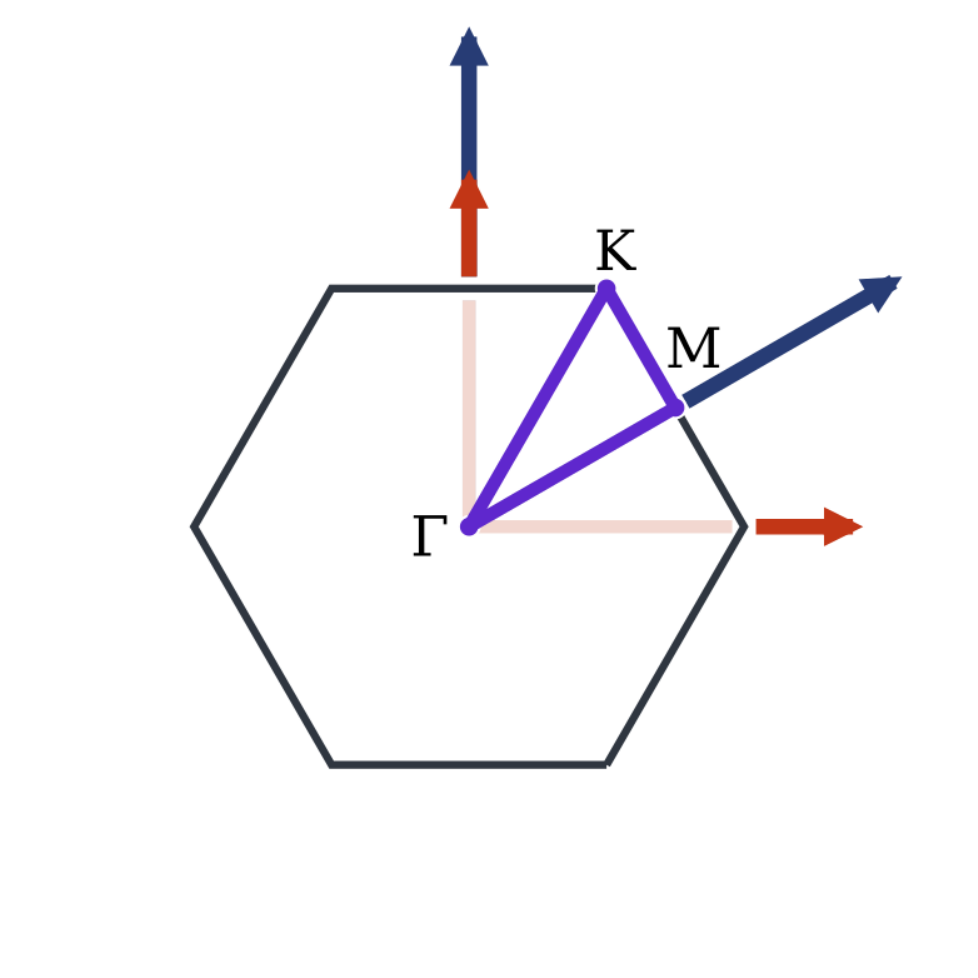}
			\caption{Graphene $\bf k$-path in the Brillouin zone.}
		\end{subfigure}
		\begin{subfigure}{0.45\textwidth}
			\includegraphics[width=\linewidth]{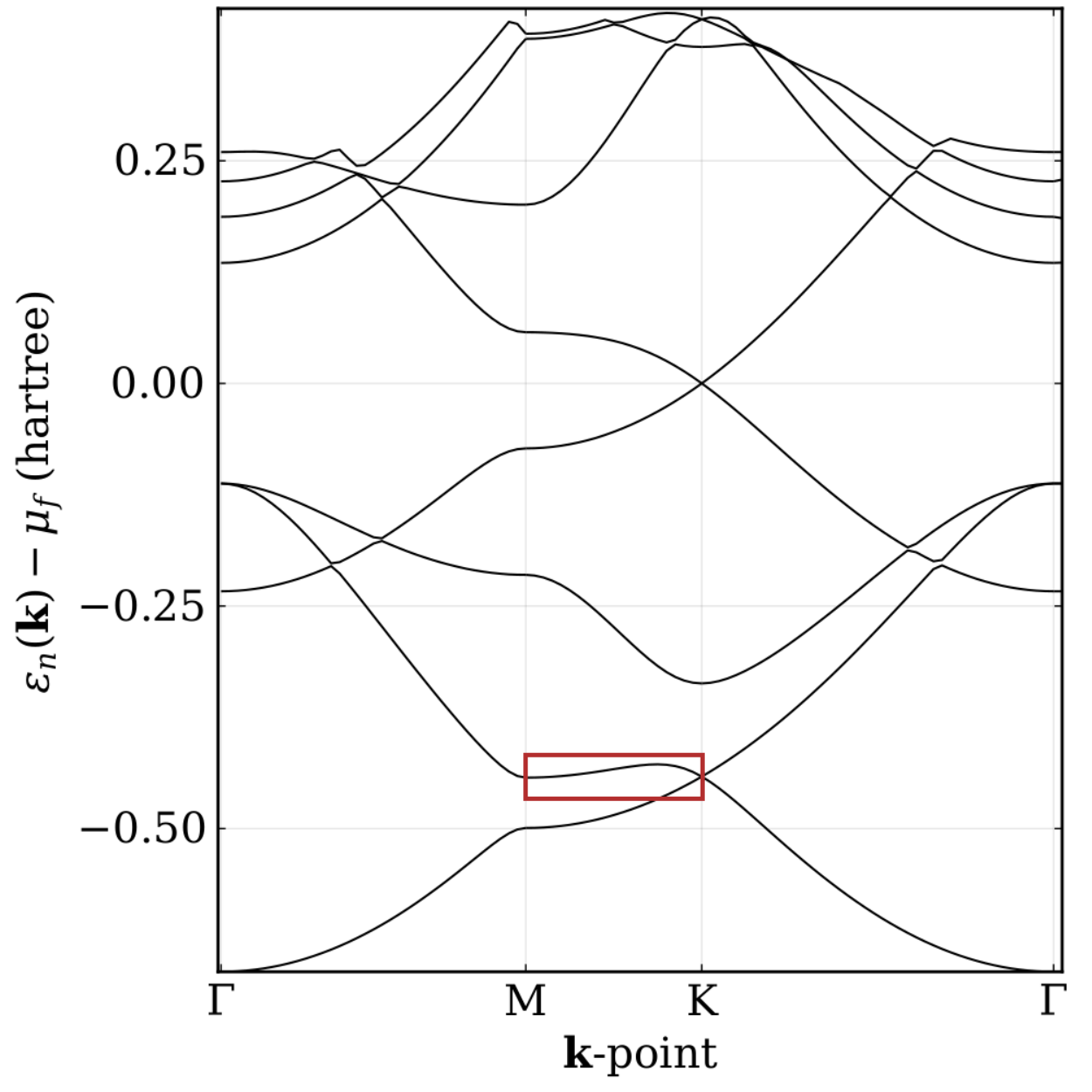}
			\caption{Graphene reference band structure along this path.}
		\end{subfigure}\hfill
		\caption{Graphene $\bf k$-path (left) and the corresponding band structure (right). The paths are automatically generated by the {\it Brillouin.jl} package based on crystallographic considerations. The red (resp. blue) arrows display the Cartesian coordinate axes (resp. the reciprocal basis vectors). All bands are shifted so that the Fermi level appears at~zero Hartree on the graph. The red box shows the part of the band diagram on which Figure \ref{fig:graphene_band_and_derivatives} focuses.}
		\label{fig:kpath_and_band_diagram_1}
	\end{figure}

\begin{figure}[ht]
		\begin{subfigure}{0.45\textwidth}
			\includegraphics[width=\linewidth]{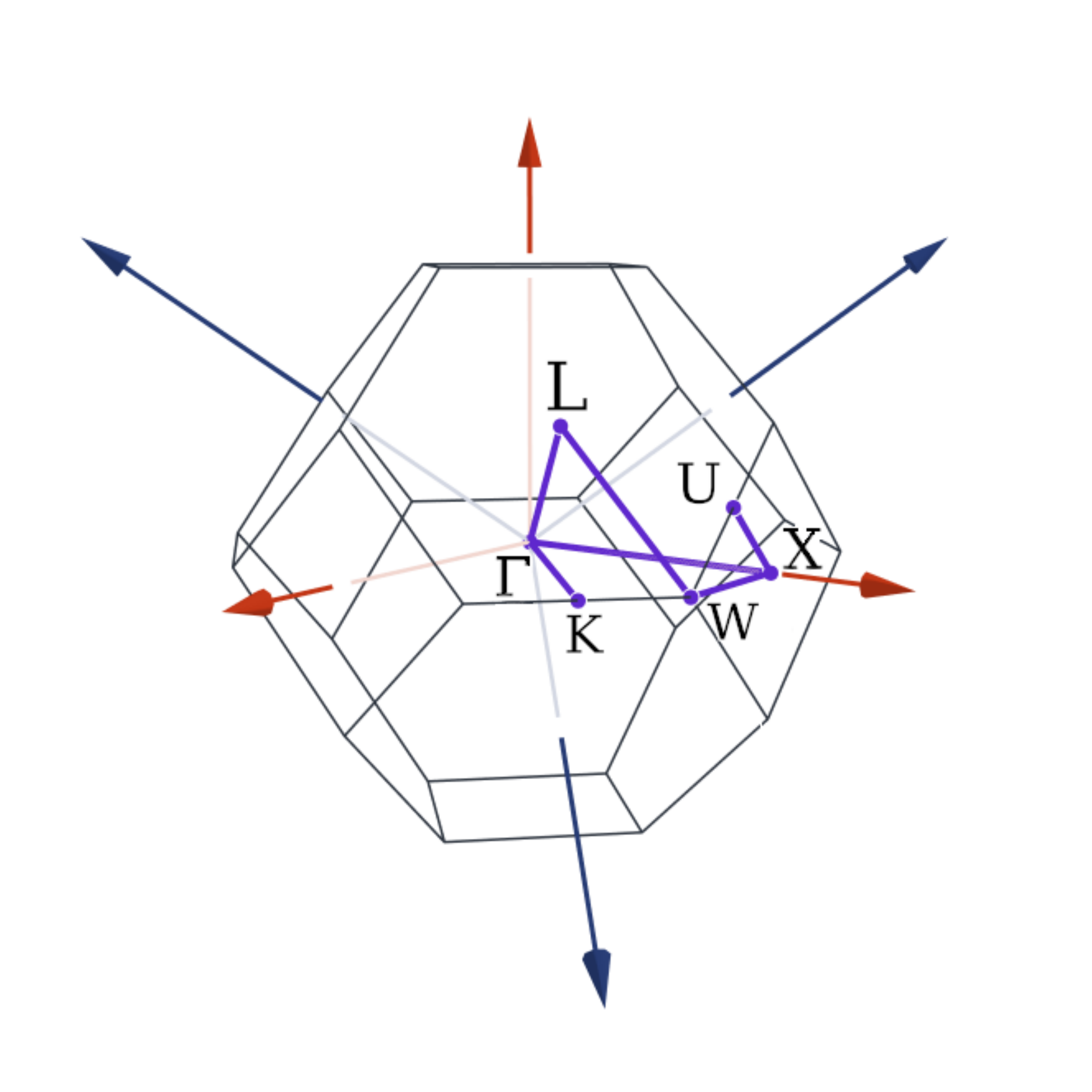}
			\caption{Silicon $\bf k$-path in the Brillouin zone.}
		\end{subfigure}
		\begin{subfigure}{0.45\textwidth}
			\includegraphics[width=\linewidth]{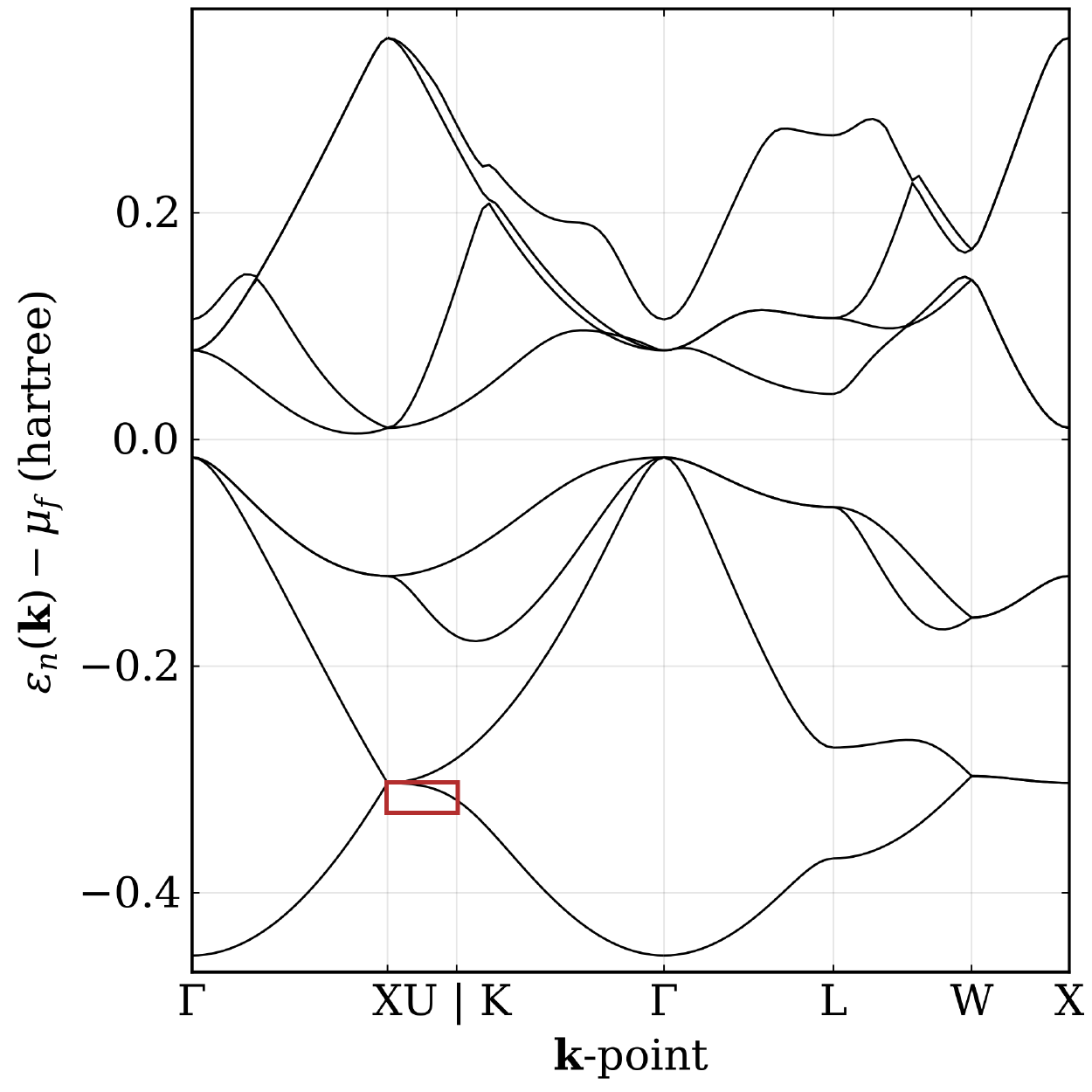}
			\caption{Silicon reference band structure along this path.}
		\end{subfigure}\hfill
		\caption{FCC crystalline silicon $\bf k$-path (left) and the corresponding band structure (right). The paths are automatically generated by the {\it Brillouin.jl} package based on crystallographic considerations. Red and blue arrows display the Cartesian coordinate axes and the reciprocal basis vectors respectively. All bands are again shifted so that the Fermi level appears at zero Hartree on the graph, and the red box shows the part of the band diagram on which Figure \ref{fig:silicon_band_and_derivatives} focuses.}
		\label{fig:kpath_and_band_diagram_2}
	\end{figure}

	We begin by computing a reference ground-state effective potential using a Kohn-Sham DFT self-consistent field procedure with cutoff energy ${\rm E_c^{\rm ref}}=30$ Ha and a fine Monkhorst-Pack $\bf k$--point grid with $12$ points in each sampled dimension. The SCF tolerance is set to machine (double) precision. We choose to work with a PBE functional \cite{perdew1996generalized} that is standard in solid state electronic structure computations and Hartwigsen-Goedecker-Teter-Hutter separable dual-space Gaussian pseudopotentials \cite{hartwigsen1998relativistic}. In a second step,  the obtained effective potential is used to construct the reference Hamiltonian $\mH_{\bf k}^{\rm E_c^{\rm ref}}$ whose eigenvalues will serve as reference data. 
	The same potential is used to construct the Hamiltonian operator $\mH_{\bf k}^{\rm E_c}$ defined through \eqref{eq:Galerkin_2}, and the modified Hamiltonian operator $\widetilde{\mH}_{\bf k}^{\rm E_c}$ defined through Equation \eqref{eq:Galerkin_3} for a custom set of $\bf k$--points and for a fixed ${\rm E_c} \ll {\rm E_c^{\rm ref}}$. The chosen $\bf k$-points are located on the band-structure paths automatically generated by the $Brillouin.jl$ package using the crystallography based method introduced in \cite{hinuma2017band}. For reference, the $\bf k$-paths in the Brillouin zone of graphene and FCC silicon and the corresponding band structures are displayed in Figures \ref{fig:kpath_and_band_diagram_1} and \ref{fig:kpath_and_band_diagram_2} respectively.
	
	\noindent {\bf Regularity of energy bands as a function of blow-up function singularity}
	
	\begin{figure}[ht]
		\centering
		\begin{subfigure}{0.49\textwidth}
			\includegraphics[width=0.9\linewidth]{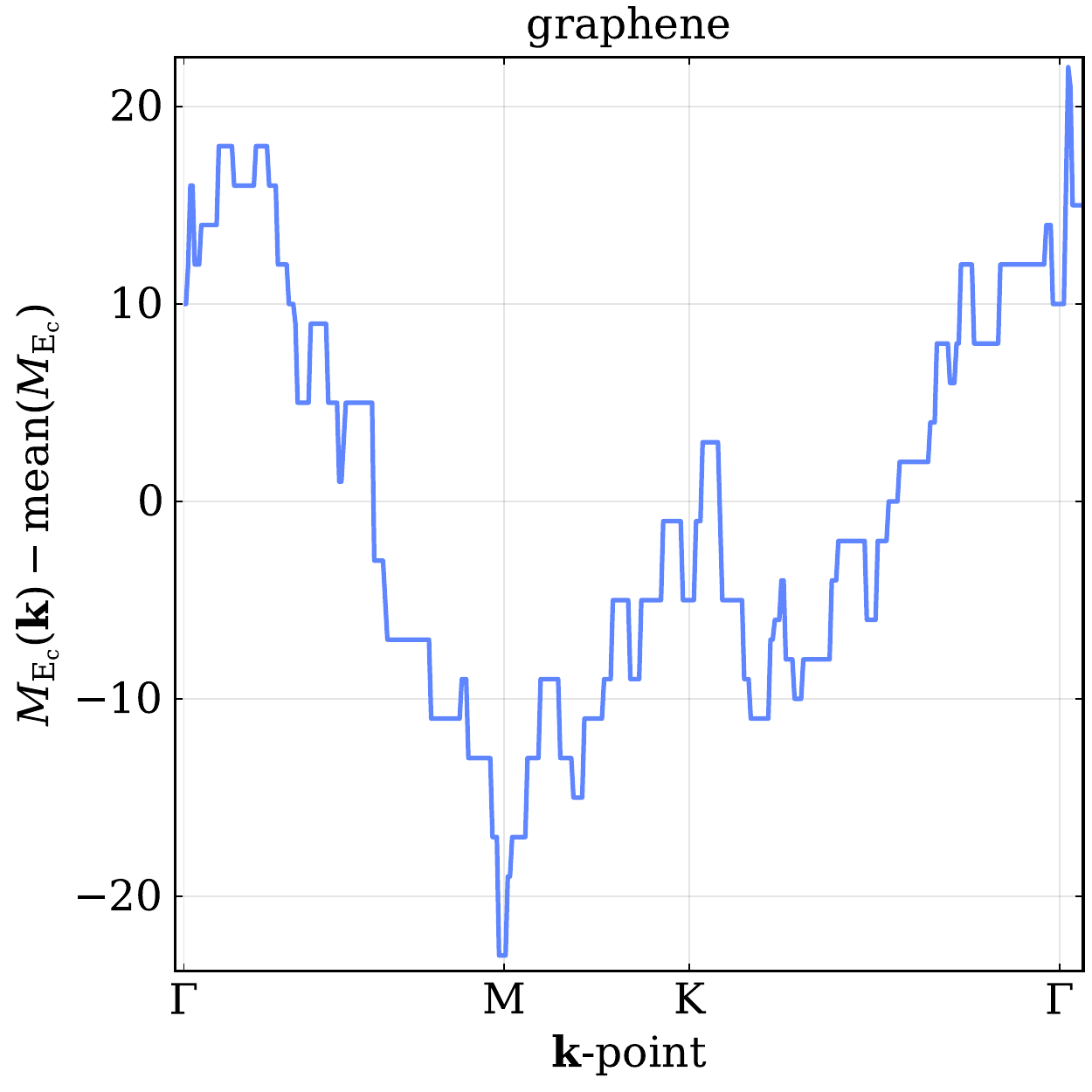}
		\end{subfigure}\hfill
		\begin{subfigure}{0.49\textwidth}
			\includegraphics[width=0.9\linewidth]{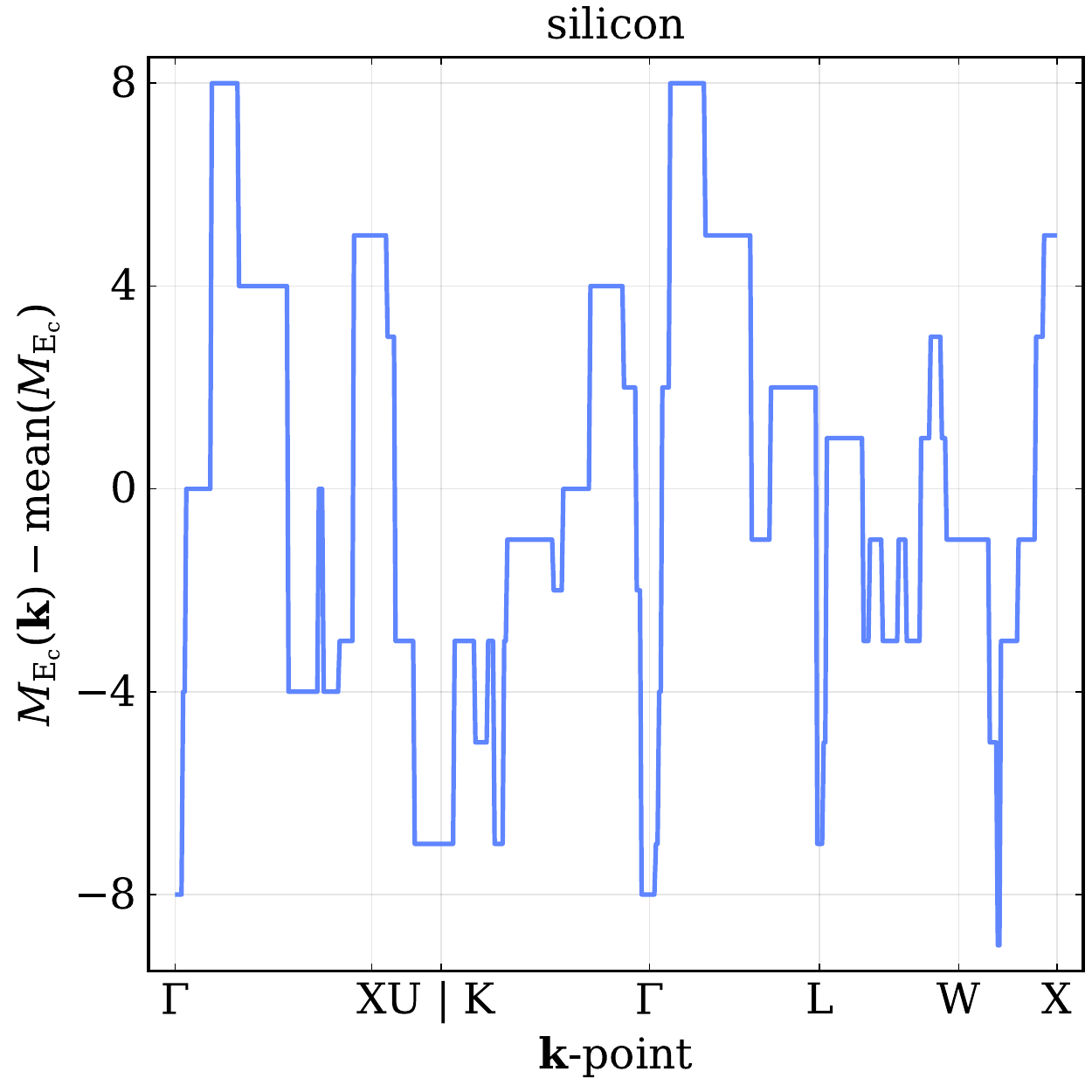}
		\end{subfigure}
		\caption{Distance to the mean value of $M_{\rm E_c}(\cdot)$ along standard $\mathbf{k}$-path with a small ${\rm E_c}=5$ Ha, showing the abrupt changes in the cardinality of the $\bf k$-dependent plane-wave basis set $\mathcal{B}_{\bold{k}}^{\rm E_c}$. The test cases are (left) a single layer of graphene and (right) FCC crystalline silicon. The $\bold{k}$-path is automatically generated by \textit{DFTK} using the \textit{Brillouin.jl} package.}
		\label{fig:M_EC}
	\end{figure}
	
\begin{figure}[ht]
		\centering
		\begin{subfigure}{0.33\textwidth}
			\includegraphics[width=\linewidth]{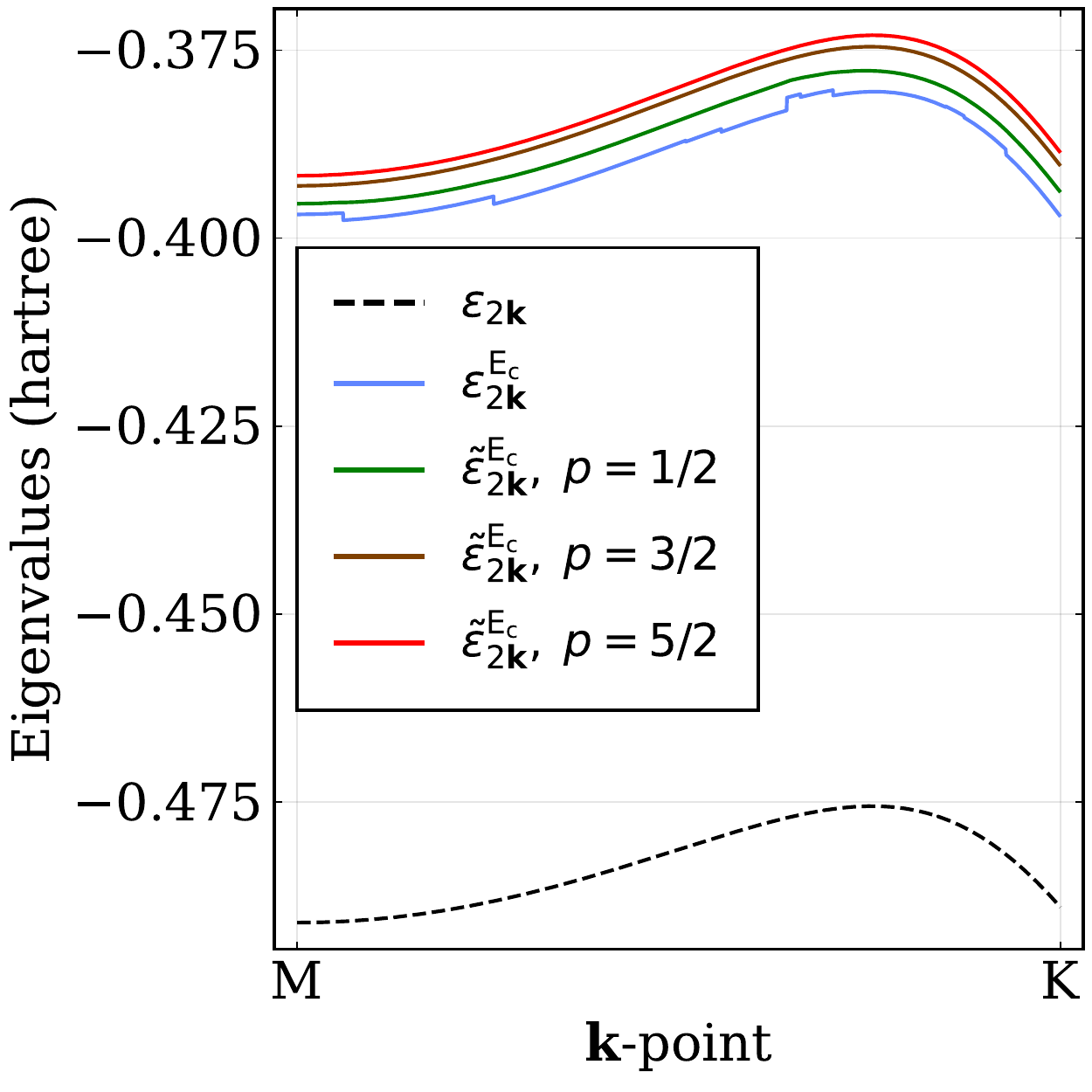}
			\caption{Energy band}
			\label{fig:graphene_first_band}
		\end{subfigure}\hfill
		\begin{subfigure}{0.33\textwidth}
			\includegraphics[width=\linewidth]{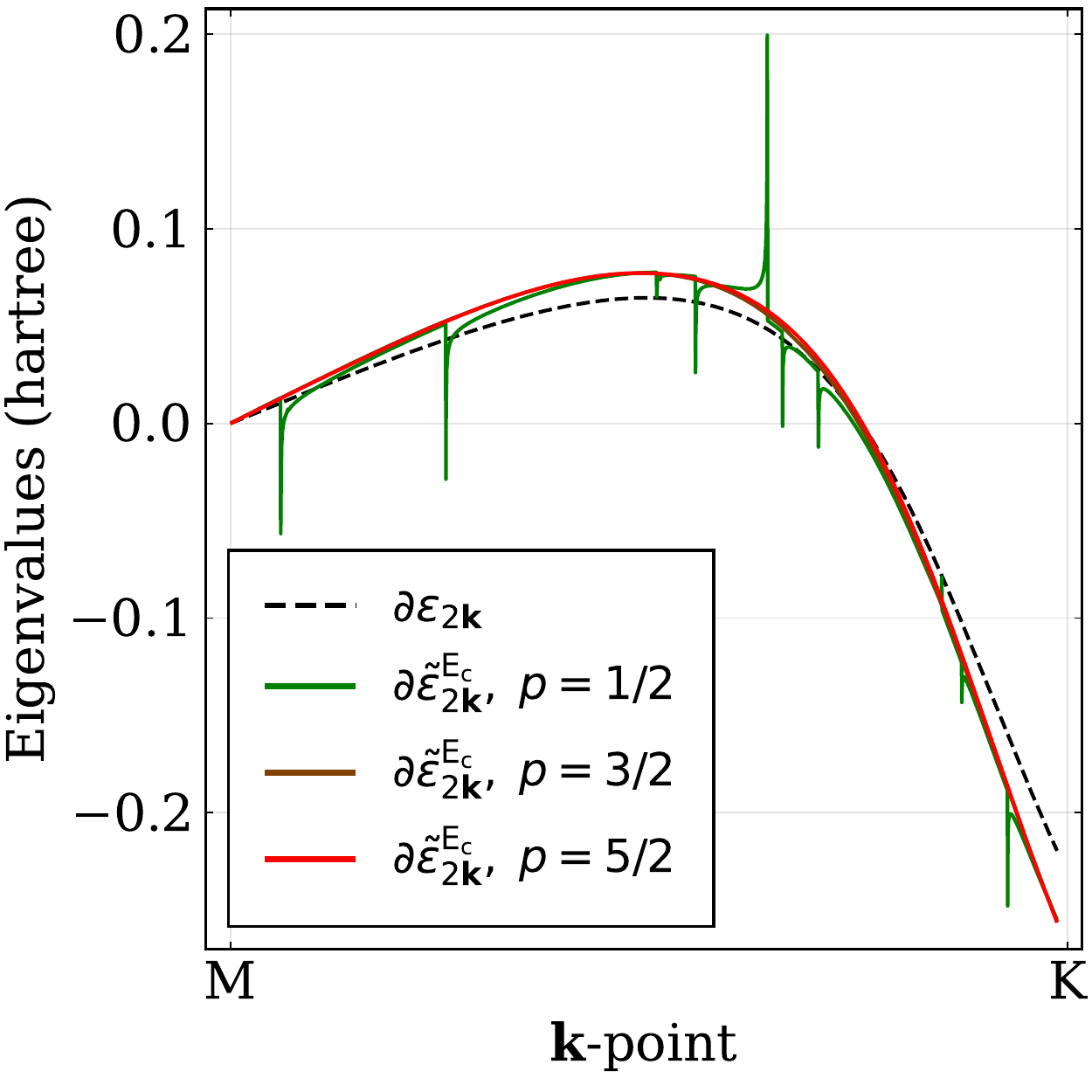}
			\caption{First derivative}
		\end{subfigure}\hfill
		\begin{subfigure}{0.33\textwidth}
			\includegraphics[width=\linewidth]{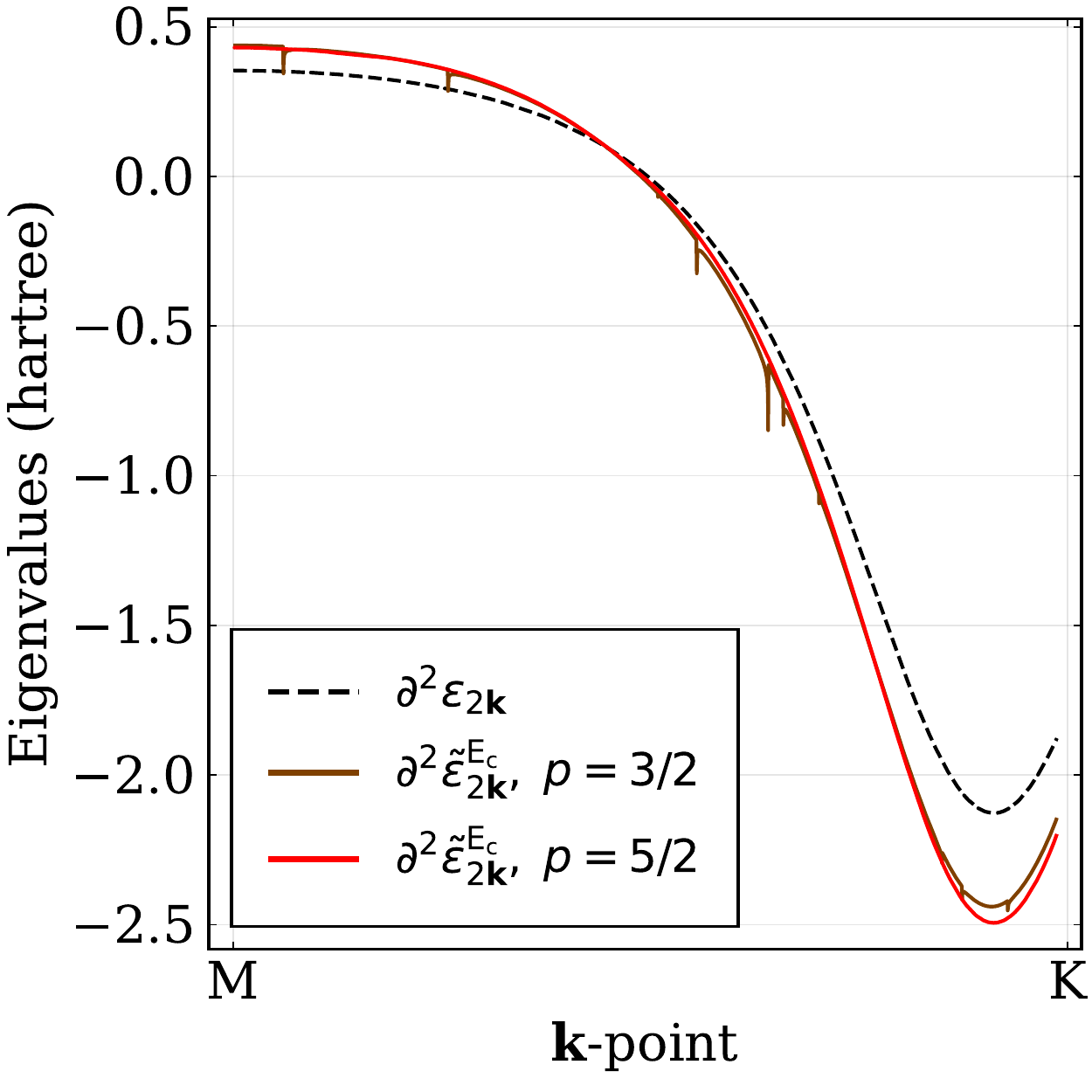}
			\caption{Second derivative}
		\end{subfigure}
		\caption{Comparison of the first and second derivatives of the second band of graphene between points {\bf M} and {\bf K} of the band-structure for the $\bold{k}$-dependent and modified discretization schemes.}
		\label{fig:graphene_band_and_derivatives}
	\end{figure}

\begin{figure}[ht]
		\centering
		\begin{subfigure}{0.33\textwidth}
			\includegraphics[width=\linewidth]{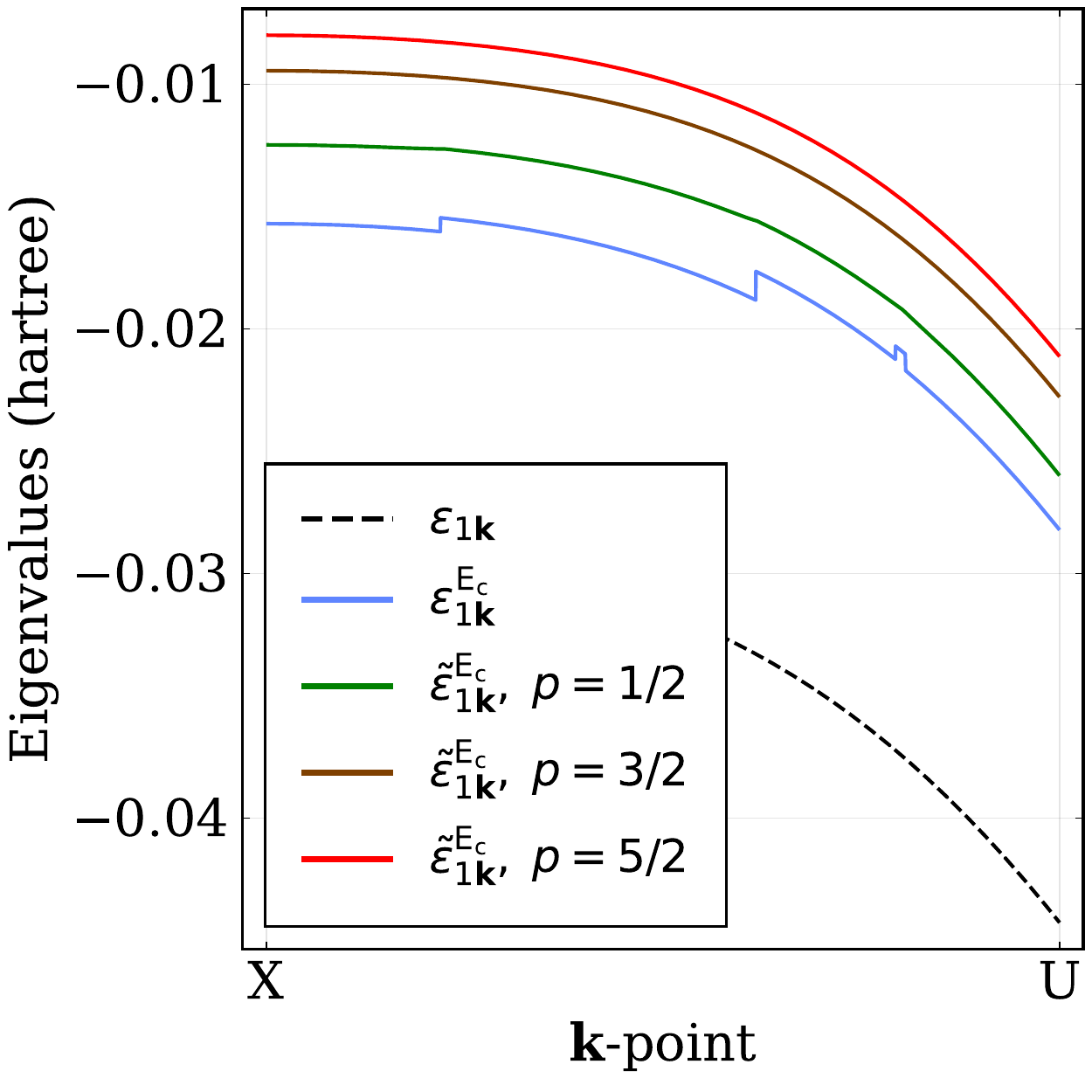}
			\caption{Energy band}
			\label{fig:silicon_first_band}
		\end{subfigure}\hfill
		\begin{subfigure}{0.33\textwidth}
			\includegraphics[width=\linewidth]{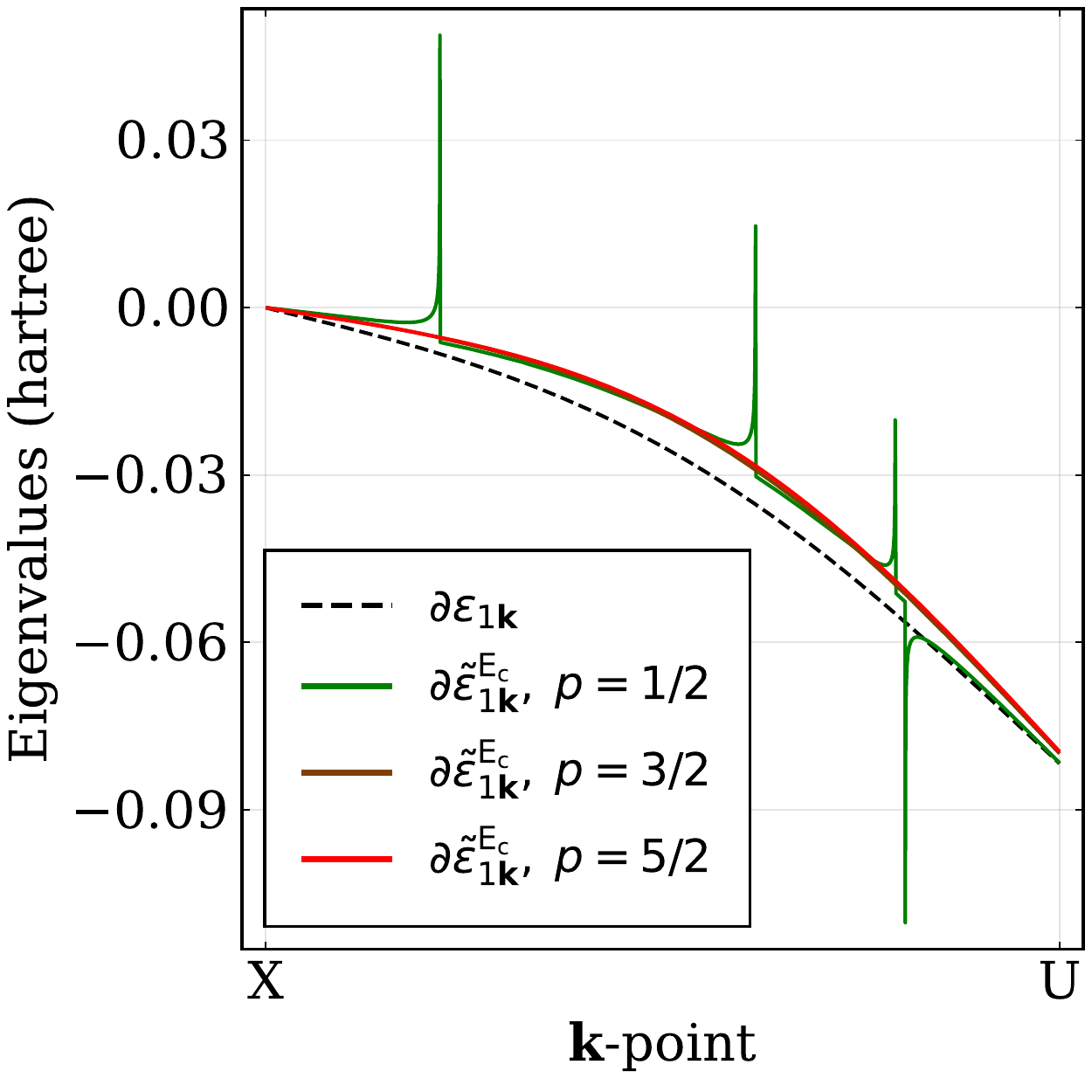}
			\caption{First derivative}
		\end{subfigure}\hfill
		\begin{subfigure}{0.33\textwidth}
			\includegraphics[width=\linewidth]{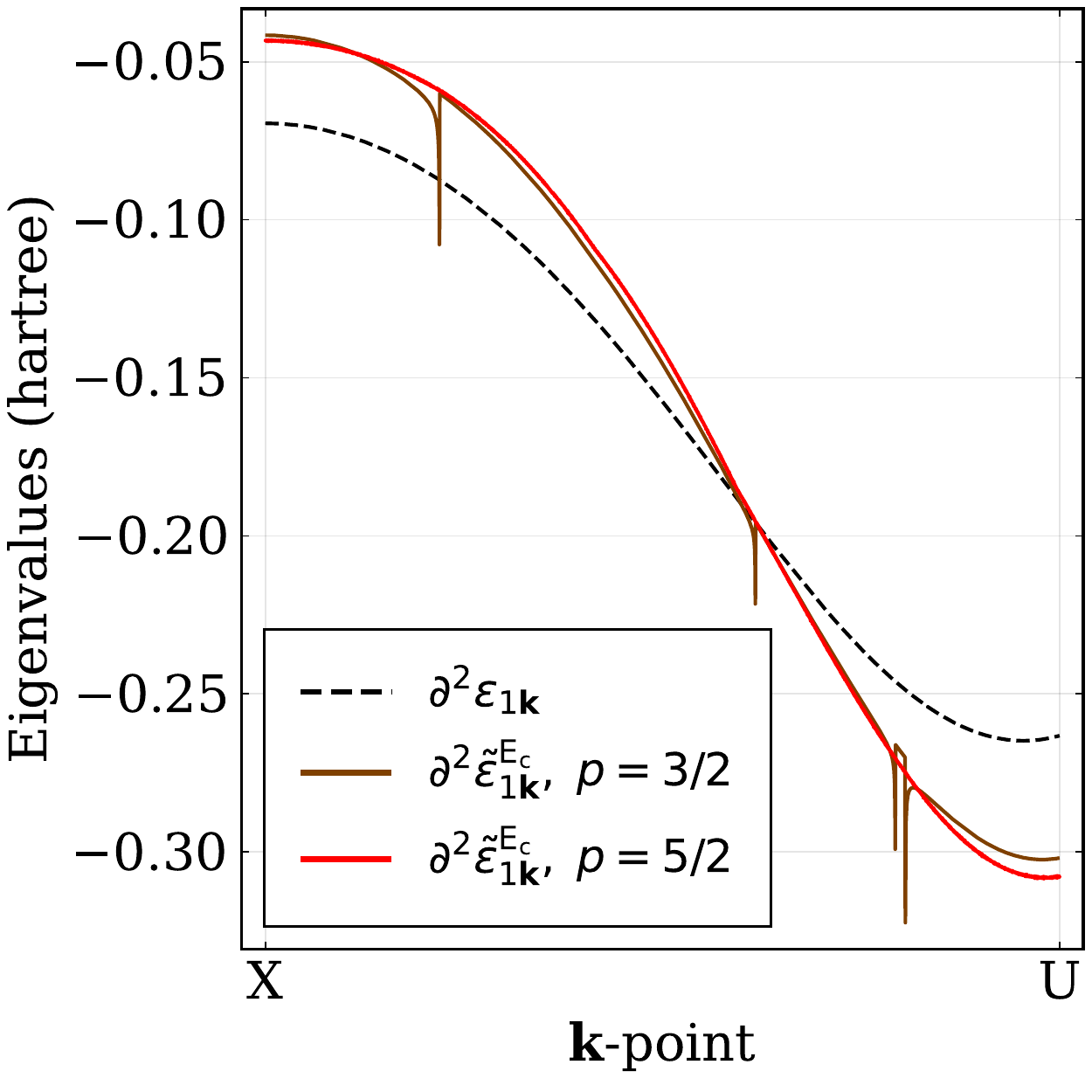}
			\caption{Second derivative}
		\end{subfigure}
		\caption{Comparison of the first and second derivatives of the first band of silicon between points {\bf X} and {\bf U} of the band-structure for the $\bold{k}$-dependent and modified discretization schemes.}
		\label{fig:silicon_band_and_derivatives}
	\end{figure}

	We choose a very low ${\rm E_c}=5$ Ha in order to clearly highlight the expected irregularities of the energy bands of the standard Hamiltonian operator $\mH^{\rm E_c}_{\bold{k}}$. Figure \ref{fig:M_EC} displays the abrupt changes in the size of the $\bf k$-dependent plane-wave basis $\mathcal{B}_{\bold{k}}^{\rm E_c}$ along the band-structure paths of graphene and FCC silicon for this choice of ${\rm E_c}$. As in the 1D case, the energy bands produced by the $\bold{k}$-dependent Galerkin discretization \eqref{eq:Galerkin_2} are highly irregular, as we read from Figures \ref{fig:graphene_first_band} and \ref{fig:silicon_first_band}. On the other hand, the modified energy bands produced by the Galerkin discretization \eqref{eq:Galerkin_3} (also displayed in Figures \ref{fig:graphene_band_and_derivatives} and \ref{fig:silicon_band_and_derivatives}) appear to be smooth in accordance with the choice of blow-up function $\mathscr{G}$ and in agreement with the theoretical result of Theorem \ref{prop:2}.

		\begin{figure}[h!]
		\centering
		\begin{subfigure}{0.49\textwidth}
			\includegraphics[width=0.9\linewidth]{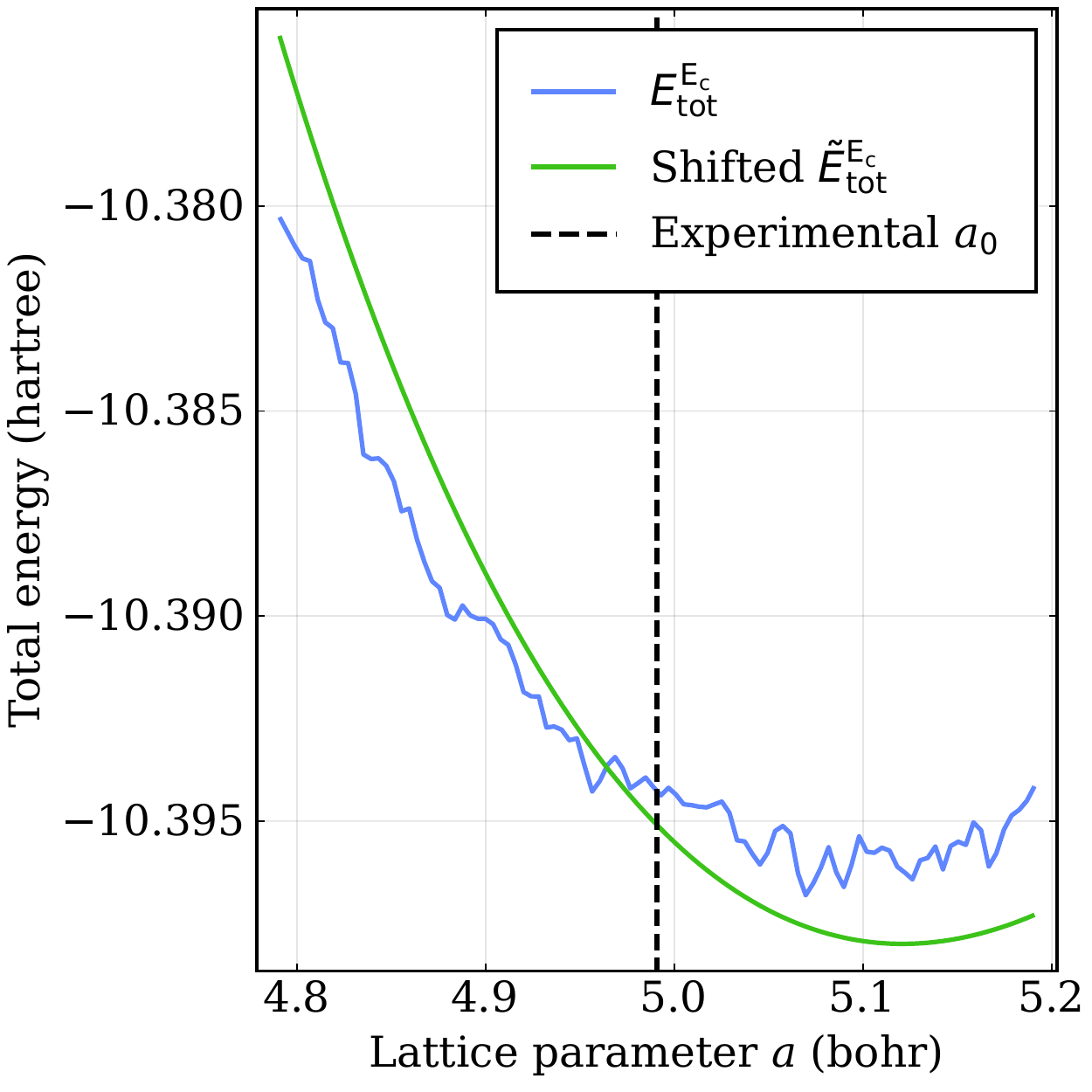}
			\caption{Graphene}
		\end{subfigure}\hfill
		\begin{subfigure}{0.49\textwidth}
			\includegraphics[width=0.9\linewidth]{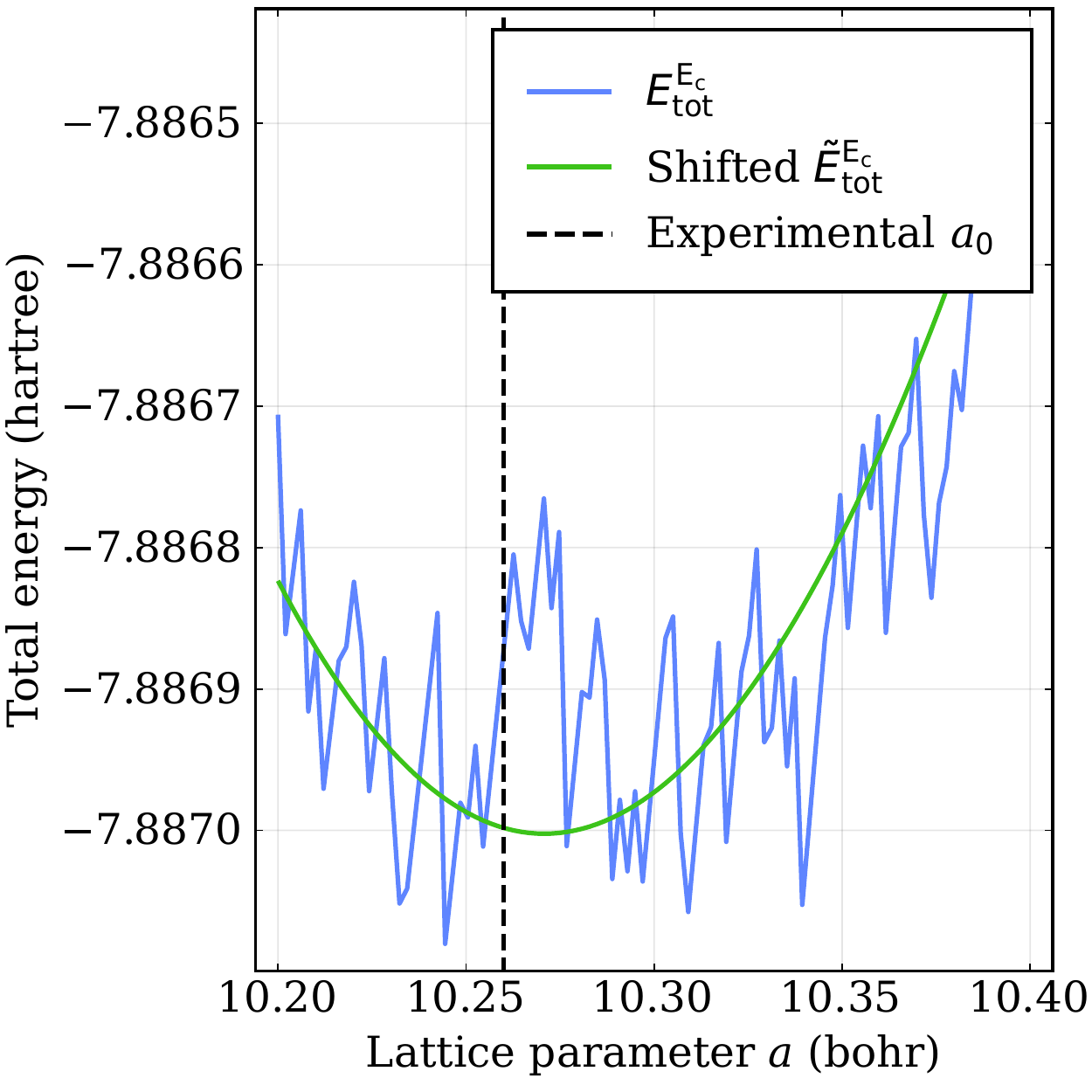}
			\caption{FCC crystalline silicon}
		\end{subfigure}
		\caption{Energy per unit volume of graphene and FCC silicon as a function of the lattice parameter $a$ for the $\bold{k}$-dependent and modified discretization schemes (Equations \eqref{eq:Galerkin_2} and \eqref{eq:Galerkin_3} respectively). The blow-up function $\mathscr{G}$ has a blow-up of order $|\cdot|^{-\frac{3}{2}}$ and the cutoff energy is set to $E_c=5$ Ha. For the sake of legibility, the total energy of the system for the modified Hamiltonian is shifted to the mean value of the standard Hamiltonian total energy over the sample of parameters $a$. The empirical value $a_0$ of the equilibrium lattice parameter is also indicated.} \label{fig:band_vs_lattice_cste}
	\end{figure}

	Let us remark that for small ${\rm E_c}$'s, the low-energy eigenfunctions $u_{n,{\bold k}}$ have non-negligible projections on Fourier modes $e_{\bold{G}}\in \mathcal{B}^{\rm E_{\rm c}}_{\bf k}$ with `modified' kinetic energy that is much higher than the standard kinetic energy, this artificially higher kinetic energy being due to the fact that the blow-up function $\mathscr{G}$ diverges from the $x \mapsto x^2$ curve when $x$ goes to $1$. As a result, the modified-Hamiltonian energy bands appear significantly over-estimated in comparison to the approximate energy bands produced by the standard Hamiltonian operator $\mH^{\rm E_c}_{\bold{k}}$. On the other hand, for large ${\rm E_c}$'s, the $e_{\bold{G}}$'s in $\mathcal{B}^{\rm E_{\rm c}}_{\bf k}$ with over-estimated kinetic energies do not contribute much to the Fourier expansions of the low-energy eigenfunctions $u_{n,{\bold k}}$. This suggests that the proper range of application of the modified-operator approach is typically in the regime where ${\rm E_c}$ is not too small so that the modified discretization scheme matches the accuracy of the standard Galerkin discretization while ensuring the targeted regularity. \vspace{4mm}
	
	\noindent {\bf Regularity of energy bands as a function of the volume of the unit cell}
	
	Another possible application of the modified-operator approach concerns the computations of energy bands as a function of the volume of the unit cell, that can be used to estimate the macroscopic volumetric mass density of a crystalline material or its bulk modulus~\cite[Chapter 5.6]{kaxiras_2003}. The size of the $\bf k$-dependent discretization basis at a given ${\bold k}$-point depends on the unit cell volume through its associated reciprocal lattice, so that within the standard Galerkin approximation, the energy per unit volume is a rough function of the crystalline parameters. For both graphene and FCC silicon, the unit cell is parameterized by a single lattice parameter $a$. Figure \ref{fig:band_vs_lattice_cste} displays the energy of graphene and FCC silicon per unit volume as a function of $a$ around the experimental value $a_0$ of the equilibrium lattice parameter. In both cases, we observe high oscillations of the energy per unit volume. We see that the modified-operator approach produces much smoother energy curves.

	\section{Proofs of the Main Results} \label{sec:7}

	We will begin with the proof of the {\it a priori} error estimate Theorem \ref{prop:1}.  \vspace{4mm}
	
	\begin{proof}[\textbf{Proof of Theorem \ref{prop:1}}]~
		
		The lower bound in Inequality \eqref{eq:Error_2} follows directly from Remark \ref{rem:fibers_2} so we need only prove the upper bound. Also, since both $\widetilde{\varepsilon}_n^{\rm E_c}$ and $\varepsilon_n$ are $\mathbb{L}^*$-periodic functions on $\mathbb{R}^d$, it suffices to establish the error estimate~\eqref{eq:Error_2} only for $\bold{k} \in \Omega^*$.

		From Definition \ref{def:basis_2} of the $\bold{k}$-dependent basis set $\mathcal{B}_{\bold{k}}^{\rm E_c}$ (see also Remark \ref{rem:basis_2}), we deduce the existence of some ${\rm E_c'} >0$ such that $n={\rm dim}\; {\rm Y}_{n}^{\bold{k}}\leq M_{\rm {E_c'}}^-$ where we recall that $M_{\rm {E_c'}}^-$ denotes the minimal size of the basis $\mathcal{B}_{\bold{k}}^{\rm {E_c'}}$ over all $\bold{k} \in \mathbb{R}^d$. Consequently, for each $\bold{k} \in \Omega^*$ there also exists ${\rm E_c^{\bold{k}}} \geq {\rm {E_c'}}$ such that the space
		\begin{align*}
			\mY_{n}^{\bold{k}, \rm E_c}:= \{\Psi:=\Pi_{\bold{k}, {\rm E_c^{\bold{k}}}}\Phi \colon \Phi \in  \mY_{n}^{\bold{k}} \},
		\end{align*}
		is $n$-dimensional. We set ${\rm \widetilde{E_c}}:= \sup_{\bold{k} \in \Omega^*} \rm E_c^{\bold{k}}$.

		Using now the min-max theorem, we deduce that for all ${\rm E_c} \geq {\rm \widetilde{E_c}}$ and any $\bold{k} \in \Omega^*$ it holds that
		\begin{align*}
			\widetilde{\varepsilon}_{n,{\bold k}}^{\rm 4E_c} &\leq  \max_{\substack{\Psi \in \mY_{n}^{\bold{k}, \rm {E_c}}\\[0.2em] \Vert \Psi \Vert_{L^2_{\rm per}(\Omega)}=1 } }  \left(\Psi, \widetilde{\mH}_{\bold{k}}^{\mathscr{G}, \rm 4E_c} \Psi \right)_{L^2_{\rm per}(\Omega) }\\
			&= \max_{\substack{\Phi \in \mY_{n}^{\bold{k}}\\[0.2em] \Vert \Pi_{\bold{k}, {\rm {E_c}}}\Phi \Vert_{L^2_{\rm per}(\Omega)} =1}} \left(\Phi, \left(\Pi_{\bold{k}, {\rm {E_c}}} \widetilde{\mH}_{\bold{k}}^{\mathscr{G}, 4\rm E_c} \Pi_{\bold{k}, {\rm {E_c}}}\right)\Phi \right)_{L^2_{\rm per}(\Omega) }\\
			&\leq\max_{\substack{\Phi \in \mY_{n}^{\bold{k}}\\[0.2em] \Vert \Pi_{\bold{k}, {\rm {E_c}}}\Phi \Vert_{L^2_{\rm per}(\Omega)} =1}} \left(\Phi,  \left(\Pi_{\bold{k}, {\rm {E_c}}}\widetilde{\mH}_{\bold{k}}^{\mathscr{G}, \rm 4E_c}\Pi_{\bold{k}, {\rm {E_c}}} -\mH_{\bold{k}}\right)\Phi \right)_{L^2_{\rm per}(\Omega) } \\
			&+\max_{\substack{\Phi \in \mY_{n}^{\bold{k}}\\[0.2em] \Vert \Pi_{\bold{k}, {\rm {E_c}}}\Phi \Vert_{L^2_{\rm per}(\Omega)} =1}} \Big( \underbrace{\left(\Phi,  \mH_{\bold{k}}\Phi\right)_{L^2_{\rm per}(\Omega) }  -\varepsilon_n(\bold{k})}_{\leq \;0} + \varepsilon_n(\bold{k})\Vert \Phi \Vert^2_{L^2_{\rm per}(\Omega)}\Big).
		\end{align*}
		It therefore follows that for all ${\rm E_c} \geq {\rm \widetilde{E_c}}$ and for any $\bold{k} \in \Omega^*$ we have
		\begin{equation}\label{eq:prop1_1}
		\begin{split}
			\widetilde{\varepsilon}_{n,{\bold k}}^{\rm 4E_c} - \varepsilon_{n,{\bold k}} \leq & 
			\max_{\substack{\Phi \in \mY_{n}^{\bold{k}}\\[0.2em] \Vert \Pi_{\bold{k}, {\rm {E_c}}}\Phi \Vert_{L^2_{\rm per}(\Omega)} =1}} \underbrace{\left(\Phi,  \left(\Pi_{\bold{k}, {\rm {E_c}}}\widetilde{\mH}_{\bold{k}}^{\mathscr{G}, \rm 4E_c}\Pi_{\bold{k}, {\rm {E_c}}} -\mH_{\bold{k}}\right)\Phi \right)_{L^2_{\rm per}(\Omega) }}_{:= (\rm I)}\\[0.2em]
			+& \hspace{-2mm}\max_{\substack{\Phi \in \mY_{n}^{\bold{k}}\\[0.2em] \Vert \Pi_{\bold{k}, {\rm {E_c}}}\Phi \Vert_{L^2_{\rm per}(\Omega)} =1}}  \hspace{-6mm}\varepsilon_n(\bold{k})\underbrace{\Vert \Pi^{\perp}_{\bold{k}, {\rm E_c}} {\Phi} \Vert^2_{L^2_{\rm per}(\Omega)}}_{:= (\rm II)}.
			\end{split}
		\end{equation}
		
		Let us first simplify the term (I) for an arbitrary $\Phi \in \mY_{n}^{\bold{k}}$ and $\bold{k} \in \Omega^*$. We  begin by rewriting the term (I) as
		\begin{equation}\label{eq:prop1_2}
		\begin{split}
			{(\rm I)}&= \underbrace{\left(\Phi,  \left(\Pi_{\bold{k}, {\rm {E_c}}}\widetilde{\mH}_{\bold{k}}^{\mathscr{G}, \rm 4E_c}\Pi_{\bold{k}, {\rm {E_c}}} -{\mH^{\rm E_c}_{\bold{k}}}\right)\Phi \right)_{L^2_{\rm per}(\Omega) }}_{:= \rm I(A)}\\[0.2em]
			&+ \underbrace{\left(\Phi,  \left(\mH_{\bold{k}}^{\rm E_c} -{\mH_{\bold{k}}}\right)\Phi \right)_{L^2_{\rm per}(\Omega) }}_{:= \rm I(B)}.
			\end{split}
		\end{equation}
		
		We claim that the term I(A) is identically zero. Indeed, recalling Definition~\ref{def:basis_2} of the basis $\mathcal{B}_{\bold{k}}^{\rm E_c}$ as well as the definitions of the exact fiber and modified Hamiltonian matrix given by Equations \eqref{eq:fiber} and \eqref{eq:fiber_mod} respectively, and using the fact that $\Pi_{\bold{k}, {\rm {4E_c}}}\Pi_{\bold{k}, {\rm {E_c}}}=\Pi_{\bold{k}, {\rm {E_c}}}=\Pi_{\bold{k}, {\rm {E_c}}}\Pi_{\bold{k}, {\rm {4E_c}}}$ we deduce that 
		\begin{align*}
			&\left(\Pi_{\bold{k}, {\rm {E_c}}}\widetilde{\mH}_{\bold{k}}^{\mathscr{G}, \rm 4E_c}\Pi_{\bold{k}, {\rm {E_c}}} -{\mH^{\rm E_c}_{\bold{k}}}\right)\Phi\\[0.5em]
			= \Pi_{\bold{k}, {\rm {E_c}}}&\Bigg(\sum_{\substack{\bold{G} \in \mathbb{L}^*\\[0.2em] \frac{1}{2}\vert\bold{k}+ \bold{G}\vert^2 < {\rm E_c}}} \widehat{\Phi}_{\bold{G}} \left(4{\rm E_c} \mathscr{G} \left(\frac{\vert\bold{k}+ \bold{G}\vert}{\sqrt{8\rm E_c}}\right) - \frac{1}{2}\vert\bold{G} + \bold{k}\vert^2\right) e_{\bold{G}}\Bigg).
		\end{align*}
		Since $\frac{1}{2}\vert \bold{k} +\bold{G}\vert^2 < {\rm E_c} \implies \frac{\vert \bold{k} + \bold{G} \vert}{\sqrt{8 \rm E_c} } < \frac{1}{2}$, we obtain from Definition \ref{def:g} of $\mathscr{G}$ that 
		\begin{align}\label{eq:prop1_3}
			\left(\Pi_{\bold{k}, {\rm {E_c}}}\widetilde{\mH}_{\bold{k}}^{\mathscr{G}, \rm 4E_c}\Pi_{\bold{k}, {\rm {E_c}}} -{\mH^{\rm E_c}_{\bold{k}}}\right)\Phi=0,
		\end{align}
		which implies the claimed result.
		
		The simplification of the term I(B) is classical but we perform it for the sake of completeness. For ease of exposition, in the sequel we will use ${\rm C}>0$ to denote a generic constant whose value may change from step to step but that remains independent of $\bold{k} \in \mathbb{R}^d$, ${\rm E_c}>0$, and $\Phi \in {\rm Y_n^{\bold{k}}}$.

		Let us begin by noting that since $\mH_{\bold{k}}^{\rm E_c}= \Pi_{\bold{k}, {\rm {E_c}}}\mH_{\bold{k}}\Pi_{\bold{k}, {\rm {E_c}}}$ by definition, Equation \eqref{eq:fiber} implies that
		\begin{align}\nonumber
			{\rm I(B)}&= -2\left(\Phi,  \left(\Pi_{\bold{k}, {\rm E_c}}^{\perp}\mH_{\bold{k}}\Pi_{\bold{k}, {\rm E_c}} \right)\Phi \right)_{L^2_{\rm per}(\Omega) } \hspace{0mm}-\left(\Phi,  \left(\Pi_{\bold{k}, {\rm E_c}}^{\perp}\mH_{\bold{k}}\Pi_{\bold{k}, {\rm E_c}}^{\perp}\right)\Phi \right)_{L^2_{\rm per}(\Omega) }\\ \label{eq:prop1_4}
			&\leq 2\Vert \Pi_{\bold{k}, {\rm E_c}}^{\perp}V\Pi_{\bold{k}, {\rm E_c}}\Phi\Vert_{L^2_{\rm per}(\Omega)}\Vert \Pi_{\bold{k}, {\rm E_c}}^{\perp}\Phi\Vert_{L^2_{\rm per}(\Omega)} \hspace{0mm}- \min\{\varepsilon_1(\bold{k}), 0\}\Vert \Pi_{\bold{k}, {\rm E_c}}^{\perp}\Phi\Vert_{L^2_{\rm per}(\Omega)}^2,
		\end{align}
		
		We simplify this last estimate term-by-term. Since $V \in H^{r}_{\rm per}(\Omega)$, we deduce that each exact eigenfunction $u_{n, \bold{k}},~ n \in \mathbb{N}^*$ is an element of $H^{r+2}_{\rm per}(\Omega)$ (see, e.g., \cite{cances2010numerical}). It follows that $\Phi \in H^{r+2}_{\rm per}(\Omega)$ so that we can write 
		\begin{align}\nonumber
			\Vert \Pi_{\bold{k}, {\rm E_c}}^{\perp}\Phi\Vert_{L^2_{\rm per}(\Omega)}^2= \hspace{-2mm}\sum_{\substack{\bold{G} \in \mathbb{L}^*\\ \frac{1}{2} | \bold{k} + \bold{G}|^2 \geq {\rm E_c}} } \hspace{-2mm}\vert \widehat{\Phi}_{\bold{G}} \vert^2 &\leq  \left(\frac{1}{2{\rm E_c}}\right)^{r+2}\hspace{-2mm}\sum_{\substack{\bold{G} \in \mathbb{L}^*\\ \frac{1}{2} | \bold{k} + \bold{G}|^2 \geq {\rm E_c}} } \left(1+ \vert\bold{k}+\bold{G}\vert^2\right)^{r+2}\vert \widehat{\Phi}_{\bold{G}} \vert^2\\[0.5em] \label{eq:prop1_4prime}
			&\leq  \left(\frac{\rm C}{{\rm E_c}}\right)^{r+2}\sum_{\substack{\bold{G} \in \mathbb{L}^*\\ \frac{1}{2} | \bold{k} + \bold{G}|^2 \geq {\rm E_c}} } \left(1+\vert\bold{G}\vert^2\right)^{r+2}\vert \widehat{\Phi}_{\bold{G}} \vert^2\\
			&\leq \left(\frac{\rm C}{{\rm E_c}}\right)^{r+2}\Vert \Phi\Vert_{H^{r+2}_{\rm per}(\Omega)}^2, \nonumber
		\end{align}
		where the second line follows from the fact that $\bold{k} \in \Omega^*$ so that, in particular, $\vert \bold{k} +\bold{G}\vert < \text{diam}(\Omega^*)+ \vert \bold{G}\vert ~\forall \bold{G} \in \mathbb{L}^*$.
		
		In order to simplify the term in Inequality \eqref{eq:prop1_4} involving the potential $V$, we will use similar tactics. As a first step, we make use of the fact that $V \in L^{\infty}_{\rm per}(\Omega)$ is a multiplicative operator so that we may write
		\begin{align}\nonumber
			\Vert \Pi_{\bold{k}, {\rm E_c}}^{\perp}V\Pi_{\bold{k}, {\rm E_c}}\Phi\Vert^2_{L^2_{\rm per}(\Omega)}&= \hspace{-3mm}\sum_{\substack{\bold{G}'\in \mathbb{L}^* \\ \frac{1}{2}\vert\bold{k}+\bold{G}'\vert^2\geq {\rm E_c} } }\Bigg\vert\sum_{\substack{\bold{G} \in \mathbb{L}^*\\ \frac{1}{2} | \bold{k} + \bold{G}|^2 < {\rm E_c}} } \widehat{V}_{\bold{G}'-\bold{G}} \widehat{\Phi}_{\bold{G}} \Bigg\vert^2 \\[1em] \nonumber
			&\leq \hspace{-3mm}\sum_{\substack{\bold{G}'\in \mathbb{L}^* \\ \frac{1}{2}\vert\bold{k}+\bold{G}' \vert^2\geq {\rm E_c} } }\Bigg\vert\sum_{\substack{\bold{G} \in \mathbb{L}^*\\ \frac{1}{2} | \bold{k} + \bold{G}|^2 < {\rm E_c}} } \hspace{-4mm}1 \hspace{2mm} \hspace{-1mm}\sum_{\substack{\bold{G} \in \mathbb{L}^*\\ \frac{1}{2} | \bold{k} + \bold{G}|^2 < {\rm E_c}} } \hspace{-4mm}\vert\widehat{V}_{\bold{G}'-\bold{G}} \vert^2\vert \widehat{\Phi}_{\bold{G}} \vert^2 \Bigg\vert,\\[1em] \nonumber
			&\leq \sum_{\substack{\bold{G}'\in \mathbb{L}^* \\ \frac{1}{2}\vert\bold{k}+\bold{G}'\vert^2\vert\geq {\rm E_c} } } \hspace{-4mm} {\rm C}{\rm E_c^{\frac{d}{2}}} \hspace{1mm} \hspace{-4mm}\sum_{\substack{\bold{G} \in \mathbb{L}^*\\ \frac{1}{2} | \bold{k} + \bold{G}|^2 < {\rm E_c}} } \hspace{-4mm} \vert\widehat{V}_{\bold{G}'-\bold{G}} \vert^2\vert \widehat{\Phi}_{\bold{G}} \vert^2 \\[1em] \nonumber
			&= {\rm C}{\rm E_c^{\frac{d}{2}}} \hspace{-3mm}\sum_{\substack{\bold{G} \in \mathbb{L}^*\\ \frac{1}{2} | \bold{k} + \bold{G}|^2 < {\rm E_c}} } \hspace{-4mm} \vert \widehat{\Phi}_{\bold{G}} \vert^2 \hspace{-3mm}\sum_{\substack{\bold{G}'\in \mathbb{L}^* \\ \frac{1}{2}\vert\bold{k}+\bold{G}' \vert^2\geq {\rm E_c} } } \vert\widehat{V}_{\bold{G}'-\bold{G}} \vert^2\\[1em] 
			&= {\rm C}{\rm E_c^{\frac{d}{2}}} \hspace{-3mm}\sum_{\substack{\bold{G} \in \mathbb{L}^*\\ \frac{1}{2} | \bold{k} + \bold{G}|^2 < {\rm E_c}} } \hspace{-4mm} \vert \widehat{\Phi}_{\bold{G}} \vert^2 \hspace{-3mm}\sum_{\substack{\bold{R}\in \mathbb{L}^* \\ \frac{1}{2}\vert\bold{k}+\bold{G}+\bold{R} \vert^2\geq {\rm E_c} } } \vert\widehat{V}_{\bold{R}} \vert^2. \label{eq:prop_1_5}
		\end{align}
		where the second step follows from the Cauchy-Schwarz inequality and the third step follows by bounding the number of lattice points inside a $d$-dimensional ball of radius $\sqrt{2\rm E_c}$ centered at $\bold{k} \in \mathbb{R}^d$.
		Using now a similar calculation to the one carried out to arrive at Inequality \eqref{eq:prop1_4prime}, we deduce that
		\begin{align*}
			\sum_{\substack{\bold{R}\in \mathbb{L}^* \\ \frac{1}{2}\vert\bold{k}+\bold{G} + \bold{R}\vert^2\geq {\rm E_c} } } \hspace{-2mm}\vert\widehat{V}_{\bold{R}} \vert^2 &\leq \left(\frac{1}{2{\rm E_c}}\right)^{r}\sum_{\substack{\bold{R}\in \mathbb{L}^*\\ \frac{1}{2}\vert\bold{k}+\bold{G} + \bold{R}\vert^2\geq {\rm E_c} } } \hspace{-4mm}\left(1+ \vert\bold{k}+\bold{G} + \bold{R}\vert^2\right)^{r}\vert \widehat{V}_{\bold{R}}\vert^2\\
			&\leq \left(\frac{1}{2{\rm E_c}}\right)^{r}\sum_{\substack{\bold{R}\in \mathbb{L}^*\\ \frac{1}{2}\vert\bold{k}+\bold{G} + \bold{R}\vert^2\geq {\rm E_c} } } \hspace{-4mm}\left(1+ 2\vert \bold{k}+\bold{R}\vert^2 +2\vert\bold{G}\vert^2\right)^{r}\vert \widehat{V}_{\bold{R}}\vert^2\\
			&\leq  \left(\frac{\rm C}{{\rm E_c}}\right)^{r}\sum_{\substack{\bold{R}\in \mathbb{L}^*\\ \frac{1}{2}\vert\bold{k}+\bold{G} + \bold{R}\vert^2\geq {\rm E_c} } } \hspace{-4mm}\left(1+ \vert \bold{k}+\bold{R}\vert^2 \right)^{r} \vert \widehat{V}_{\bold{R}}\vert^2+\left(1+ \vert\bold{G}\vert^2\right)^{r}\vert \widehat{V}_{\bold{R}}\vert^2\\
			&\leq \left(\frac{\rm C}{{\rm E_c}}\right)^{r}\left(\Vert V\Vert^2_{H^r_{\rm per}(\Omega) }+\left(1+ \vert \bold{G}\vert^2 \right)^{r}\Vert V\Vert^2_{L^2_{\rm per}(\Omega) }\right).
		\end{align*}
		Plugging in this last expression in Inequality \eqref{eq:prop_1_5} easily allows us to deduce that
			\begin{align}\label{eq:prop_1_6}
			\Vert \Pi_{\bold{k}, {\rm E_c}}^{\perp}V\Pi_{\bold{k}, {\rm E_c}}\Phi\Vert^2_{L^2_{\rm per}(\Omega)} \leq \left(\frac{\rm C}{{\rm E_c}}\right)^{r-\frac{d}{2}}\Vert V\Vert^2_{H^r_{\rm per}(\Omega) }\Vert \Phi\Vert^2_{H^{r+2}_{\rm per}(\Omega) } .
		\end{align}
		

		One final, similar calculation allows us to simplify also the term (II) in Inequality~\eqref{eq:prop1_1}. Combining now Estimates~\eqref{eq:prop1_4prime}-\eqref{eq:prop_1_6} with Inequalities \eqref{eq:prop1_1}-\eqref{eq:prop1_4}, we deduce that for all ${\rm E_c} \geq {\rm \widetilde{E_c}}$ and any $\bold{k} \in \Omega^*$ it holds that
		\begin{align*}
			\widetilde{\varepsilon}_{n,{\bold k}}^{\rm 4E_c} - \varepsilon_{n,{\bold k}}  \leq \hspace{-3mm} \max_{\substack{\Phi \in \mY_{n}^{\bold{k}}\\[0.2em] \Vert \Pi_{\bold{k}, {\rm {E_c}}}\Phi \Vert_{L^2_{\rm per}(\Omega)} =1}}\Bigg(&\left(\frac{{\rm C}}{{\rm E_c}}\right)^{{r+1-\frac{d}{4}}}\Vert V\Vert_{H^r_{\rm per}(\Omega)}\Vert \Phi\Vert^2_{H^{{r+2}}_{\rm per}(\Omega)} \\
			-  &\left(\frac{\rm C}{{\rm E_c}}\right)^{r+2} \hspace{-2mm}\min\{\varepsilon_1(\bold{k}), 0\}\Vert \Phi\Vert^2_{H^{r+2}_{\rm per}(\Omega)}\Bigg)\\[1.2em]
			+ \max_{\substack{\Phi \in \mY_{n}^{\bold{k}}\\[0.2em] \Vert \Pi_{\bold{k}, {\rm {E_c}}}\Phi \Vert_{L^2_{\rm per}(\Omega)} =1}}&\left(\frac{\rm C}{{\rm E_c}}\right)^{r+2} \varepsilon_m(\bold{k})\Vert \Phi\Vert^2_{H^{r+2}_{\rm per}(\Omega)}.
		\end{align*}
		Collecting terms, we obtain an estimate of the form
		\begin{align*}
			\widetilde{\varepsilon}_{n,{\bold k}}^{\rm 4E_c} - \varepsilon_{n,{\bold k}}  \leq  \left(\frac{\rm C}{\rm E_c}\right)^{r+1-\frac{d}{4}} \hspace{-2mm}\max_{\substack{\Phi \in \mY_{n}^{\bold{k}}\\[0.2em] \Vert \Pi_{\bold{k}, {\rm {E_c}}}\Phi \Vert_{L^2_{\rm per}(\Omega)} =1}} \Vert \Phi\Vert^2_{H^{r+2}_{\rm per}(\Omega)}.
		\end{align*}
		
		To conclude, we notice that we may write 
		\begin{align*}
			\max_{\substack{\Phi \in \mY_{n}^{\bold{k}}\\[0.2em] \Vert \Pi_{\bold{k}, {\rm {E_c}}}\Phi \Vert_{L^2_{\rm per}(\Omega)} =1}} \Vert \Phi\Vert^2_{H^{r+2}_{\rm per}(\Omega)} \leq \max_{\substack{\Phi \in \mY_{n}^{\bold{k}}\\[0.2em] \Vert\Phi \Vert_{L^2_{\rm per}(\Omega)} =1}} \Vert \Phi\Vert^2_{H^{r+2}_{\rm per}(\Omega)} \hspace{-2mm}\underbrace{\max_{\substack{\Phi \in \mY_{n}^{\bold{k}}\\[0.2em] \Vert \Pi_{\bold{k}, {\rm {E_c}}}\Phi \Vert_{L^2_{\rm per}(\Omega)} =1}} \Vert \Phi\Vert^2_{L^{2}_{\rm per}(\Omega)}}_{:=(\rm III)},
		\end{align*}
		and it is well-known that there exists an upper bound ${\rm C}_{M, {\rm E_c}}\geq0$ for the term (III), provided that the basis~$\mathcal{B}_{\bold{k}}^{\rm E_c}$ satisfies the so-called approximation property (which it does) and that ${\rm E_c} $ is larger than some threshold ${\rm \widehat{E_c}} \geq 0$. Defining appropriate constants and setting the discretization cutoff ${\rm E_c^*}$ large enough thus completes the proof.
	\end{proof}
	
	\vspace{4mm}
	
	We can now turn our attention to the more technical proof of Theorem \ref{prop:2} which characterizes precisely the regularity of the modified energy bands $\{\widetilde{\varepsilon}_{n}^{\rm E_c}\}_{n \in \mathbb{N}^*}$. As stated in Section \ref{sec:5}, the proof of Theorem \ref{prop:2} will require the use of Lemma \ref{lem:contin} which we now prove.  \vspace{4mm}

		\begin{proof}[\textbf{Proof of Lemma \ref{lem:contin}}]~
			
			We begin by claiming that thanks to the assumptions on the matrix ${\rm C}_{n}$, for every $\Upsilon >0$ there exists a natural number $N(\Upsilon) \in \mathbb{N}$ such that for all $n \geq N(\Upsilon)$ and any $\lambda \in \mathbb{C}$ such that $\vert\lambda\vert < \Upsilon$, the inverse matrix $\left({\rm C}_{n}- \lambda\right)^{-1}$ exists. Indeed, this is simply a consequence of the fact that for any invertible matrix $E$, all its eigenvalues are lower bounded in magnitude by~$ \Vert E^{-1} \Vert_2^{-1}$.
			

			Therefore, using a well-known determinant identity for block matrices we deduce that for all $n \geq N(\Upsilon)$ and any $\lambda \in \mathbb{C}$ such that $\vert\lambda\vert < \Upsilon$ it holds that
			\begin{align}\label{eq:proof_new_2}
				\text{\rm det}\left({\mH}_{n}- \lambda\right) &= \text{\rm det}\left({\rm C}_{n}- \lambda\right)\; \text{\rm det}\left({\rm A}_{n}-\lambda -{\rm B}_{n} \left({\rm C}_{n}-\lambda\right)^{-1}\widetilde{{\rm B}_{n}}\right).
			\end{align}
			Equation \eqref{eq:proof_new_2} implies that for any $\Upsilon >0$ sufficiently large there exists a natural number $N(\Upsilon) \in \mathbb{N}$ such that for all $n \geq N(\Upsilon)$
			\begin{align}\nonumber
				&\lambda \in \mathbb{C} \text{ such that }\vert\lambda\vert < \Upsilon \text{ is an eigenvalue of } {\mH}_{n} \quad \\ \label{eq:proof_new_7}
				&\iff \\ \nonumber
				&\lambda \in \mathbb{C} \text{ such that }\vert\lambda\vert < \Upsilon \text{ is a solution to } \text{\rm det}\left({\rm A}_{n}-\lambda -{\rm B}_{n} \left({\rm C}_{n}-\lambda\right)^{-1}\widetilde{{\rm B}_{n}}\right)=0.
			\end{align}

			In view of the relation \eqref{eq:proof_new_7}, we are motivated to define for any $\Upsilon >0$ and any $n \geq N(\Upsilon)$, the open disk $ {\mathbb{B}}_{\Upsilon}(0):= \left\{z \in \mathbb{C} \colon \vert z\vert < \Upsilon\right\}$ and the non-linear function $\mathcal{g}_{\Upsilon, n} \colon {\mathbb{B}}_{\Upsilon}(0) \rightarrow \mathbb{C}$ given by
			\begin{align*}
				\mathcal{g}_{\Upsilon, n}(\lambda):=\text{\rm det}\left({\rm A}_{n}-\lambda -{\rm B}_{n} \left({\rm C}_{n}- \lambda\right)^{-1}\widetilde{{\rm B}_{n}}\right) \qquad \forall \lambda \in {\mathbb{B}}_{\Upsilon}(0).
			\end{align*}
			
			We are now interested in studying the zeros of the function $\mathcal{g}_{\Upsilon, n}$ in the asymptotic regime $n \to \infty$. Our goal is to show that the sequence of functions $\{\mathcal{g}_{\Upsilon, n}\}_{n \in \mathbb{N}}$ satisfy the hypotheses of Hurwitz's theorem from complex analysis (see, e.g., \cite[Chapter VII, Theorem 2.5]{MR503901}). 
			
			We begin by establishing that for any $\Upsilon >0$ and any $n \geq N(\Upsilon)$, the non-linear function $\mathcal{g}_{\Upsilon, n}$ is holomorphic on the open disk~$\mathbb{B}_{\Upsilon}(0)$. To do so, observe that for any $\Upsilon >0$, all $n\geq N(\Upsilon)$ and any $\lambda \in {\mathbb{B}}_{\Upsilon}(0) $ we can define the matrix
			\begin{align*}
				{\rm Z}_{n}(\lambda):={\rm A}_{n}-\lambda -{\rm B}_{n} \left({\rm C}_{n}- \lambda\right)^{-1}\widetilde{{\rm B}_{n}}.
			\end{align*}
			
			Recall that the natural number $N(\Upsilon)$ was chosen so that all eigenvalues of ${\rm C}_{n}$ are strictly larger in magnitude than $\Upsilon$ and consequently, $\left({\rm C}_{n}- \lambda\right)^{-1}$ exists for all $n \geq N(\Upsilon)$ and all $\lambda \in {\mathbb{B}}_{\Upsilon}(0) \subset \mathbb{C}$. In view of the assumptions on the sub-matrices ${\rm A}_{n}, {\rm B}_{n}, \widetilde{{\rm B}_{n}}$ and ${\rm C}_{n}$, it therefore follows that the matrix $Z_{n}(\lambda)$ exists and is bounded for all $n \geq N(\Upsilon)$, and has a power series expansion in the disk ${\mathbb{B}}_{\Upsilon}(0)$ of the form
			\begin{align}\label{eq:proof_new_6}
				{\rm Z}_{n}(\lambda)=  {\rm A}_{n}-\lambda -\sum_{q=0}^{\infty}{\rm B}_{n} \lambda^{q}\left({\rm C}_{n}\right)^{-q-1}\widetilde{{\rm B}_{n}}=\sum_{q=0}^{\infty} \lambda^q \;M_{q, n},
			\end{align}
			where each ${M}_{q, n}$ is a square matrix of dimension $p= \text{dim } {\rm A}$. 
			
			Equation \eqref{eq:proof_new_6} implies that each entry of the matrix $ {\rm Z}_{n}(\lambda)$ is itself a holomorphic function of $\lambda$ with a power series expansion valid in $\mathbb{B}_{\Upsilon}(0)$. Moreover, since $\mathcal{g}_{\Upsilon, n}(\lambda)= \text{\rm det}\left({\rm Z}_{n}(\lambda)\right)$ for each $\lambda \in \mathbb{B}_{\Upsilon}(0)$, and the determinant is a polynomial of the entries of the underlying matrix, we deduce that for any $\Upsilon >0$ and all $n\geq N(\Upsilon)$ the non-linear function $\mathcal{g}_{\Upsilon, n}$ is indeed holomorphic on~${\mathbb{B}}_{\Upsilon}(0) \subset \mathbb{C}$.

			Next, we claim that on any compact set $K\subset {\mathbb{B}}_{\Upsilon}(0)$, the sequence of non-linear functions $\{\mathcal{g}_{\Upsilon, n}\}_{\ell \in \mathbb{N}}$ converges uniformly to the characteristic polynomial of ${\rm A}_{n}$, i.e., 
			\begin{align}\label{eq:proof_new_3}
				\lim_{n \to \infty} \sup_{\lambda \in K} \left\vert \text{\rm det}\left({\rm A}_{n}-\lambda\right)- \mathcal{g}_{\Upsilon, n}(\lambda) \right\vert = 0.
			\end{align}
			
			To prove that Equation \eqref{eq:proof_new_3} indeed holds, we appeal to a known determinant inequality for differences of matrices (see \cite[Theorem 2.12]{MR2421470}): For any two matrices $ E, F \in \mathbb{C}^{p\times p}$, it holds that
			\begin{align}\label{eq:proof_new_4}
				\vert \text{\rm det}(E) -\text{\rm det}(E+F) \vert \leq p \Vert F \Vert_{2}\; \max \left\{\Vert E \Vert_{2}, \Vert E+F\Vert_2\right\}^{p-1}.
			\end{align}
			
			Applying Inequality \eqref{eq:proof_new_4} to our situation yields that for any $\Upsilon >0$, any compact set $K\subset {\mathbb{B}}_{\Upsilon}(0)$, all $\lambda \in K $ and all~$n \geq N(\Upsilon)$ it holds that
			\begin{align*}
				&\vert \text{\rm det}\left({\rm A}_{n} - \lambda\right) -\mathcal{g}_{\Upsilon, n}(\lambda) \vert  \\
				\leq  ~p &\big\Vert {\rm B}_{n} \left({\rm C}_{n}- \lambda\right)^{-1}\widetilde{{\rm B}_{n}}\big\Vert_2 \max \Big\{ \left\Vert {\rm A}_{n}-\lambda\right\Vert_2, \big\Vert {\rm A}_{n}-\lambda -{\rm B}_{n} \left({\rm C}_{n}- \lambda\right)^{-1}\widetilde{{\rm B}_{n}}\big\Vert_2 \Big\}^{p-1}.
			\end{align*}
			
			Using now the assumptions on the sub-matrices ${\rm A}_{n}, {\rm B}_{n},\widetilde{{\rm B}_{n}}$ and ${\rm C}_{n}$, we deduce that for any $\Upsilon >0$ and any compact set $K\subset {\mathbb{B}}_{\Upsilon}(0)$, there exists a constant $\Lambda_{\Upsilon, K}$ (which depends also on $p$) such that for all $n \geq N(\Upsilon)$ and all $\lambda \in K$ it holds that
			\begin{align*}
				\left\vert \text{\rm det}\left({\rm A}_{n} - \lambda\right) -\mathcal{g}_{\Upsilon, n}(\lambda) \right\vert \leq \Lambda_{\Upsilon, K}\; \Vert {\rm C}_{n}^{-1}\Vert_2,
			\end{align*}
			from which Equation \eqref{eq:proof_new_3} now readily follows. Let us also emphasize here that since $\lim_{n \to \infty} {\rm A}_{n}={\rm A} $, we have in fact shown that the sequence of non-linear functions $\{\mathcal{g}_{\Upsilon, \ell}\}_{\ell \in \mathbb{N}}$ converges uniformly to the characteristic polynomial of ${\rm A}$ on any compact set $K \subset {\mathbb{B}}_{\Upsilon}(0)$.
			
			In order to complete our analysis, we observe that the characteristic polynomial of ${\rm A}$ is an \emph{entire} function which is not identically zero on any open subset of $\mathbb{C}$. As a consequence, Hurwitz's theorem can be applied: For any $\Upsilon >0$ and any non-empty open, connected set $U$ such that $\overline{U} \subset {\mathbb{B}}_{\Upsilon}(0)$ and $\text{det}({\rm A}-\lambda)\neq 0~ \forall\lambda \in \partial U$,  there exists $N(\Upsilon)\in \mathbb{N}$ such that for all $n\geq N(\Upsilon)$, the non-linear function $\mathcal{g}_{\Upsilon, n}$ and the characteristic polynomial $\text{det}({\rm A}-\bullet)$ of the matrix~${\rm A} $ have the same number of zeros in $U$ counting multiplicity. In particular,

			\begin{enumerate} 
				
				\item[(i)] for all $\Upsilon >0$ sufficiently large, there exists $\widetilde{N}(\Upsilon) \in \mathbb{N}$ such that for all $n\geq \widetilde{N}(\Upsilon)$ the non-linear function $\mathcal{g}_{\Upsilon, n}$ has \underline{exactly} $p= \text{dim } { \rm A}$ zeros counting multiplicity with magnitude strictly smaller than $\Upsilon$;

				\item[(ii)] picking a fixed $\Upsilon$ large enough so that ${\mathbb{B}}_{\Upsilon}(0)$ contains all roots of the characteristic polynomial of ${\rm A}$, for every $\rho>0$ sufficiently small, there exists $N(\rho)\in \mathbb{N}$ such that for any eigenvalue $\lambda^{\rm A}$ of the matrix ${\rm A}$ with algebraic multiplicity $Q \in \mathbb{N}^*$, and all $n \geq N(\rho)$, the non-linear function $\mathcal{g}_{\Upsilon, n}$ has \underline{exactly} $Q$ zeros (counting multiplicity) in the open disk $\mathbb{B}_{\rho}\big(\lambda^{\rm A}\big) \subset \mathbb{C}$ ;
			\end{enumerate}
			The proof now follows easily by making use of Relation \eqref{eq:proof_new_7}.

		\end{proof}
	
	\vspace{4mm}
	
	\begin{proof}[\textbf{Proof of Theorem \ref{prop:2}}]~
		
		For clarity of exposition, we will divide this proof into three portions: We will first consider the regularity of the approximate energy bands away from crossings and away from changes in the cardinality of the $\bold{k}$-dependent basis sets. This will allows us to deduce, as a corollary, the regularity of the approximate energy bands at crossings but under the assumption that the cardinality of the $\bold{k}$-dependent basis set does not change. Lastly, we will consider the regularity of the approximate energy bands in the neighborhood of points where the cardinality of the $\bold{k}$-dependent basis sets may change.
		
		For the remainder of this proof we recall the setting of the $\bold{k}$-dependent modified Galerkin discretization \eqref{eq:Galerkin_3}, we select ${\rm E_c}>0$ such that $M_{\rm E_c}^->0$, and we pick some index $n \in \{1, \ldots, M_{\rm E_c}^- \}$ and some point $\bold{k}_0 \in \mathbb{R}^d$.

		\vspace{4mm}
	\noindent	\textbf{Case one:}~ we assume that $\bold{k}_0 \in \mathbb{R}^d$ satisfies
		\begin{align*}
			&\text{For all } \bold{G} \in \mathbb{L}^* \text{ it holds that } \vert \bold{k}_0+ \bold{G}\vert^2 \neq 2 {\rm E_c} \quad \text{and}\\[0.5em]
			&\widetilde{\varepsilon}^{\rm E_c}_n (\bold{k}_0) \neq \widetilde{\varepsilon}^{\rm E_c}_{\widetilde{n}} (\bold{k}_0) \quad \forall {\widetilde{n}} \in \{1, \ldots, M_{\rm E_c}^-\} \quad \text{ with }  \quad {\widetilde{n}} \neq n.
		\end{align*}
		In other words, we assume that there is no change in the cardinality of the basis set at $\bold{k}_0$ and that there are no band crossings at $(\bold{k}_0, \widetilde{\varepsilon}^{\rm E_c}_n(\bold{k}_0))$.
		
		\vspace{2mm}
		
		We claim that in this case, the approximate energy band $\widetilde{\varepsilon}^{\rm E_c}_n$ is of class $\mathscr{C}^m$ in a neighborhood of $\bold{k}_0$, i.e., $\widetilde{\varepsilon}^{\rm E_c}_n$ has the same local regularity at $\bold{k}_0$ as the blow-up function $\mathscr{G}$ does on the interval $(0, 1) \subset \mathbb{R}$. Indeed, this is a straightforward application of the implicit function theorem: We notice that the dimensions of the modified-Hamiltonian matrices $\widetilde{\mH}_{\bold{k}}^{\mathscr{G}, \rm E_c}$ do not change in a sufficiently small neighborhood of $\bold{k}_0$ and the dependence of this matrix on $\bold{k}$ in such a neighborhood is of class $\mathscr{C}^m$ (thanks to the regularity properties of $\mathscr{G}$). Since, additionally, the eigenvalue $\widetilde{\varepsilon}^{\rm E_c}_n (\bold{k}_0)$ is simple, it can be shown that the assumptions of the implicit function theorem hold, and therefore by a classical argument (see, e.g.,\cite[Theorem 5.3]{MR2744852}) it follows that the approximate energy band $\widetilde{\varepsilon}^{\rm E_c}_n$ is indeed of class $\mathscr{C}^m$ in a neighborhood of $\bold{k}_0$ as claimed. 
		
		Let us remark here that a similar argument involving the implicit function theorem yields $\mathscr{C}^m$ regularity, as a function of $\bold{k} \in \mathbb{R}^d$, of the (normalized) eigenfunction $\widetilde{u}^{\rm E_c}_{n, \bold{k}}$ associated with the eigenvalue $\widetilde{\varepsilon}^{\rm E_c}_{n, \bold{k}}$ at $\bold{k}=\bold{k}_0$. A detailed argument can, for instance, be found in \cite[Chapter 9, Theorem 8]{MR2356919}. This additional fact will of use in the sequel.

		\vspace{4mm}
	\noindent	\textbf{Case two:}~ we assume that $\bold{k}_0 \in \mathbb{R}^d$ satisfies
			\begin{align*}
			&\text{For all } \bold{G} \in \mathbb{L}^* \text{ it holds that } \vert \bold{k}_0+ \bold{G}\vert^2 \neq 2 {\rm E_c} \quad \text{and}\\[0.5em]
			&\exists ~{\widetilde{n}} \in \{1, \ldots, M_{\rm E_c}^-\} \text{ with } {\widetilde{n}} \neq n \colon \quad \widetilde{\varepsilon}^{\rm E_c}_n (\bold{k}_0) = \widetilde{\varepsilon}^{\rm E_c}_{\widetilde{n}} (\bold{k}_0).
		\end{align*}
		In other words, we assume that there is no change in the cardinality of the basis set at $\bold{k}_0$ but there is a band crossing at $\bold{k}_0, \widetilde{\varepsilon}^{\rm E_c}_n(\bold{k}_0))$.
		
		\vspace{2mm}
		
		
		We claim that in this case, the approximate energy band $\widetilde{\varepsilon}^{\rm E_c}_n$ is either Lipschitz continuous in a neighborhood of $\bold{k}_0$ if $m\geq 1$ or of class $\mathscr{C}^{0}$ otherwise. To this end, we notice that the dimensions of the matrix $\widetilde{\mH}_{\bold{k}}^{\mathscr{G}, \rm E_c}$ do not change in a sufficiently small neighborhood of $\bold{k}_0$ and the dependence of this matrix on $\bold{k}$ in such a neighborhood is of class $\mathscr{C}^m, ~m\geq 0$ thanks to the regularity properties of $\mathscr{G}$. If $m=0$, then it follows from a classical argument (see \cite[Appendix~V, Page 363]{MR0387634}) that all approximate energy bands $\{\widetilde{\varepsilon}_n^{\rm E_c}\}_{n=1}^{M_{\rm E_c}^-}$ are continuous at $\bold{k}_0$. If, on the other hand, $m\geq 1$, then the claimed Lipschitz continuity follows from the min-max theorem.

		\vspace{4mm}  
	\noindent	\textbf{Case three:}~ we assume that $\bold{k}_0 \in \mathbb{R}^d$ satisfies
		\begin{align*}
			\text{There exists } \bold{G} \in \mathbb{L}^* \text{ such that  } \vert \bold{k}_0+ \bold{G}\vert^2 = 2 {\rm E_c}.
		\end{align*}
			In other words, we assume that there \emph{is} a change in the cardinality of the basis set at $\bold{k}_0$.
			
			\vspace{2mm}
			

		We claim that in this case, there are two possibilities: if there is no band crossing at $(\bold{k}_0, \widetilde{\varepsilon}^{\rm E_c}_n(\bold{k}_0))$, then the approximate energy band $\widetilde{\varepsilon}^{\rm E_c}_n$ is of class $\mathscr{C}^m$ in a neighborhood of $\bold{k}_0$, i.e., $\widetilde{\varepsilon}^{\rm E_c}_n$ has the same local regularity at $\bold{k}_0$ as the rate of blow-up of the function $\mathscr{G}(x)$ in the limit $x \to 1$. If, on the other hand, there is a band crossing at $(\bold{k}_0, \widetilde{\varepsilon}_n(\bold{k}_0))$, then the approximate energy band $\widetilde{\varepsilon}^{\rm E_c}_n$ is Lipschitz continuous at $\bold{k}_0$ if $m\geq 1$ and of class $\mathscr{C}^0$ otherwise. \vspace{2mm}

		We begin by defining the non-empty sets
		\begin{align*}
			{\mS}^{\rm E_c^-}_{\bold{k}_0}:=& \left \{\bold{G} \in \mathbb{L}^* \colon \frac{1}{2}\vert \bold{k}_0+ \bold{G}\vert^2 < {\rm E_c}\right \} \quad \text{and}\\
			{\mS}^{\rm E_c}_{\bold{k}_0}:=& \left \{\bold{G} \in \mathbb{L}^* \colon \frac{1}{2}\vert \bold{k}_0+ \bold{G}\vert^2 = {\rm E_c}\right \} \quad \text{with}\quad  M:= {\rm dim} \;{\mS}_{\bold{k}_0}^{\rm E_c} < \infty.
		\end{align*}
		
		Clearly there exists $\overline{\delta} > 0$ sufficiently small such that $\forall ~\bold{k} \in \mathbb{R}^d$ with $\vert \bold{k}- \bold{k}_0 \vert < \overline{\delta}$, the $\bold{k}$-dependent basis sets satisfy
		\begin{align*}
			\left\{ e_\bold{G} \colon \bold{G} \in {\mS}^{\rm E_c^-}_{\bold{k}_0} \right\}  \subseteq		~\mathcal{B}_{\rm \bold{k}}^{\rm E_c} ~ \subseteq \left\{ e_\bold{G} \colon \bold{G} \in {\mS}^{\rm E_c^-}_{\bold{k}_0} \cup {\mS}^{\rm E_c}_{\bold{k}_0}\right\} .
		\end{align*}
		Let us therefore fix some $\delta \leq \overline{\delta}$ and consider the open ball $\mathbb{B}_{\delta}(\bold{k}_0)$ of radius $\delta$ centered at $\bold{k}_0$.  We will study the behavior of sequences of $\bold{k}$-points in this open ball that converge to $\bold{k}_0$.

		To do so, we consider a specific decomposition of the open ball $\mathbb{B}_{\delta}(\bold{k}_0)$ into sectors $\{\Omega_j\}_{j=1}^M$ defined as follows: First, for every $\widetilde{\bold{G}} \in {\mS}^{\rm E_c}_{\bold{k}_0}$ we define the open set
		\begin{align}\label{eq:Cone}
			{\mS}_{\widetilde{\bold{G}}}= \left\{\bold{k} \in \mathbb{B}_{\delta}(\bold{k}_0) \colon ~ \frac{1}{2}\vert \bold{k} + \widetilde{\bold{G}}\vert^2 < {\rm E_c}\right\}.
		\end{align}	
		It is now easy to see that there are exactly two cases:
		\begin{enumerate}
			\item For all $\bold{k} \in {\mS}_{\widetilde{\bold{G}}}$, the Fourier mode $e_{\widetilde{\bold{G}}}$ is an element of the $\bold{k}$-dependent basis set $\mathcal{B}_{\rm \bold{k}}^{\rm E_c}$.
			
			\item For all $\bold{k} \in \mathbb{B}_{\delta}(\bold{k}_0)\setminus {\mS}_{\widetilde{\bold{G}}}$, the Fourier mode $e_{\widetilde{\bold{G}}}$ is \underline{not} an element of the $\bold{k}$-dependent basis set $\mathcal{B}_{\rm \bold{k}}^{\rm E_c}$.
			
		\end{enumerate}
		
		\vspace{2mm}

		Next, we label the elements of ${\mS}^{\rm E_c}_{\bold{k}_0}$ as $\bold{G}_1, \bold{G}_2, \ldots, \bold{G}_M$. It follows that there exist sets $\Omega_j, ~j\in \{1, \ldots, 2^M\}~\subset~\mathbb{B}_{\delta}(\bold{k}_0)$ such that $\mathbb{B}_{\delta}(\bold{k}_0) = \cup_{j=1}^{2^M} \Omega_j$ with
		\begin{equation}\label{eq:prop_2_1}
			\begin{split}
				&\Omega_1 := \mathbb{B}_{\delta}(\bold{k}_0) \setminus \Big(\bigcup_{\bold{G} \in {\mS}^{\rm E_c}_{\bold{k}_0}} \mS_{\bold{G}}\Big) \quad \text{and}\\
				&\forall j \in \{2, \ldots, 2^M\}, ~\exists~ L \leq M, ~ J:=\{j_1, j_2, \ldots, j_L\}\subset \{1, \ldots, M\} ~ \text{such that} \\ 
		&	\hphantom{\forall j \in asdaasd}\Omega_j=\left(\bigcap_{\ell \in J} \mS_{\bold{G}_{\ell}}\right) \setminus \left(\bigcup_{\ell \in \{1, \ldots, M\}\setminus J} \mS_{\bold{G}_{\ell}}\right).
			\end{split}
		\end{equation}
		A visual example of the above decomposition is displayed in Figure \ref{fig:decomp}. Note that we allow for the possibility of some $\Omega_j, ~j\in \{1, \ldots, 2^M\}$ to be empty. Two observations should now be made. 
		\begin{figure}[ht]
			\centering
			\includegraphics[width=0.75\textwidth, trim={0cm, 2cm, 0cm, 0cm},clip=true]{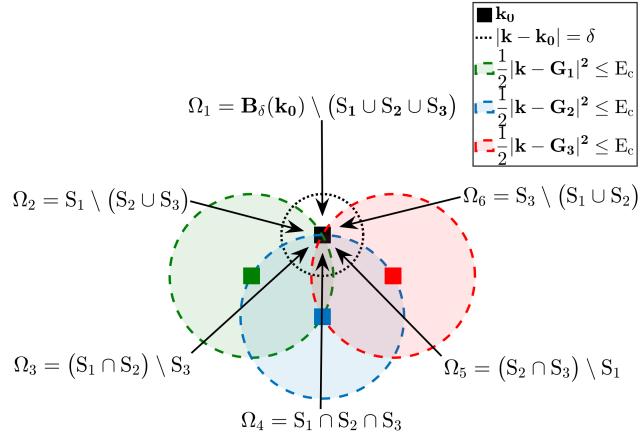} 
			\caption{An example of Decomposition \eqref{eq:prop_2_1} introduced above for a hexagonal lattice. For clarity, we show only the decomposition corresponding to three elements of ${\mS}^{\rm E_c}_{\bold{k}_0}$.}
			\label{fig:decomp}
		\end{figure}

		\vspace{3mm}
		\noindent \textbf{Observation one:} the set $\Omega_{1}$ has the following property: For all $\bold{k} \in \Omega_1$  it holds that $ \mathcal{B}_{\bold{k}}^{\rm E_c}=\mathcal{B}_{\rm \bold{k}_0}^{\rm E_c}$. 
		
		\vspace{3mm}

	\noindent	\textbf{Observation two:} each set $\Omega_{j}, ~j >1$ has the following property: For all $\bold{k} \in \Omega_j$ it holds that
		\begin{align*}
			\mathcal{B}_{\bold{k}}^{\rm E_c} = \; \mathcal{B}_{\rm \bold{k}_0}^{\rm E_c} \cup \{e_{\bold{G}_{j_1}},  e_{\bold{G}_{j_2}}, \ldots, e_{\bold{G}_{j_L}}\}.
		\end{align*}

		Let us now fix some $j \in \{1, \ldots, 2^M\}$ and consider the set $\Omega_{j}$. Our goal is to study the convergence of the approximate, modified energy band $\widetilde{\varepsilon}_n^{\rm E_c}(\cdot)$ (and its derivatives) along sequences $\left\{\bold{k}_{\ell}\right\}_{\ell \in \mathbb{N}^*}$ contained in $\Omega_j$ that converge to $\bold{k}_0$. The reason we have introduced the decomposition $\{\Omega_j\}_{j=1}^{2^M}$ and we restrict ourselves, as a first step, to sequences contained in a fixed $\Omega_j$ is because for each fixed $\Omega_j$ and all $\bold{k} \in \Omega_j$, we have precise knowledge of the $\bold{k}$-dependent basis set thanks to \textbf{Observation two}.

		Obviously, if $\Omega_j$ is an empty set, then there is nothing to study so we assume without loss of generality that $\Omega_j$ is non-empty. Additionally, thanks to \textbf{Observation one} above, the choice $\Omega_{j}=\Omega_1$ is already covered by the proof of \textbf{Case one} of our proof so we may assume that $j >1$. It now follows from \textbf{Observation two} that for all $\bold{k} \in \Omega_j$,  we have the decomposition
		\begin{align}\label{eq:prop_2_2a}
			\mathcal{B}_{\rm \bold{k}}^{\rm E_c}=& ~\mathcal{B}_{\rm \bold{k}_0}^{\rm E_c} \cup~ \widetilde{\mathcal{B}_{\Omega_j}}^{\rm E_c} \quad \text{and} \quad \mathcal{B}_{\rm \bold{k}_0}^{\rm E_c} \cap~ \widetilde{\mathcal{B}_{\Omega_j}}^{\rm E_c}= \emptyset, \qquad \text{where}\\
		\label{eq:prop_2_2b}
			\widetilde{\mathcal{B}_{\Omega_j}}^{\rm E_c}:=&~ \left\{e_{\bold{G}}\colon ~\bold{G} \in \{\bold{G}_{j_1}, \bold{G}_{j_2},\ldots, \bold{G}_{j_L}\}\right\} \quad \text{and} \quad \widetilde{\mX}_{\bold{k}_0, \Omega_{j}}^{\rm E_c}:= \text{span}\;\widetilde{\mathcal{B}_{\Omega_j}}^{\rm E_c}.
		\end{align}
		Let us remark here that in Equation \eqref{eq:prop_2_2a}, the fact that the set intersection is empty follows from the fact that if the Fourier mode~$e_{\bold{G}_k} \in \widetilde{\mathcal{B}_{\Omega_j}}^{\rm E_c}$, then $\bold{G}_k \in \mS_{\bold{k}_0}^{\rm E_c}$ so that $e_{\bold{G}_k} \notin \mathcal{B}_{\rm \bold{k}_0}^{\rm E_c} $ by Definition \ref{def:basis_2}.

		From Equations \eqref{eq:prop_2_2a} and \eqref{eq:prop_2_2b} we can now deduce that $ \forall \bold{k}\in \Omega_j$ it holds that
		\begin{align*}
			\mX_{\bold{k}}^{\rm E_c}:=& \text{\rm span}\; \mathcal{B}_{\rm \bold{k}}^{\rm E_c} = \text{span}\;\mathcal{B}_{\rm \bold{k}_0}^{\rm E_c}  \oplus \text{span}\; \widetilde{\mathcal{B}_{\Omega_j}}^{\rm E_c}= \mX_{\bold{k}_0}^{\rm E_c}  \oplus\widetilde{\mX}_{\bold{k}_0, \Omega_{j}}^{\rm E_c}.
		\end{align*}
		
		We denote by $\widetilde{\Pi}_{\Omega_j, \rm E_c} \colon L^2_{\rm per}(\Omega) \rightarrow L^2_{\rm per}(\Omega)$ the $L^2_{\rm per}$-orthogonal projection operator onto $\widetilde{\mX}_{\bold{k}_0, \Omega_{j}}^{\rm E_c}$. It follows that for all $\bold{k} \in \Omega_j$, the modified Hamiltonian matrix $\widetilde{\mH}_{\bold{k}}^{\mathscr{G}, \rm E_c}$ admits the following block representation:

		\begin{equation}\label{eq:prop_2_3}
			\arraycolsep=1.4pt\def\arraystretch{4} \widetilde{\mH}_{\bold{k}}^{\mathscr{G}, \rm E_c}=
			\left[ 
			\begin{array}{c|c}
				\hspace{2mm} \Pi_{\bold{k}_0, \rm E_c}\left(\widetilde{\mH}_{\bold{k}}^{\mathscr{G}, \rm E_c}\right) \Pi_{\bold{k}_0, \rm E_c} \hspace{2mm}& \Pi_{\bold{k}_0, \rm E_c}\left(\widetilde{\mH}_{\bold{k}}^{\mathscr{G}, \rm E_c}\right) \widetilde{\Pi}_{\Omega_j, \rm E_c} \\[1em]
				\hline
				\widetilde{\Pi}_{\Omega_j, \rm E_c} 	\left(\widetilde{\mH}_{\bold{k}}^{\mathscr{G}, \rm E_c}\right) \Pi_{\bold{k}_0, \rm E_c}& \hspace{2mm} \widetilde{\Pi}_{\Omega_j, \rm E_c} \left(\widetilde{\mH}_{\bold{k}}^{\mathscr{G}, \rm E_c}\right)\widetilde{\Pi}_{\Omega_j, \rm E_c} \hspace{2mm}
			\end{array}\right]= \arraycolsep=6pt\def\arraystretch{2} \left[
			\begin{array}{c|c}
				\mathcal{A}_{\bold{k}} & \mathcal{B}_{\bold{k}} \\
				\hline
				\mathcal{B}^*_{\bold{k}} & \mathcal{C}_{\bold{k}}
			\end{array}
			\right].
		\end{equation}
		We will now consider the convergence properties of each of the sub-matrices appearing in Equation \eqref{eq:prop_2_3}. 
		
		\vspace{2.5mm}
	\noindent	\textbf{Convergence properties of sub-matrix $\mathcal{A}_{\bold{k}}$.}~
		
		For the sub-matrix $\mathcal{A}_{\bold{k}}$, using Equation \eqref{eq:fiber_mod} we have that for all $\bold{k}\in \Omega_j$ and any $\Phi \in \mX_{\bold{k}_0}^{\rm E_c}$ it holds that
		\begin{align}\label{eq:prop_2_3b}
			\mathcal{A}_{\bold{k}}\Phi 	=\Pi_{\bold{k}_0, \rm E_c}\left(\widetilde{\mH}_{\bold{k}}^{\mathscr{G}, \rm E_c}\right)  \Phi = \hspace{2mm}\Pi_{\bold{k}_0, \rm E_c} \hspace{-4 mm}\sum_{\substack{\bold{G} \in \mathbb{L}^*\\[0.2em] \vert \bold{k}_0 +\bold{G}\vert^2 < 2\rm E_c} }\widehat{\Phi}_{\bold{G}} \left({\rm E_c} \mathscr{G} \left(\frac{\vert \bold{G}+ \bold{k}\vert }{\sqrt{2 \rm E_c}}\right) +V\right) e_{\bold{G}}.
		\end{align}
		Since the blow-up function $\mathscr{G}$ is continuous on the interval $(0, 1)$, we see immediately that for any sequence $\left\{\bold{k}_\ell\right\}_{\ell \in \mathbb{N}^*} \subset \Omega_j$ such that $\lim_{\ell \to \infty}\bold{k}_\ell =\bold{k}_0$ it holds that
		\begin{align}
			\lim_{\ell \to\infty}	\left\Vert 	\mathcal{A}_{\bold{k}_\ell}\Phi -   \widetilde{\mH}_{\bold{k}_0}^{\mathscr{G}, \rm E_c}\Phi \right\Vert_{L^2_{\rm per} (\Omega)} =0 \qquad \text{and thus} \qquad \lim_{\ell \to \infty} \left\Vert \mathcal{A}_{\bold{k}_\ell} -\widetilde{\mH}_{\bold{k}_0}^{\mathscr{G}, \rm E_c}\right\Vert_2=0 \label{eq:prop_2_4}
		\end{align}
		Additionally, if the blow-up function $\mathscr{G}$ is of class $\mathscr{C}^m$ on $(0, 1)$ for $m>0$, then Equation \eqref{eq:prop_2_3b} also allows us to deduce that the sub-matrix $\mathcal{A}_{\bold{k}}, ~\bold{k} \in \mathbb{R}^d$ is continuously differentiable up to order $m$ at~$\bold{k}=\bold{k}_0$.
		
		\vspace{2.5mm}
		
	\noindent \textbf{Convergence properties of sub-matrices $\mathcal{B}_{\bold{k}}$ and $\mathcal{B}_{\bold{k}}^*$.}~
		
		Notice that the blow-up function $\mathscr{G}$ does not appear in the off-diagonal blocks $\mathcal{B}$ and $\mathcal{B}^*$ and that the effective potential $V$ is independent of $\bold{k}$ and $\bold{k}_0$.  Consequently, a similar argument as the one used for the sub-matrix $\mathcal{A}_{\bold{k}}$ yields that for any sequence $\left\{\bold{k}_\ell\right\}_{\ell \in \mathbb{N}^*} \subset \Omega_j$ such that $\lim_{\ell \to \infty}\bold{k}_\ell =\bold{k}_0$ we have
		\begin{align}  \label{eq:prop_2_5}
			\lim_{\ell \to \infty} \left \Vert \mathcal{B}_{\bold{k}_\ell} - \mathcal{B}_{\bold{k}_0}\right\Vert_2=0 \qquad \text{and} \qquad \lim_{\ell \to \infty} \left \Vert\mathcal{B}^*_{\bold{k}_\ell} -  \mathcal{B}^*_{\bold{k}_0}\right\Vert_2=0.
		\end{align}
		where $\mathcal{B}_{\bold{k}_0}, \mathcal{B}_{\bold{k}_0}^*$ are fixed, rectangular matrices that are independent of the specific choice of sequence $\left\{\bold{k}_\ell\right\}$, although they depend of course on the chosen sector $\Omega_j$. 
		
		\vspace{3mm}
	\noindent \textbf{Convergence properties of sub-matrix $\mathcal{C}_{\bold{k}}$.}~
		
		We use once again Equation \eqref{eq:fiber_mod} to deduce that for all  $\bold{k} \in \Omega_j$ and any $\Phi \in \widetilde{\mX}_{\bold{k}_0, \Omega_{j}}^{\rm E_c}$ it holds that
		\begin{align}\label{eq:prop_2_6}
			\mathcal{C}_{\bold{k}}\Phi= 	\widetilde{\Pi}_{\Omega_j, \rm E_c} \left(\widetilde{\mH}_{\bold{k}}^{\mathscr{G}, \rm E_c}\right) \Phi &=\widetilde{\Pi}_{\Omega_j, \rm E_c} \hspace{-2 mm}\sum_{{\bold{G} \in \{\bold{G}_{j_1}, \ldots,\bold{G}_{j_{\ell}} \} \subset \mathbb{L}^*}} \hspace{-2 mm}\widehat{\Phi}_{\bold{G}} \left({\rm E_c} \mathscr{G} \left(\frac{\vert \bold{G}+\bold{k}\vert }{\sqrt{2 \rm E_c}}\right) +V\right) e_{\bold{G}}.
		\end{align}
		
		Consider now a sequence $\left\{\bold{k}_{\ell}\right\}_{j \in \mathbb{N}^*} \subset \Omega_j$ such that $\lim_{\ell \to \infty}\bold{k}_\ell=\bold{k}_0$.   Recalling that $\{\bold{G}_{j_1}, \ldots, \bold{G}_{j_{\ell}}\} \subset \mS_{\bold{k}_0}^{\rm E_c}$ and the corresponding Fourier modes $e_{\bold{G}_{j_1}}, \ldots, e_{\bold{G}_{j_\ell}}$ are elements of $\widetilde{\mathcal{B}_{\Omega_j}}^{\rm E_c}$, and using Equation \eqref{eq:Cone} we see that
		\begin{align*}
			\vert \bold{G}_{j_k}+ \bold{k}_\ell\vert &\leq \sqrt{2\rm E_c} \quad \text{ for all }\ell\in \mathbb{N}^*, ~j_1, \ldots, j_{\ell} \quad \text{and}\\
			\lim_{\ell \to \infty} \vert \bold{G}_{j_k}+ \bold{k}_\ell\vert  &= \sqrt{2\rm E_c} \quad \text{ for all }  ~j_1, \ldots, j_{\ell}.
		\end{align*}
		Since, on the one hand the blow-up function $\mathscr{G}(x)$ has a singularity at $x=1$, and on the other hand the effective potential $V$ is independent of $\{\bold{k}_\ell\}_{\ell \in \mathbb{N}^*} \subset \Omega_j $ we infer that for all $\ell \in \mathbb{N}^*$, we can write the matrix $\mathcal{C}_{\bold{k}_\ell}$ in the form
		\begin{align}\label{eq:prop_2_60}
			\mathcal{C}_{\bold{k}_\ell } = \mathcal{D}_{\bold{k}_\ell} + \mathcal{N}_{\bold{k}_0},
		\end{align}
		where $\mathcal{D}_{\bold{k}_\ell}$ and $\mathcal{N}_{\bold{k}_0}$ are both square matrices of dimension $\text{dim}\widetilde{\mX}_{\bold{k}_0, \Omega_{j}}^{\rm E_c}$, and the matrix $\mathcal{D}_{\bold{k}_\ell }$ is diagonal with entries that all diverge to $+\infty$ in the limit $\ell \to \infty$ while the entries of $\mathcal{N}_{\bold{k}_0}$ are independent of $\ell$.  A particular consequence of this is that the matrix $\mathcal{D}_{\bold{k}_{\ell}}$ is invertible for $\ell$ sufficiently large.
		
		We now claim that for all $p \in \{0, \ldots, m\}$ it holds that
		\begin{align*}
			\lim_{\ell \to \infty}\frac{	\left \Vert \mathcal{D}_{\bold{k}_\ell }^{-1} \right\Vert_2}{\vert \bold{k}_0 - \bold{k}_{\ell}\vert^p} =0.
		\end{align*}
		To see this, we recall Equation \eqref{eq:prop_2_6} and the fact that $\mathcal{D}_{\bold{k}_\ell}$ is diagonal so that it suffices to show that $\forall\bold{G}_{j_k} \in \{\bold{G}_{j_1}, \ldots, \bold{G}_{j_{\ell}}\}$ and any $p \in \{0, \ldots, m\}$ it holds that
		\begin{align*}
			\lim_{\ell \to \infty} \; \frac{1}{\vert \bold{k}_0 - \bold{k}_{\ell}\vert^p} \cdot \frac{1}{ \mathscr{G} \left(\frac{\vert \bold{G}_{j_k}+ \bold{k}_\ell\vert }{\sqrt{2 \rm E_c}}\right)}  &=0 \quad \text{or equivalently}\\[0.5em]
			\lim_{\ell \to \infty}  \; \vert \bold{k}_0 - \bold{k}_{\ell}\vert^p \cdot { \mathscr{G} \left(\frac{\vert \bold{G}_{j_k}+ \bold{k}_\ell\vert }{\sqrt{2 \rm E_c}}\right)}  &=+\infty.
		\end{align*}
		Using simple algebra, one can show that for all $\bold{G}_{j_k} \in \{\bold{G}_{j_1}, \ldots, \bold{G}_{j_{\ell}}\}$ and any $p \in \{0, \ldots, m\}$ we have
		\begin{align*}
			\lim_{\ell \to \infty} \; \vert \bold{k}_0 - \bold{k}_{\ell}\vert^p \cdot { \mathscr{G} \left(\frac{\vert \bold{G}_{j_k}+ \bold{k}_\ell\vert }{\sqrt{2 \rm E_c}}\right)}  =+\infty\quad \iff \quad \lim_{x \to 1^{-}}  \; (1-x)^j \cdot { \mathscr{G} \left(x\right)}= +\infty.
		\end{align*}
		But this latter condition is satisfied by the blow-up function $\mathscr{G}$ by assumption (see Definition \ref{def:g}). We therefore conclude that for all $p \in \{0, \ldots, m\}$ it holds that
		\begin{align}
			\lim_{\ell \to \infty}	\frac{\mathcal{D}_{\bold{k}_\ell}^{-1}}{\vert \bold{k}_0 - \bold{k}_{\ell}\vert^p } =0 =	\lim_{\ell \to \infty}	\frac{\mathcal{C}_{\bold{k}_\ell}^{-1}}{\vert \bold{k}_0 - \bold{k}_{\ell}\vert^p } \qquad \text{in the matrix 2-norm topology.} \label{eq:prop_2_7}
		\end{align}
		
		\vspace{4mm}
		
		Consider again a sequence $\left\{\bold{k}_{\ell}\right\}_{\ell \in \mathbb{N}^*} \subset \Omega_j$ such that $\lim_{\ell \to \infty}\bold{k}_\ell=\bold{k}_0$. Having understood the convergence properties of the sub-blocks of the modified Hamiltonian matrix $\widetilde{\mH}_{\bold{k}}^{\mathscr{G}, \rm E_c}$ , we will now study the convergence of the approximate energy band $\widetilde{\varepsilon}_{n}^{\rm E_c}$ and its derivatives up to order $m$ as functions of the sequence $\left\{\bold{k}_{\ell}\right\}_{\ell \in \mathbb{N}^*}$. 
		
		
		\vspace{4mm}
		
	\noindent\textbf{Continuity of energy bands.}~
		
		Recall that $M_{\rm E_c}(\bold{k})$ denotes the dimension of the matrix $\widetilde{\mH}_{\bold{k}}^{\mathscr{G}, \rm E_c}$ at $\bold{k}\in \mathbb{R}^d$. Thanks to the definition of the set $\Omega_j$, we see that for each element of the sequence $ \{\bold{k}_{\ell}\}_{\ell \in \mathbb{N}^*}$, the dimension of $\widetilde{\mH}_{\bold{k}_{\ell}}^{\mathscr{G}, \rm E_c}$ remains constant, i.e.,$M_{\rm E_c}(\bold{k}_{\ell})=M \in \mathbb{N}^*$. Consequently, we may apply Lemma \ref{lem:contin} to the modified Hamiltonian matrices $\left\{\widetilde{\mH}_{\bold{k}_{\ell}}^{\mathscr{G}, \rm E_c}\right\}_{\ell \in \mathbb{N}^*}$. Indeed, thanks to the convergence properties of the sub-matrices $\{\mathcal{A}_{\bold{k}_{\ell}}\}_{\ell \in \mathbb{N}^*}, \{\mathcal{B}_{\bold{k}_{\ell}}\}_{\ell \in \mathbb{N}^*}, \{\mathcal{B}^*_{\bold{k}_{\ell}}\}_{\ell \in \mathbb{N}^*}$ and $\{\mathcal{C}_{\bold{k}_{\ell}}\}_{\ell \in \mathbb{N}^*}$ established above (and taking a subsequence, if necessary, to ensure the invertibility of all $\mathcal{C}_{\bold{k}_{\ell}}$), we see that the assumptions of Lemma \ref{lem:contin} are satisfied. Denoting therefore, $p=\text{\rm dim } \mathcal{A}_{\bold{k}_0}$ and recalling that $\widetilde{\mH}_{\bold{k}_{0}}^{\mathscr{G}, \rm E_c}=\mathcal{A}_{\bold{k}_0}$, we deduce that for each $q \in \{1, \ldots, p\}$ it holds that
			\begin{align}\label{eq:new_hassan_0}
				\lim_{\ell \to \infty} \widetilde{\varepsilon}_q^{\rm E_c}(\bold{k}_{\ell})   &= \widetilde{\varepsilon}_q^{\rm E_c}(\bold{k}_{0})\quad \text{and} \\ 	
				\lim_{\ell \to \infty} \widetilde{\varepsilon}^{\rm E_c}_{p+1} (\bold{k}_{\ell})&= \lim_{\ell \to \infty} \widetilde{\varepsilon}^{\rm E_c}_{p+2} (\bold{k}_{\ell})= \cdots =\lim_{\ell \to \infty} \widetilde{\varepsilon}^{\rm E_c}_{M} (\bold{k}_{\ell})= \infty. \nonumber
			\end{align}

			In order to conclude the continuity of the bounded energy bands $\{\widetilde{\varepsilon}^{\rm E_c}_{q}\}_{q \in \{1, \ldots, p\}}$, it suffices to recall that we have considered a sequence $\{\bold{k}_{\ell}\}_{\ell \in \mathbb{N}^*} \subset \Omega_j$ for some $j \in \{2, \ldots, 2^M\}$. But since $\Omega_j$ was chosen arbitrarily and there are only a finite number of possible choices for $\Omega_j$, we conclude that Equation \eqref{eq:new_hassan_0} holds for \underline{any} sequence $\{\bold{k}_{\ell}\}_{\ell \in \mathbb{N}^*} \subset \mathbb{B}_{\delta}(\bold{k})$. It follows that all bounded, modified energy bands $\{\widetilde{\varepsilon}^{\rm E_c}_{q}\}_{q=1}^p$ are continuous at $\bold{k}=\bold{k}_0$ as claimed. Noting that $p=\text{\rm dim } \mathcal{A}_{\bold{k}_0} \leq M_{\rm E_c}^-$ completes the proof of continuity.

			If the blow-up function $\mathscr{G}$ satisfies Properties (1)-(4) from Definition \ref{def:g} only for $m=0$, then we are done. Hence, we may assume that $m \geq 1$ and that all eigenvalues of $\mathcal{A}_{\bold{k}_0}=\widetilde{\mH}^{\mathscr{G}, \rm E_c}_{\bold{k}_0}$ are simple. We study next the regularity of the derivatives of the bounded, modified energy bands.

		\vspace{4mm}
		
		\textbf{First order differentiability of energy bands.}~
		
		We begin by studying the convergence of the eigenvectors associated with the bounded energy bands $\widetilde{\varepsilon}_q^{\rm E_c}(\bold{k}_{\ell}), ~q \in \{1, \ldots, p\}$. The primary tool we will use for this study will be the Schur complement associated with the block decomposition~\eqref{eq:prop_2_3} of the modified Hamiltonian matrix $\widetilde{\mH}_{\bold{k}_\ell}^{\mathscr{G}, \rm E_c}$.
		
		Let $q \in \{1, \ldots, p\}$ be the index of a bounded energy band. A straightforward calculation using the block decomposition \eqref{eq:prop_2_3} shows that for any $\bold{k}_{\ell}\in \Omega_j$ it holds that
		\begin{equation}\label{eq:prop_2_8}
			\Pi_{\bold{k}_0, \rm E_c}	\left(\widetilde{\varepsilon}_q^{\rm E_c} (\bold{k}_\ell)\;\widetilde{u}^{{\rm E_c}}_{q, \bold{k}_\ell }\right)= \mathcal{A}_{\bold{k}_\ell}  \widetilde{u}^{{\rm E_c}}_{q, \bold{k}_\ell }- \mathcal{B}_{\bold{k}_\ell} \left(\mathcal{C}_{\bold{k}_\ell}  - \widetilde{\varepsilon}_q^{\rm E_c} (\bold{k}_\ell)\right)^{-1}\mathcal{B}_{\bold{k}_\ell}^*\widetilde{u}^{{\rm E_c}}_{q, \bold{k}_\ell },
		\end{equation}
		where $\widetilde{u}^{{\rm E_c}}_{q, \bold{k}_\ell}$ denotes the $q^{\rm th}$ normalized eigenfunction of the modified Hamiltonian matrix~$\widetilde{\mH}^{\mathscr{G}, \rm E_c}_{\bold{k}_\ell}$. Additionally, thanks to Equations \eqref{eq:prop_2_5} and~\eqref{eq:prop_2_7}, we deduce from Equation \eqref{eq:prop_2_8} that 
		\begin{align}\label{eq:prop_2_9}
			\lim_{\ell \to \infty} \left(\Pi_{\bold{k}_0, \rm E_c}	\left(\widetilde{\varepsilon}_q^{\rm E_c} (\bold{k}_\ell)\widetilde{u}^{{\rm E_c}}_{q, \bold{k}_\ell}\right)- \mathcal{A}_{\bold{k}_\ell}  \widetilde{u}^{{\rm E_c}}_{q, \bold{k}_\ell}\right)=0.
		\end{align}

		Next, observe that since the sequence $\left\{\Pi_{\bold{k}_0, \rm E_c}\; \widetilde{u}^{{\rm E_c}}_{q, \bold{k}_{\ell}} \right\}_{{\ell \in \mathbb{N}^*}}$ is bounded, it possesses a convergent subsequence, which we also write as $\left\{\Pi_{\bold{k}_0, \rm E_c}\; \widetilde{u}^{{\rm E_c}}_{q, \bold{k}_{\ell}}\right\}_{{\ell \in \mathbb{N}^*}}$. We can then deduce from Equations \eqref{eq:prop_2_4}, \eqref{eq:new_hassan_0}, and \eqref{eq:prop_2_9} that
		\begin{align*}
			\widetilde{\varepsilon}_q^{\rm E_c} (\bold{k}_0)\; \left(\lim_{\ell \to \infty}  \Pi_{\bold{k}_0, \rm E_c}	\widetilde{u}^{{\rm E_c}}_{q, \bold{k}_{\ell}}\right)&=   \mathcal{A}_{\bold{k}_{0}} \; \left(\lim_{\ell \to \infty}   \Pi_{\bold{k}_0, \rm E_c}	\widetilde{u}^{{\rm E_c}}_{q, \bold{k}_{\ell}}\right).
		\end{align*}
		But this implies that either $\underset{\ell \to \infty}{\lim}  \Pi_{\bold{k}_0, \rm E_c}	\widetilde{u}^{{\rm E_c}}_{q, \bold{k}_{\ell}}$ is (up to normalization) the eigenvector $\widetilde{u}^{{\rm E_c}}_{q, \bold{k}_0}$ associated with the eigenvalue $\widetilde{\varepsilon}_q^{\rm E_c} (\bold{k}_0)$ or $\underset{\ell \to \infty}{\lim} \Pi_{\bold{k}_0, \rm E_c}	\widetilde{u}^{{\rm E_c}}_{q, \bold{k}_{\ell}}=0$. Suppose on the contrary that $\underset{\ell \to \infty}{\lim}  \Pi_{\bold{k}_0, \rm E_c}\widetilde{u}^{{\rm E_c}}_{q, \bold{k}_{\ell}}=0$ and note that the block decomposition \eqref{eq:prop_2_3} implies that for all $\bold{k}_{\ell}\in \Omega_j$ we have
		\begin{equation}\label{eq:prop_2_8b}
			\Pi^{\perp}_{\bold{k}_0, \rm E_c}	\left(\widetilde{\varepsilon}_q^{\rm E_c} (\bold{k}_\ell)\;\widetilde{u}^{{\rm E_c}}_{q, \bold{k}_\ell }\right)= \mathcal{B}^*_{\bold{k}_\ell} \widetilde{u}^{{\rm E_c}}_{q, \bold{k}_\ell }+ \mathcal{C}_{\bold{k}_\ell}\widetilde{u}^{{\rm E_c}}_{q, \bold{k}_\ell }.
		\end{equation}
		We can now take the limit $\ell \to \infty$ on both sides of Equation \eqref{eq:prop_2_8b}. But $\underset{\ell \to \infty}{\lim} \widetilde{\varepsilon}_q^{\rm E_c} (\bold{k}_\ell)=\widetilde{\varepsilon}_q^{\rm E_c} (\bold{k}_0)< \infty$ and $ \underset{\ell \to \infty}{\lim} \Vert\Pi^{\perp}_{\bold{k}_0, \rm E_c}\; \widetilde{u}^{{\rm E_c}}_{q, \bold{k}_{\ell}} \Vert=1$ while all eigenvalues of the matrix $\mathcal{C}_{\bold{k}_{\ell}}$ diverge to $+\infty$ in the limit $\ell \to \infty$. Consequently, we must have 
		\begin{align}\label{eq:hassan_01}
			\underset{\ell \to \infty}{\lim}\Pi_{\bold{k}_0, \rm E_c}\; \widetilde{u}^{{\rm E_c}}_{q, \bold{k}_{\ell}}= \underset{\ell \to \infty}{\lim} \widetilde{u}^{{\rm E_c}}_{q, \bold{k}_{\ell}}= \widetilde{u}^{{\rm E_c}}_{q, \bold{k}_0}.
		\end{align}
		
		Moreover, similar to the argument for sequential continuity of the approximate, modified energy bands, we conclude from the fact that $\{\bold{k}_{\ell}\}_{\ell \in \mathbb{N}^*} \subset \Omega_j$ and only a finite number of possibilities exist for the choice of $j \in \{2, \ldots, 2^M\}$, that Equation~\eqref{eq:hassan_01} also holds for  \underline{any} sequence $\{\bold{k}_{\ell}\}_{\ell \in \mathbb{N}^*} \subset \mathbb{B}_{\delta}(\bold{k})$, which proves sequential continuity of the normalized eigenfunctions associated with all bounded energy bands.

		Equipped with the convergence properties of the eigenvectors associated with the bounded energy bands, we can now consider the first order derivatives of the bounded, modified energy band  $\widetilde{\varepsilon}_q^{\rm E_c}(
		\bold{k})$, $q \in \{1, \ldots, p\}$ at $\bold{k}=\bold{k}_0$. Thanks once again to Equations~\eqref{eq:prop_2_5} and~\eqref{eq:prop_2_7}, we deduce from the Schur-type decomposition \eqref{eq:prop_2_8} that 
		\begin{align*}
			\lim_{\ell \to \infty} \frac{\left(\Pi_{\bold{k}_0, \rm E_c}	\left(\widetilde{\varepsilon}_q^{\rm E_c} (\bold{k}_\ell)\;\widetilde{u}^{{\rm E_c}}_{q, \bold{k}_\ell}\right)- \mathcal{A}_{\bold{k}_\ell}  \; \widetilde{u}^{{\rm E_c}}_{q, \bold{k}_\ell}\right)}{ \vert \bold{k}_{\ell} - \bold{k}_{0}\vert}=0.
		\end{align*}
		Adding and subtracting the term $\widetilde{\varepsilon}_q^{\rm E_c} (\bold{k}_0) \; \Pi_{\bold{k}_0, \rm E_c}\widetilde{u}^{{\rm E_c}}_{q, \bold{k}_{\ell}}$, using the fact that $\widetilde{\varepsilon}_q^{\rm E_c} (\bold{k}_0)\widetilde{u}^{{\rm E_c}}_{q, \bold{k}_0}= \mathcal{A}_{\bold{k}_0}\widetilde{u}^{{\rm E_c}}_{q, \bold{k}_0}$, and taking the inner product with the eigenfunction~$\widetilde{u}^{{\rm E_c}}_{q, \bold{k}_0}$ then yields 
		\begin{align}\label{eq:prop_2_10}
			\lim_{\ell \to \infty} \frac{\left(\left(\widetilde{\varepsilon}_q^{\rm E_c} (\bold{k}_\ell)-\widetilde{\varepsilon}_q^{\rm E_c} (\bold{k}_0)+ \mathcal{A}_{\bold{k}_0}-\mathcal{A}_{\bold{k}_\ell}  \right)\widetilde{u}^{{\rm E_c}}_{q, \bold{k}_0}, \Pi_{\bold{k}_0, \rm E_c}\;\widetilde{u}^{{\rm E_c}}_{q, \bold{k}_\ell}\right)_{L^2_{\rm per}(\Omega)}}{\vert \bold{k}_{\ell}  - \bold{k}_0\vert}= 0.
		\end{align}
		
		Next, recall that  $\mathcal{A}_{\bold{k}}$ is $m$-times continuously differentiable at $\bold{k}=\bold{k}_0$ and 
		denote by ${\rm d}\mathcal{A}_{\bold{k}_0}\colon \mathbb{R}^d \rightarrow \mathbb{R}^{p \times p}$ the total derivative of $\mathcal{A}_{\bold{k}}$ at $\bold{k}=\bold{k}_0$. Adding and subtracting the term $ {\rm d}\mathcal{A}_{\bold{k}_0}[\bold{k}_{\ell}-\bold{k}_0]$, i.e., ${\rm d}\mathcal{A}_{\bold{k}_0}$ acting on the vector $\bold{k}_{\ell}-\bold{k}_0$, we can deduce from Equation \eqref{eq:prop_2_10} that
		\begin{align}\label{eq:prop_2_10b}
			&\lim_{\ell \to \infty} \frac{\left(\left(\widetilde{\varepsilon}_q^{\rm E_c} (\bold{k}_\ell)-\widetilde{\varepsilon}_q^{\rm E_c} (\bold{k}_0)- {\rm d}\mathcal{A}_{\bold{k}_0}[\bold{k}_{\ell}-\bold{k}_0] \right)\widetilde{u}^{{\rm E_c}}_{q, \bold{k}_0}, \Pi_{\bold{k}_0, \rm E_c}\;\widetilde{u}^{{\rm E_c}}_{q, \bold{k}_\ell}\right)_{L^2_{\rm per}(\Omega)}}{\vert \bold{k}_{\ell}  - \bold{k}_0\vert}=0.
		\end{align}
		Adding and subtracting the term $\left({\rm d}\mathcal{A}_{\bold{k}_0}[\bold{k}_{\ell}-\bold{k}_0]\widetilde{u}^{{\rm E_c}}_{q, \bold{k}_0}, \widetilde{u}^{{\rm E_c}}_{q, \bold{k}_0}\right)_{L^2_{\rm per}(\Omega)}$ and using simple algebra, it can be deduced~that
		\begin{align}\label{eq:prop_2_12}
			&\lim_{\ell \to \infty} \frac{\left(\widetilde{\varepsilon}_q^{\rm E_c} (\bold{k}_\ell)-\widetilde{\varepsilon}_q^{\rm E_c} (\bold{k}_0)\right) -\left({\rm d}\mathcal{A}_{\bold{k}_0}[\bold{k}_{\ell}-\bold{k}_0]\widetilde{u}^{{\rm E_c}}_{q, \bold{k}_0}, \widetilde{u}^{{\rm E_c}}_{q, \bold{k}_0}\right)_{L^2_{\rm per}(\Omega)}}{\vert \bold{k}_{\ell} - \bold{k}_{0}\vert}=0.
		\end{align}

		Similar to the argument for sequential continuity of the approximate, modified energy bands, we conclude from the fact that $\{\bold{k}_{\ell}\}_{\ell \in \mathbb{N}^*} \subset \Omega_j$ and only a finite number of possibilities exist for the choice of $j \in \{2, \ldots, 2^M\}$, that Equation~\eqref{eq:prop_2_12} also holds for  \underline{any} sequence $\{\bold{k}_{\ell}\}_{\ell \in \mathbb{N}^*} \subset \mathbb{B}_{\delta}(\bold{k})$. Thus, the total derivative of the approximate, modified energy bands exists at~$\bold{k}_0$ and is given by
		\begin{align}\label{eq:band_der}
			{\rm d}\widetilde{\varepsilon}_{q}^{\rm E_c}( \bold{k}_0)= \left({\rm d}\mathcal{A}_{\bold{k}_0}\widetilde{u}^{{\rm E_c}}_{q, \bold{k}_0}, \widetilde{u}^{{\rm E_c}}_{q, \bold{k}_0}\right)_{L^2_{\rm per}(\Omega)}.
		\end{align}
		As a consequence, $\widetilde{\varepsilon}_q^{\rm E_c}$ is of class $\mathscr{C}^1$ at $\bold{k}_0$ as claimed. Noting that we picked an arbitrary $q \in \{1\, \ldots, p\}$ where $p=\text{\rm dim } \mathcal{A}_{\bold{k}_0} \leq M_{\rm E_c}^-$ completes the proof of differentiability of order one.

		If the blow-up function $\mathscr{G}$ satisfies Properties (1)-(4) from Definition \ref{def:g} only for $m\leq 1$, then we are done. Hence, we may assume that $m \geq 2$. 
		
		\vspace{4mm}

	\noindent \textbf{Higher order differentiability of energy bands.}~
		
		Imitating the procedure carried out for the case $m=1$, we will first make use of order one differentiability of the energy band $\widetilde{\varepsilon}_{q, \bold{k}_0}^{\rm E_c}, ~ q\in \{1, \ldots, p\}$, to establish order one differentiability of the corresponding eigenfunction~$\widetilde{u}^{\rm E_c}_{q, \bold{k}_0}$.
		
		Let $q \in \{1, \ldots, p\}$ be the index of a bounded energy band. As a first remark, let us recall that by construction, the approximation space $\mX_{\bold{k}}^{\rm E_c}$ is identical for all $\bold{k} \in \Omega_j$. Since the sequence $\{\bold{k}_{\ell}\}_{\ell \in \mathbb{N}^*} \subset \Omega_j$, it can readily be deduced from the definition of the modified Hamiltonian matrix $\widetilde{\mH}_{\bold{k}}^{\mathscr{G}, \rm E_c}$ given by Equation \eqref{eq:fiber_mod} and the regularity properties of the blow-up function $\mathscr{G}$ that the modified Hamiltonian matrix $\widetilde{\mH}_{\bold{k}}^{\mathscr{G}, \rm E_c}$ is $m$-times continuously differentiable at any $\bold{k}=\bold{k}_{\ell}$. Moreover, we have assumed that all eigenvalues of $\mathcal{A}_{\bold{k}_0}=\widetilde{\mH}^{\mathscr{G}, \rm E_c}_{\bold{k}_0}$ are simple. Therefore, as discussed in \textbf{Case one} of the current proof, the implicit function theorem can be used to prove that for $\ell$ sufficiently large, the energy band $\widetilde{\varepsilon}_{q}^{\rm E_c}(\bold{k})$ and the associated (normalized) eigenfunction $\widetilde{u}^{\rm E_c}_{q, \bold{k}}$ are $m$-times continuously differentiable (as a function of $\bold{k} \in \mathbb{R}^d$) at any $\bold{k}=\bold{k}_{\ell}$. Without loss of generality, we may assume that this is the case for all $\ell \in \mathbb{N}^*$.
		
		Next, let us recall the Schur-type decomposition \eqref{eq:prop_2_8} which offers an expression for the eigenvalue $\widetilde{\varepsilon}_{q}^{\rm E_c}(\bold{k})$ in terms of the block decomposition and Schur complement of the modified Hamiltonian matrix $\widetilde{\mH}_{\bold{k}}^{\mathscr{G}, \rm E_c}$. Taking partial derivatives with respect to the $i^{\rm th}$ component of  $\bold{k}= (\bold{k}_1, \ldots, \bold{k}_d) \in \mathbb{R}^d$ of this equation yields that for any $\ell \in \mathbb{N}^*$ it holds that
		\begin{align}\nonumber
			&\left(\mathcal{A}_{\bold{k}_{\ell}} - \widetilde{\varepsilon}_q^{\rm E_c} (\bold{k}_\ell) - \mathcal{B}_{\bold{k}_\ell} \left(\mathcal{C}_{\bold{k}_\ell}  - \widetilde{\varepsilon}_q^{\rm E_c} (\bold{k}_\ell)\right)^{-1}\mathcal{B}_{\bold{k}_\ell}^*\right) \; \partial_{i}\Pi_{\bold{k}_0, \rm E_c} \widetilde{u}^{{\rm E_c}}_{q, \bold{k}_\ell } \\=  -&\Bigg(\partial_{i}\mathcal{A}_{\bold{k}_{\ell}} - \partial_{i}\widetilde{\varepsilon}_q^{\rm E_c} (\bold{k}_\ell) \label{eq:prop_2_13}  \\
			+& \mathcal{B}_{\bold{k}_\ell} \left(\mathcal{C}_{\bold{k}_\ell}  - \widetilde{\varepsilon}_q^{\rm E_c} (\bold{k}_\ell)\right)^{-1} \left(\partial_{i} \mathcal{C}_{\bold{k}_\ell} -\partial_{i}\widetilde{\varepsilon}_q^{\rm E_c} (\bold{k}_\ell) \right)\left(\mathcal{C}_{\bold{k}_\ell}  - \widetilde{\varepsilon}_q^{\rm E_c} (\bold{k}_\ell)\right)^{-1}\mathcal{B}_{\bold{k}_\ell}^*\Bigg) \; \Pi_{\bold{k}_0, \rm E_c} \widetilde{u}^{{\rm E_c}}_{q, \bold{k}_\ell }, \nonumber
		\end{align}
		where we have used the fact that the sub-matrices $\mathcal{B}_{\bold{k}}$ and $\mathcal{B}_{\bold{k}}^*$ do not change for different choices of $\bold{k} \in \Omega_j$ while the sub-matrices $\mathcal{A}_{\bold{k}}$ and $\mathcal{C}_{\bold{k}}$ are $m$-times continuously differentiable by construction for all $\bold{k}=\bold{k}_{\ell} \in \Omega_j$. 
		
		Let now $\left(\lambda_{\mathcal{A}_{\bold{k}_0}}, v_{\bold{k}_0}\right) = \left( \widetilde{\varepsilon}^{\rm E_c}_{q}(\bold{k}_{0}), \widetilde{u}^{{\rm E_c}}_{q, \bold{k}_{0}}\right)$ denote the unique, normalized eigenpair of the matrix $\mathcal{A}_{\bold{k}_0} \in \mathbb{R}^{p \times p}$ such that $\underset{\ell \to \infty}{\lim} \widetilde{\varepsilon}^{\rm E_c}_{q}(\bold{k}_{\ell}) = \lambda_{\mathcal{A}_{\bold{k}_0}}$ and $\underset{\ell \to \infty}{\lim}\Pi_{\bold{k}_0, \rm E_c}\widetilde{u}^{{\rm E_c}}_{q, \bold{k}_{\ell}}=  v_{\bold{k}_0}$. Thanks to the regularity properties of the sub-matrix $\mathcal{A}_{\bold{k}_0}$, we can once again deduce that both $\lambda_{\mathcal{A}_{\bold{k}}}$, and $v_{\bold{k}}$ are $m$-times continuously differentiable at $\bold{k}=\bold{k}_0$. Our goal now is to use Equation \eqref{eq:prop_2_13} to demonstrate that $\underset{\ell \to \infty}{\lim} \partial_{i}\widetilde{u}^{{\rm E_c}}_{q, \bold{k}_\ell } = \partial_{i} v_{\bold{k}_0}$.
		
		To this end, we claim that in fact
		\begin{align}\label{eq:prop_2_14}
			\lim_{\ell \to \infty}   \mathcal{B}_{\bold{k}_\ell} \left(\mathcal{C}_{\bold{k}_\ell}  - \widetilde{\varepsilon}_q^{\rm E_c} (\bold{k}_\ell)\right)^{-1} \left(\partial_{i} \mathcal{C}_{\bold{k}_\ell} -\partial_{i}\widetilde{\varepsilon}_q^{\rm E_c} (\bold{k}_\ell) \right)\left(\mathcal{C}_{\bold{k}_\ell}  - \widetilde{\varepsilon}_q^{\rm E_c} (\bold{k}_\ell)\right)^{-1}\mathcal{B}_{\bold{k}_\ell}^* =0.
		\end{align}
		Indeed, Equation \eqref{eq:prop_2_14} is a consequence of the properties of the blow-up function $\mathscr{G} \colon \mathbb{R}\rightarrow \mathbb{R}$ given by Definition \ref{def:g} as can easily be verified using a similar calculation as the one used to arrive at Equation \eqref{eq:prop_2_7}.

		Taking limits on both sides of Equation \eqref{eq:prop_2_13} and using the convergence properties we have proven thus far, we obtain that
		\begin{equation}\label{eq:prop_2_17}
		\begin{split}
			\lim_{\ell \to \infty} \left(\mathcal{A}_{\bold{k}_{\ell}} - \widetilde{\varepsilon}_q^{\rm E_c} (\bold{k}_\ell) - \mathcal{B}_{\bold{k}_\ell} \left(\mathcal{C}_{\bold{k}_\ell}  - \widetilde{\varepsilon}_q^{\rm E_c} (\bold{k}_\ell)\right)^{-1}\mathcal{B}_{\bold{k}_\ell}^*\right)& \; \partial_{i} \Pi_{\bold{k}_0, \rm E_c} \widetilde{u}^{{\rm E_c}}_{q, \bold{k}_\ell }\\
			=  -&\Big(\partial_{i}\mathcal{A}_{\bold{k}_0} - \partial_{i}\widetilde{\varepsilon}_q^{\rm E_c} (\bold{k}_0)\Big) \; \widetilde{u}^{{\rm E_c}}_{q, \bold{k}_0 }. 
			\end{split}
		\end{equation}

		A direct calculation now yields that the right hand side of the above equation is $L^2_{\rm per}$-orthogonal to $\text{span}\{\widetilde{u}^{{\rm E_c}}_{q, \bold{k}_0 }\}$. Moreover, thanks to the convergence properties of the sub-matrices $\mathcal{A}_{\bold{k}_{\ell}}, \mathcal{B}_{\bold{k}_{\ell}}, \mathcal{B}_{\bold{k}_{\ell}}^*$ and $\mathcal{C}_{\bold{k}_{\ell}}$, we also deduce that for~$\ell$ sufficiently large it holds that
		\begin{equation*}
			\left(\mathcal{A}_{\bold{k}_{\ell}} - \widetilde{\varepsilon}_q^{\rm E_c} (\bold{k}_\ell) - \mathcal{B}_{\bold{k}_\ell} \left(\mathcal{C}_{\bold{k}_\ell}  - \widetilde{\varepsilon}_q^{\rm E_c} (\bold{k}_\ell)\right)^{-1}\mathcal{B}_{\bold{k}_\ell}^*\right) \hspace{2mm} \text{is invertible on}\hspace{2mm}  \left\{ \widetilde{u}^{{\rm E_c}}_{q, \bold{k}_{0}}\right\}^{\perp} \subset \mX^{\rm E_c}_{\bold{k}_0}. \end{equation*}
		Consequently, Equation \eqref{eq:prop_2_17} yields that
		\begin{align}\nonumber
			\Pi_{\widetilde{u}^{{\rm E_c}}_{q}(\bold{k}_0)}^{\perp}&\Big(\lim_{\ell \to \infty} \partial_{i}\Pi_{\bold{k}_0, \rm E_c} \widetilde{u}^{{\rm E_c}}_{q, \bold{k}_\ell }\Big)\\
			=  -&\left(\Pi_{\widetilde{u}^{{\rm E_c}}_{q}(\bold{k}_0)}^{\perp}\left(\mathcal{A}_{\bold{k}_{0}} - \widetilde{\varepsilon}_q^{\rm E_c}(\bold{k}_0)\right)\Pi_{\widetilde{u}^{{\rm E_c}}_{q}(\bold{k}_0)}^{\perp} \right)^{-1}\left(\partial_{i}\mathcal{A}_{\bold{k}_0} - \partial_{i}\widetilde{\varepsilon}_q^{\rm E_c} (\bold{k}_0)\right) \; \widetilde{u}^{{\rm E_c}}_{q, \bold{k}_0 }, \label{eq:final_hassan} 
		\end{align}
		where $\Pi_{\widetilde{u}^{{\rm E_c}}_{q}(\bold{k}_0)}^{\perp}$ is the $L^2_{\rm per}$-orthogonal projection operator onto $\{\widetilde{u}^{{\rm E_c}}_{q, \bold{k}_0 }\}^{\perp} \subset \mX^{\rm E_c}_{\bold{k}_0}$. 
		
		To proceed to the conclusion, we need to demonstrate that 
		\begin{align*}
			\left(\mathbb{I}-\Pi^{\perp}_{\widetilde{u}^{{\rm E_c}}_{q}(\bold{k}_0)}\right)\left(\lim_{\ell \to \infty} \partial_{i}\Pi_{\bold{k}_0, \rm E_c} \widetilde{u}^{{\rm E_c}}_{q, \bold{k}_\ell }\right)=0 \qquad \text{and}\qquad  \lim_{\ell \to \infty} \partial_{i}\Pi^{\perp}_{\bold{k}_0, \rm E_c} \widetilde{u}^{{\rm E_c}}_{q, \bold{k}_\ell }=0.
		\end{align*}
		We first focus on the latter limit. To this end, we recall Equation \eqref{eq:prop_2_8b}, which yields that for all $\ell \in \mathbb{N}^*$ it holds that
		\begin{align}\label{eq:prop_2_20}
			\Pi^{\perp}_{\bold{k}_0, \rm E_c}	\left(\widetilde{\varepsilon}_q^{\rm E_c} (\bold{k}_\ell)\;\widetilde{u}^{{\rm E_c}}_{q, \bold{k}_\ell }\right)= \mathcal{B}^*_{\bold{k}_\ell} \widetilde{u}^{{\rm E_c}}_{q, \bold{k}_\ell }+ \mathcal{C}_{\bold{k}_\ell}\widetilde{u}^{{\rm E_c}}_{q, \bold{k}_\ell }.
		\end{align}
		We have already demonstrated that $\underset{\ell \to \infty}{\lim}\Pi^{\perp}_{\bold{k}_0, \rm E_c}	\widetilde{u}^{{\rm E_c}}_{q, \bold{k}_\ell }=0$. Recalling therefore the decomposition $\mathcal{C}_{\bold{k}_{\ell}}=\mathcal{D}_{\bold{k}_{\ell}}+ \mathcal{N}_{\bold{k}_0}$ introduced priori to Equation Equation \eqref{eq:prop_2_7}, we see that we must have
		\begin{align}\label{eq:prop_2_21}
			\lim_{\ell \to \infty} \Vert \mathcal{D}_{\bold{k}_\ell} \Pi^{\perp}_{\bold{k}_0, \rm E_c} \widetilde{u}^{{\rm E_c}}_{q, \bold{k}_\ell }\Vert_{L^2_{\rm per}(\Omega)}=0.
		\end{align}
		Taking now partial derivatives with respect to the $i^{\rm th}$ component of $\bold{k}= (\bold{k}_1, \ldots, \bold{k}_d) \in \mathbb{R}^d$ of Equation \eqref{eq:prop_2_20}, taking limits, and simplifying terms that obviously goes to zero now, we obtain that
		\begin{align*}
			\lim_{\ell \to \infty} \partial_{i} \Pi^{\perp}_{\bold{k}_0, \rm E_c}\widetilde{u}^{{\rm E_c}}_{q, \bold{k}_\ell } = -\lim_{\ell \to \infty} \left( \mathcal{C}_{\bold{k}_\ell} - \widetilde{\varepsilon}_q^{\rm E_c} (\bold{k}_\ell) \right)^{-1}\left( \partial_i \mathcal{D}_{\bold{k}_\ell} - \partial_i \widetilde{\varepsilon}_q^{\rm E_c} (\bold{k}_\ell) \right)\Pi^{\perp}_{\bold{k}_0, \rm E_c}\widetilde{u}^{{\rm E_c}}_{q, \bold{k}_\ell }=0,
		\end{align*}
		where the last equality follows again from the properties of the blow-up function $\mathscr{G}\colon \mathbb{R}\rightarrow \mathbb{R}$ defined through Definition \ref{def:g}. We conclude that $\lim_{\ell \to \infty} \partial_{i}\Pi^{\perp}_{\bold{k}_0, \rm E_c} \widetilde{u}^{{\rm E_c}}_{q, \bold{k}_\ell }=0$ as claimed.
		
		
		It remains to prove that $\big(\mathbb{I}-\Pi^{\perp}_{\widetilde{u}^{{\rm E_c}}_{q}(\bold{k}_0)}\big)\big(\lim_{\ell \to \infty} \partial_{i}\Pi_{\bold{k}_0, \rm E_c} \widetilde{u}^{{\rm E_c}}_{q, \bold{k}_\ell }\big)=0$. To this end, we note that due to the normalization of $\widetilde{u}^{{\rm E_c}}_{q, \bold{k}_\ell }$, for all $\ell \in \mathbb{N}^*$ it holds that
		\begin{align*}
			\left( \Pi_{\bold{k}_0, \rm E_c} \widetilde{u}^{{\rm E_c}}_{q, \bold{k}_\ell }, \partial_{i} \Pi_{\bold{k}_0, \rm E_c} \widetilde{u}^{{\rm E_c}}_{q, \bold{k}_\ell }\right)_{L^2_{\rm per}(\Omega)}&=\frac{1}{2}  \partial_{i} \left\Vert \Pi_{\bold{k}_0, \rm E_c}^{\perp} \widetilde{u}^{{\rm E_c}}_{q, \bold{k}_\ell }\right\Vert^2_{L^2_{\rm per}(\Omega)}\\[1em]
			\implies \left(\mathbb{I}-\Pi^{\perp}_{\widetilde{u}^{{\rm E_c}}_{q}(\bold{k}_0)}\right)\left(\lim_{\ell \to \infty} \partial_{i}\Pi_{\bold{k}_0, \rm E_c} \widetilde{u}^{{\rm E_c}}_{q, \bold{k}_\ell }\right) &= \frac{1}{2}  \lim_{\ell \to \infty} \partial_{i} \left\Vert \Pi_{\bold{k}_0, \rm E_c}^{\perp} \widetilde{u}^{{\rm E_c}}_{q, \bold{k}_\ell }\right\Vert^2_{L^2_{\rm per}(\Omega)}=0.
		\end{align*}
		Collecting these convergence results and recalling Equation \eqref{eq:final_hassan}, we see that we have in fact shown that
		\begin{equation}\label{eq:prop_2_22}
			\lim_{\ell \to \infty} \partial_{i}\widetilde{u}^{{\rm E_c}}_{q, \bold{k}_\ell }=  -\left(\Pi_{\widetilde{u}^{{\rm E_c}}_{q}(\bold{k}_0)}^{\perp}\left(\mathcal{A}_{\bold{k}_{0}} - \widetilde{\varepsilon}_q^{\rm E_c}(\bold{k}_0)\right)\Pi_{\widetilde{u}^{{\rm E_c}}_{q}(\bold{k}_0)}^{\perp} \right)^{-1}\left(\partial_{i}\mathcal{A}_{\bold{k}_0} - \partial_{i}\widetilde{\varepsilon}_q^{\rm E_c} (\bold{k}_0)\right) \; \widetilde{u}^{{\rm E_c}}_{q, \bold{k}_0 }.
		\end{equation}
		
		Similar to the argument for sequential continuity of the approximate, modified energy bands, we conclude from the fact that $\{\bold{k}_{\ell}\}_{\ell \in \mathbb{N}^*} \subset \Omega_j$ and only a finite number of possibilities exist for the choice of $j \in \{2, \ldots, 2^M\}$, that Equation~\eqref{eq:prop_2_22} also holds for  \underline{any} sequence $\{\bold{k}_{\ell}\}_{\ell \in \mathbb{N}^*} \subset \mathbb{B}_{\delta}(\bold{k})$. It is now straightforward to conclude from Equation~\eqref{eq:prop_2_22} that $\underset{\ell \to \infty}{\lim} \partial_{i}\widetilde{u}^{{\rm E_c}}_{q, \bold{k}_\ell } = \partial_{i} v_{\bold{k}_0}$ as claimed since $\partial_{i} v_{\bold{k}_0}$ by definition also satisfies the equation
		\begin{align*}
			\left(\mathcal{A}_{\bold{k}_{0}} - \widetilde{\varepsilon}_q^{\rm E_c}(\bold{k}_0)\right)\partial_{i} v_{\bold{k}_0}=  -\left(\partial_{i}\mathcal{A}_{\bold{k}_0} - \partial_{i}\widetilde{\varepsilon}_q^{\rm E_c} (\bold{k}_0)\right) \; \widetilde{u}^{{\rm E_c}}_{q, \bold{k}_0 }.
		\end{align*}
		As a consequence, $\widetilde{u}_{q, \bold{k}}^{\rm E_c}$ is of class $\mathscr{C}^1$ at $\bold{k}=\bold{k}_0$ as claimed. Noting that we picked an arbitrary $q \in \{1\, \ldots, p\}$ where $p=\text{\rm dim } \mathcal{A}_{\bold{k}_0} \leq M_{\rm E_c}^-$ completes the proof of differentiability of order one of the bounded energy band eigenfunctions. 
		
		By making use of this first order differentiability, we can perform a similar demonstration involving limits of finite-difference approximations of second order derivatives in order to establish $\mathscr{C}^2$ regularity of the modified energy band $\widetilde{\varepsilon}_q^{\rm E_c} $ at $\bold{k}_0$. For the sake of brevity, we omit the details of these (and higher order differentiability) demonstrations.

	\end{proof}

	\section*{Acknowledgements}
	
	The authors thank Antoine Levitt and Fran\c cois Gygi for useful comments and fruitful discussions. The first two authors also thank IPAM where portions of this work were completed (March-June 2022). This project has received funding from the European Research Council (ERC) under the European Union's Horizon 2020 research and innovation program (grant agreement No. 810367).

	\bibliographystyle{siam.bst}    
	\bibliography{refs.bib}

\end{document}